\documentclass[11pt,leqno]{article}

\usepackage{amsmath}
\usepackage{amssymb}
\usepackage{dsfont}
\usepackage{mathrsfs}
\usepackage{epsfig}
\usepackage{float}
\usepackage{color}

\numberwithin{equation}{section}

\newtheorem{Theorem}{Theorem}
\newtheorem{Corollary}[Theorem]{Corollary}
\newtheorem{Lemma}[Theorem]{Lemma}
\newtheorem{Proposition}[Theorem]{Proposition}

\newtheorem{remark}[Theorem]{Remark}

\newtheorem{example}[Theorem]{Example}

\numberwithin{Theorem}{section}

\bibstyle{plain}

\newcommand{\al}{\alpha}
\newcommand{\ka}{\kappa}
\newcommand{\R}{{\mathbf R}}
\newcommand{\ds}{\displaystyle}
\newcommand{\e}{\varepsilon}

\newcommand{\lra}{\longrightarrow}
\newcommand{\ra}{\rightarrow}
\newcommand{\p}{\partial}
\newcommand{\la}{\lambda}
\newcommand{\g}{\gamma}
\newcommand{\ov}{\overline}
\newcommand{\C}{{\mathbf C}}

\newcommand{\N}{{\mathbf N}}

\newcommand{\Sp}{{\mathbb S}}
\newcommand{\1}{{\mathds 1}}

\newcommand{\si}{\sigma}

\newcommand{\de}{\delta}

\newcommand{\ph}{\varphi}

\newcommand{\Ao}{{\mathcal A}}
\newcommand{\Fo}{{\mathcal F}}
\newcommand{\Po}{{\mathcal P}}
\newcommand{\Xo}{{\mathscr X}}

\newcommand{\So}{{\mathcal S}}
\newcommand{\To}{{\mathcal T}}

\newcommand{\Ro}{{\mathcal R}}
\newcommand{\Lo}{{\mathcal L}}

\newcommand{\Oo}{{\mathcal O}}
\newcommand{\Qo}{{\mathcal Q}}
\newcommand{\Uo}{{\mathcal U}}
\newcommand{\Sr}{{\mathscr{S}}}

\newcommand{\wh}{\widehat}

\newcommand{\les}{\leqslant}
\newcommand{\ges}{\geqslant}

\newcommand{\beq}{\begin{equation}}
\newcommand{\eeq}{\end{equation}}

\renewcommand{\b}[1]{\textcolor{blue}{#1}}

\newcommand{\references}[1]{\theinstitutions
\begin{thebibliography}{#1}}

\newcommand{\Endrefs}{\end{thebibliography}}

\title{Non-homogeneous Gagliardo-Nirenberg inequalities in $ \R^N$ and application to a biharmonic non-linear Schr\"odinger equation}

\author{Antonio {\sc J. Fern\'andez,}  Louis {\sc Jeanjean,}  \\
Rainer {\sc Mandel} and Mihai {\sc Mari\c{s}}}

\date{}

\begin{document}

\maketitle

\abstract{

\noindent
We develop a new method, based on the Tomas-Stein inequality, to establish non-homogeneous Gagliardo-Nirenberg-type inequalities in $\R^N$.
 Then we use these inequalities to study standing waves minimizing the {\it energy} when the $L^2-$norm (the \textit{mass})
 is kept fixed
for a  fourth-order  Schr\"odinger equation with mixed dispersion.
We prove optimal results  on the existence of  minimizers in the {\it mass-subcritical } and {\it mass-critical } cases.
In the { \it mass-supercritical} case global minimizers do not exist. However, if the Laplacian and the bi-Laplacian in the equation have the same sign, we are able to show the existence of local minimizers.
The existence of those local minimizers is significantly more difficult than the study of global minimizers in the {\it mass-subcritical} and {\it mass-critical} cases. They are global in time solutions with small $H^2-$norm that do not scatter.
Such special solutions do not exist if the Laplacian and the bi-Laplacian have opposite sign.
If the mass does not exceed some threshold $ \mu _0  \in (0,+\infty)$,
 our results on "best" local minimizers are  optimal.
 
 \smallskip
 

\smallskip
 

}

\section{Introduction}
\label{intro}

We consider the biharmonic  non-linear Schr\"odinger equation with mixed dispersion 
$$
i \p _t \psi + \al \Delta ^2 \psi + \beta \Delta \psi + \g |\psi |^{2 \si} \psi = 0 \qquad \mbox{ in } \R \times \R^N, 
\leqno{\text{(BNLS)}}
$$
where $ \al  , \si > 0 $ and $ \beta, \g \in \R$, $ \g \neq 0$.
This equation has been introduced by Karpman and Shagalov in \cite{Ka}  and \cite{KaSh}
to take into account the role of small fourth-order dispersion terms in the propagation of intense laser beams in a bulk medium with Kerr
non-linearity; see also \cite{FiIlPa}. 
It has also been used to describe the motion of a vortex filament in an
incompressible fluid (\cite{FuMo}).
The equation received considerable attention since then. 

By simple scaling it is possible to get rid of the parameters $\al, \beta, \g$. 
 Indeed, if $ \beta \neq 0 $, taking $ {\psi} ( t, x) = a \tilde{\psi}  \left( \frac tc, \frac xb \right)$ where 
 $ a = \al ^{- \frac{1}{2\si}} | \beta /2|^{\frac{1}{\si} } | \g |^{- \frac{1}{2\si}}$, 
 $ b =  \al^{\frac 12}| \beta/2|^{- \frac 12 }$, and
   $ c = 4 \al |\beta| ^{-2}$,
we see that $ \psi $ solves the above equation if and only if $\tilde{\psi} $ solves (after dropping  "$\; \tilde{\; }\; $")
\beq
\label{BNLSg}
i \psi _t + \Delta ^2 \psi + 2 \epsilon \Delta \psi + \vartheta |\psi |^{2 \si} \psi = 0 \qquad \mbox{ in } \R \times \R^N, 
\eeq
where $ \epsilon = \text{sgn} (\beta ) \in \{ -1,0, 1\}$ and $ \vartheta = \text{sgn} (\g ) \in \{ -1,  1\}$.
By analogy to the usual non-linear Schr\"odinger equation, the case $ \g > 0 $ (or $\vartheta = 1$) is called {\it defocusing, } and the case $ \g < 0 $ (or $\vartheta = -1$) is called {\it focusing.}

Equations (BNLS) and  (\ref{BNLSg}) are Hamiltonian. Two important quantities are conserved by the flow associated to   (\ref{BNLSg}):
the "mass" $\| \psi ( t, \cdot ) \|_{L^2}^2 $,  and the "energy" 
$$
E( \psi ) = \int_{\R^N} |\Delta \psi | ^2\, dx - 2  \epsilon \int_{\R^N} |\nabla \psi | ^2\, dx
+ \frac{\vartheta}{\si + 1}  \int_{\R^N} |\psi |^{ 2 \si + 2} \, dx. 
$$
The natural "energy space" associated to (\ref{BNLSg}) is $ H^2( \R^N)$. 
Equation (\ref{BNLSg}) is {\it mass-critical} for $ \si = \frac{4}{N}$, and {\it energy-critical} 
when $ N \ges 5 $ and $ \si = \frac{4}{N-4}$ (this corresponds to $ 2 \si + 2 = 2^{**}$, where 
$  2^{**}= \frac{2N}{N-4}$ is the  Sobolev exponent satisfying $ \| u \|_{L^{2^{**}}} \les C\| \Delta u \|_{L^2}$ for any $ u \in H^2( \R^N)$).

\medskip

Equation (\ref{BNLSg}) admits an important class of special solutions, the {\it standing waves}.
These are solutions of the form $ \psi(t, x ) = e^{- i \omega t } u (x)$, where $ \omega \in \R$ and $ u $ is a complex-valued function.
They appear as a balance between non-linearity and dispersion and are supposed to play an important role in the dynamics. The standing wave profile $u$ satisfies the equation
\beq
\label{SW}
 \Delta ^2 u + 2 \epsilon \Delta u + \omega u + \vartheta |u |^{2 \si}u = 0 \qquad \mbox{ in }  \R^N. 
\eeq
Solutions of (\ref{SW}) are critical points of the {\it action}
\beq
\label{Sw}
\Sr_{\omega} ( u  ) := \int_{\R^N} |\Delta u |^2-  2 \epsilon  |\nabla u |^2
+ \omega |u |^2  + \frac{\vartheta}{\si + 1}   |\psi |^{ 2 \si + 2} \, dx
= E(u ) + \omega \| u \|_{L^2}^2. 
\eeq

Taking into account the Hamiltonian structure of (\ref{BNLSg}), it is natural to search for standing waves as minimizers (or local minimizers) of the energy when the $L^2-$norm  is kept fixed. 
By a standard application of the approach laid down by T. Cazenave and P.-L. Lions \cite{CaL},
the set of solutions obtained in this way is orbitally stable (see Proposition \ref{orbital}). 
Studying the behaviour of the energy with respect to the mass and  to the scaling gives an insight into possible blow-up scenarios. 
We refer to \cite{BoLe} and  \cite{BCGJ-2} for blow-up results. 

\medskip

In the case $ \epsilon = - 1$, the existence of standing waves has been investigated in several papers 
(see \cite{BCGJ}, \cite{BCGJ-2}, \cite{BoCaMoNa}, \cite{BoNa})
by using various  methods,
including minimisation of the energy at fixed mass 
(see Theorem 1.1 p. 3050 in \cite{BoCaMoNa}).
Some qualitative properties of these solutions as well as the orbital stability of the set of minimizers have also been established. 

In the case $ \epsilon =-1 $, it has been observed in \cite{BoNa}, Theorem 1.1 and in \cite{BoCaMoNa}, Theorem 1.2 that it is possible to minimize $ \int_{\R^N} |\Delta u |^2 - 2|  \nabla u |^2 + \omega |u |^2 \, dx$
under the constraint $ \int_{\R^N} |u |^{ 2 \si +  2 } \, dx = 1$ provided that $ \omega > 1$. 
Although this approach gives the existence of standing waves, it is not completely satisfactory because the considered quantities are not conserved by the flow of (\ref{BNLSg}), and consequently it does not give much information about the dynamics of (\ref{BNLSg}).

\medskip

The case $ \epsilon = 1$ (corresponding to $ \beta > 0 $ in (BNLS)) is more difficult and, as far as we know, there are no satisfactory results in the literature concerning the minimisation of the energy at fixed $L^2-$norm.
Our aim is to clarify this situation in the focusing case ($\vartheta = -1$ in (\ref{BNLSg}) or $ \g < 0 $ in (BNLS)).
In the sequel we will always assume that  $ \epsilon =  1$,  although most of our results are still valid 
if $ \epsilon = 0 $ or if $ \epsilon =-1$. Rewriting our proofs with $ \epsilon = - 1$  would give alternate proofs of some results in 
\cite{BCGJ}, \cite{BCGJ-2}, \cite{BoCaMoNa}, \cite{BoNa}.

To be more precise, we focus our attention on the minimisation problem 
$$
\label{Pula}
\begin{array}{c}
\mbox{ minimize } 
\ds {E} (u):= \int_{\R^N}|\Delta u|^2\, dx - {2}\int_{\R^N}|\nabla u|^2\, dx - \frac{1}{\si + 1}  \int_{\R^N}|u|^{2\sigma+2}\, dx
\\
\\
\quad \mbox{ in the set } 
\ds S(m ): = \left\{ u\in H^2(\R^N) \; \Big| \; \int_{\R^N}|u|^{2}\, dx=m \right\} .
\end{array}
\leqno{({\Po}_{m})}
$$
We  denote 
\beq
\label{EMIN}
E_{min}(m ) : = \inf \{ E(u ) \; | \; u \in S( m ) \}. 
\eeq

\medskip

The basic properties of the function $E_{min}$ are given in Proposition \ref{Emin}. 
In particular, we show that $E_{min}(m) $ is finite for any $ m > 0 $ if $ N \si < 4$, and 
$E_{min}(m) = -\infty $ for any $ m > 0 $ if $ N \si > 4$
(of course, this is related to the fact that (\ref{BNLSg}) is mass-critical for $ \si = \frac 4N$). 
If $ N \si = 4 $, there exists some $ k_* > 0 $ such that $E_{min}(m) $ is finite for $ m \in (0, k_*)$ and 
$E_{min}(m) = - \infty $ if $ k \ges k_*$. 
A simple scaling argument shows that we have always $E_{min}(m) \les - m $. 
If $ E_{min}(m) = -m$, the minimisation problem $(\Po_m)$ does not have solutions, and 
 all minimizing sequences converge weakly to zero.
If $ E_{min}(m) < -m$, it is shown in Theorem \ref{Global} that there exist minimizers for $(\Po_m)$ and that 
all minimizing sequences are pre-compact (after translation), which gives the orbital stability of the set of minimizers by the flow associated to (\ref{BNLSg}).
If $ 0 < \si \les \frac 4N$, there exists $ m_0 \ges 0 $ such that $ E_{min}(m ) = -m $ for $ m \in (0, m_0]$ and $ E_{min}(m ) < -m $ for $ m > m_0$. 
It is an important question whether $ m_0 = 0 $ or $ m_0 > 0 $. 
Notice that the presence of standing waves prevents scattering for (\ref{BNLSg}).
Therefore, if $ m_0 = 0 $ we cannot expect a scattering theory for solutions of  (\ref{BNLSg}) having small 
$L^2-$norm.

It is easily seen that for any $ u \in H^2( \R^N)$ with $ \| u \|_{L^2}^2 = m $ we have 
\beq
\label{factorize}
\begin{array}{l}
 E( u ) + \|u \|_{L^2}^2 = \| \Delta u + u \|_{L^2} ^2 -{\frac{1}{\si + 1} \| u \|_{L^{2 \si +2}}^{2 \si + 2}}
 = \| \Delta u + u \|_{L^2} ^2\left( 1 - \frac{1}{\si + 1} \frac{\| u \|_{L^{2 \si +2}}^{2 \si + 2}}{ \| \Delta u + u \|_{L^2} ^2} \right) 
 \\
 \\
 =  \| \Delta u + u \|_{L^2} ^2\left( 1 - \frac{m^{\si }}{\si + 1} \frac{\| u \|_{L^{2 \si +2}}^{2 \si + 2}}{\| u \|_{L^2}^{2 \si}  \| \Delta u + u \|_{L^2} ^2} \right) 
 =   \| \Delta u + u \|_{L^2} ^2\left( 1 - \frac{m^{\si }}{\si + 1}  \Qo(u)^{ 2 \si + 2} \right), 
 \end{array}
\eeq
where
$$
 \Qo(u) = \frac{\| u \|_{L^{2 \si +2}}}{\| u \|_{L^2}^{\frac {\si}{\si + 1}}  \| \Delta u + u \|_{L^2} ^\frac{1}{\si + 1}}.
$$
Let $ M = \ds \sup \{ \Qo ( u ) \; | \;  u \in H^2( \R^N), u \neq 0  \}$.
If $M$ is finite, it follows from (\ref{factorize}) that $  E( u ) + \|u \|_{L^2}^2  \ges 0 $ for any $ u $ satisfying $ \| u \|_{L^2}^2 = m $ 
provided that $m$ is small enough, so that $ 1 - \frac{m^{\si }}{\si + 1} M^{ 2 \si + 2 } \ges 0 $. 
This shows that $ E_{min}(m) \ges - m $ for sufficiently small $m$, and consequently $(\Po_m)$ does not admit minimizers for small $m$. 
If $ M = \infty$,  then for any $ m > 0 $ we may find $ u \in H^2( \R^N)$,  $ u \neq 0 $  such that $ 1 - \frac{m^{\si }}{\si + 1}  \Qo(u)^{ 2 \si + 2} < 0 $. 
Then taking $ v = \frac{\sqrt{m}}{\| u \|_{L^2} } u $ we see that $  \| v \|_ {L^2} = m $, 
$ Q( v ) = Q( u)$, and (\ref{factorize}) gives $ E(v ) + \|v \|_{L^2}^2 < 0 $, which implies $ E_{min}(m) < - m $. 
Notice that $M$ is finite if and only if the   Gagliardo-Nirenberg-type inequality 
$$
\| u \|_{L^{2 \si + 2}} \les C \| u \|_{L^2}^{\frac {\si}{\si + 1}}  \| \Delta u + u \|_{L^2} ^{\frac{1}{\si + 1}} 
$$
holds true for all $ u \in H^2 ( \R^N)$. 
Obviously, $M$ is the best possible constant in this inequality. 

We will study slightly more general inequalities, namely we will investigate whether there exists $ C > 0 $ such that 
$$
\| u \|_{L^{p}} \les C \| u \|_{L^2}^{\ka}  \| \, | D|^s u - u \|_{L^2} ^{1 - \ka} \qquad
\mbox{ for all } u \in H^s ( \R^N),  
$$
where $ p \in (2, \infty)$, $ \ka \in (0, 1 )$ and $ |D|^s $ is the Fourier integral operator given by 
$ |D|^s u = \Fo ^{-1} \left( | \cdot |^s \Fo ( u) \right)$. 
This leads us to study the boundedness  of the quotient 
$$
Q_{\ka }(u ) = \frac{\| u \|_{L^{p}}}{  \| u \|_{L^2}^{\ka}  \| \left( |D|^s - 1\right) u  \|_{L^2} ^{1 - \ka} }.
$$
when   $u \in  H^s( \R^N) \setminus \{ 0 \} $. We obtain the following result. 

\begin{Theorem}
\label{main-ineq}
Let $ N \in \N^*$,   
$p \in (2, \infty)$, $ \ka \in (0, 1)$,  and $ s > 0 $. 
Then $ Q_{\ka}$ is bounded on $ H^s( \R^N) \setminus \{ 0 \}$ 
if and only if 
\beq
\label{condition}
 \ka \ges \frac12  \quad \mbox{  and }  \quad
 \frac{N}{s}\left(\frac 12 - \frac{1}{p} \right)  \les 1 - \ka \les    \frac{N+1}{2}\left(\frac 12 - \frac{1}{p} \right).
\eeq

\end{Theorem}

To prove Theorem \ref{main-ineq} we are led to develop an original approach,
 based on the Hausdorff-Young inequality in space dimension $ N = 1$, and on the Tomas-Stein inequality in higher dimensions. 
This method is not limited to the study of $Q_{\ka}$ here above. It is much more general and can be used to prove 
non-homogeneous Gagliardo-Nirenberg inequalities of the form 
$$
\| u \|_{L^p} \les C \| P_1(D) u \|_{L^2} ^{\ka } \| P_2(D) u \|_{L^2} ^{1 - \ka } , 
$$
where  $ P_1(D)$ and $P_2(D)$ are Fourier integral operators defined by 
$
P_i (D) (u ) = \Fo ^{-1} \left( P_i ( \cdot ) \Fo({u}) \right).
$
See Remark \ref{R-generalization}. 
Some quantitative variants are also available: see Remark \ref{qualitative-generalization}.

Using Theorem \ref{main-ineq} with $ s = 2$, $ p = 2 \si + 2$, and $ \ka = \frac{\si}{\si + 1}$ we infer that the quotient $ \Qo $ in (\ref{factorize}) is bounded on $ H^2 ( \R^N) \setminus \{ 0 \} $ if and only if 
$  \max \left( 1, \frac{4}{N+1} \right) \les \si \les \frac 4N$
(see Proposition \ref{P-interesting}).

As a matter of fact, our method works in the simpler case when $ \epsilon = - 1$ in (\ref{SW}). 
Proceeding as in (\ref{factorize}) we write
$$  E( u ) =  \left( \| \Delta u \|_{L^2}^2+ 2 \| \nabla u \|_{L^2} ^2 \right) \left( 1 - \frac{\|u \|_{L^2}^{2\si }}{\si + 1}  \Qo(u)^{ 2 \si + 2} \right),$$
 where 
$$
 \Qo(u) = \frac{\| u \|_{L^{2 \si +2}}}{\| u \|_{L^2}^{\frac {\si}{\si + 1}} 
  \left( \| \Delta u \|_{L^2}^2+ 2 \| \nabla u \|_{L^2} ^2 \right) ^\frac{1}{\si + 1}} ,
$$ 
and we need to study the boundedness of the quotient $\Qo $ 
to decide whether $E_{min}(m) < 0  $ for all $ m > 0 $ or $ E_{min}(m) = 0 $ for small $m$. 
In this way  it is possible to give an alternate (and shorter) proof of { Theorem 1.1 p. 5030} in \cite{BoCaMoNa}.

\medskip

Having at hand Theorem \ref{main-ineq}, 
we establish the existence of solutions to the problem $(\Po_m)$ under optimal assumptions. 
We use  some  ideas in  \cite{MM-R}, \cite{MM} and \cite{M3}, 
but all our proofs are self-contained. 
The next Theorem summarizes our main results  on  the existence of minimizers for $(\Po _m)$.  It covers all possible situations.

\begin{Theorem}
\label{GlobalExistence} 
Let $ N \in \N^*$. Let $ E_{min}$ be as in (\ref{EMIN}). The following assertions hold true. 

\medskip

(i) If $ 0 < \si < \max \left( 1, \frac{4}{N+1} \right) $ and $ \si < \frac 4N$ we have $ - \infty < E_{min}(m) < -m $ for all $ m > 0 $.

\medskip

(ii) If $ \max \left( 1, \frac{4}{N+1} \right)  \les \si < \frac 4N$, there exists $ m_0 > 0 $ 
(given by  (\ref{explicitm0})) such that $ E_{min}(m ) = - m $ for all $ m \in (0, m_0]$ and 
$ - \infty < E_{min}(m) < -m $ for any $ m > m_0 $.

\medskip

(iii) If $ \si = \frac 4N$, let $m_0 =0 $ if $ \si < 1 $ and let $m_0 $ be as in (\ref{explicitm0}) if $ \si \ges 1$. Let $ k_*$ be as in Proposition \ref{Emin} (vi). 
Then we have $ m_0 < k_*$ and  $ E_{min}(m ) = - m $ for all $ m \in (0, m_0]$, 
$ - \infty < E_{min}(m) < -m $ for $ m \in ( m_0, k_*)$ and $ E_{min}(m) = - \infty $ for $ m \ges k_*$. 

\medskip

(iv)  If $ \si > \frac 4N$ we have $ E_{min}(m) = - \infty $ for all $ m > 0 $.

\medskip

Problem $(\Po_m)$ admits solutions whenever $ - \infty < E_{min}(m) < -m $;    moreover, 
any minimizing sequence for $(\Po_m)$ has a subsequence that converges strongly in $ H^2 ( \R^N)$ modulo translations. 
Minimizers of $(\Po_m)$ solve (\ref{SW}) for some $ \omega >  1$.

Problem $(\Po_m)$ does not admit minimizers if $ m_0 > 0 $ and  $ m \in (0, m_0)$. 

\end{Theorem}

If $ \epsilon = 1$ {and $ \vartheta = -1$,} as we assume throughout this paper, 
writing  $ \omega = 1 + c $  
 equation (\ref{SW}) becomes
\beq
\label{SWc}
 \Delta ^2 u + 2  \Delta u + (1 + c) u - |u |^{2 \si}u = 0 \qquad \mbox{ in }  \R^N. 
\eeq
As already mentioned, 
solutions of (\ref{SWc}) are critical points of the {\it action} functional 
\beq
\label{Swc}
\begin{array}{l}
\ds S_c ( u  ) := \int_{\R^N} |\Delta u |^2-  2   |\nabla u |^2
+(1 + c) |u |^2 -  \frac{1}{\si + 1}   |\psi |^{ 2 \si + 2} \, dx
\\
\\
\ds = E(u ) + (1 + c) \| u \|_{L^2}^2 
= T_c ( u ) -  \frac{1}{\si + 1} \int_{\R^N} |\psi |^{ 2 \si + 2} \, dx, 
\end{array}
\eeq
where $ T_c ( u ) := \ds \int_{\R^N} |\Delta u |^2-  2   |\nabla u |^2
+(1 + c) |u |^2 \, dx.$
A classical approach to find solutions for (\ref{SWc}) is to show that 
$ t(c) :=  \inf \{  T_c (u) \; | \; u \in H^2 ( \R^N) , \;  \int_{\R^N} |u |^{ 2 \si + 2} \, dx = 1\} $ is achieved. 
In Theorem  \ref{Tc} we prove that  minimizers for $t(c)$ exist for any $ c > 0 $ and any  $ \si \in (0, \infty)$ if $ N \les 4$, respectively any $ \si \in (0, \frac{4}{N-4}) $ if $ N \ges 5$. 
 Moreover,  if $ u $ is  a minimizer for $ t(c)$ then $ v := t(c)^{\frac{1}{2 \si}} u $ solves 
 (\ref{SWc})  and for any other solution $ w \in H^2( \R^N)$ of (\ref{SWc}) we have 
 $ S_c( v ) \les S_c(w)$ (see Proposition \ref{groundstates}); we say that $ v $ is a {\it minimum action solution } of (\ref{SWc}). 
Therefore equation (\ref{SWc})  admits minimum action solutions for any energy-subcritical $ \si $ and for any $ c > 0 $.
The next result shows that minimizers given by Theorem  \ref{GlobalExistence} are minimum action solutions for (\ref{SWc}):

\begin{Theorem}
\label{GlobalProperties} 
Assume that $ 0 < \si \les \frac 4N$.
Let $ u $ be a minimizer for problem $(\Po_m)$, as given by Theorem \ref{GlobalExistence}.
The following properties hold true:

\medskip

(i) There exists some $ c = c(u ) > 0 $ such that $ u $ is a minimum action solution for (\ref{SWc}). 
Furthermore, any minimum action solution of (\ref{SWc}) with $ c = c(u) $ is also a minimizer for  $(\Po_m)$.

\medskip

(ii) If $ m _1 < m_2$, the function $ u_ 1 $ solves $ (\Po _{m_1})$ and $ u_ 2 $ solves $ (\Po _{m_2})$,   then $ c( u_1 ) < c( u_2)$. 

\medskip
(iii)  If $ 0 < \si < \frac 4N$, we have $ c( u ) \lra \infty $ as $ m \lra \infty$.

\medskip
(iv) If $ 0 < \si < \max \left(1, \frac{4}{N+1} \right) $ and $ \si \les \frac 4N$ we have $ c( u ) \lra 0 $ as $ m \lra 0 $. 
If $u_m$ is any solution of the minimisation problem $(\Po_m)$, 
denote $ v_m = \frac{ u_m }{\sqrt{m}} = \frac{ u_m }{\| u _m \|_{L^2}}$, so that $ \| v _m \|_{L^2} = 1$. 
Then we have 
$$
\| \Delta v_m \|_{L^2} \lra 1, \quad 
\| \nabla v_m \|_{L^2} \lra 1, \quad 
\| (\Delta + 1)v_m \|_{L^2} \lra 0 \qquad \mbox{ as } m \lra 0 , 
$$
and $ \| v_m \|_{L^p} \lra 0 $ for any $ p \in (2, \infty)$ if $ N \ges 4$, respectively  for any $ p \in (2, 2^{**})$ if $ N \geq 5$.

\end{Theorem}

It is proven in Proposition \ref{tc-zero} that $ t(c) \les C \sqrt{c}$ as $ c \lra 0 $ and Corollary \ref{estimates} below shows
 that for any energy-subcritical $ \si > 0 $, we have $ t(c ) \sim c ^{1 - \frac{N\si}{4 ( \si + 1)}} $ as $ c \lra \infty$.
The behaviour of minimum energy solutions of  (\ref{SWc}) (and, in particular, the behaviour of minimizers 
for the problem $(\Po_m)$ as $ m \lra \infty$ in the case $ 0 < \si < \frac 4N$) is described in Proposition 
\ref{convergence} and Corollary \ref{estimates}: after rescaling and translation, they converge to minimizers of the functional $ K(u ) := \int_{\R^N} |\Delta u |^2 + |u |^2 \, dx $ under the constraint 
$ \int_{\R^N} |u |^{ 2 \si + 2} \, dx = 1$.

\medskip

Throughout this article, $ H^s( \R^N) = H^s( \R^N, \C)$ denotes the standard Sobolev space of complex-valued 
functions and all minimization problems that  we consider are for complex-valued functions. 
Theorem \ref{main-ineq} is  valid for complex-valued functions, too. Exactly the same statement holds if we restrict ourselves to real-valued functions. 
It turns out that all minimizers of the problem ($\Po_m$) above are of the form $ e^{ i a } \tilde{u} (x)$, where $\tilde{u} $ is  a real-valued function; see Propositions \ref{groundstates} and \ref{PomToc}.

\medskip

Let us stress that we state and prove Theorems \ref{GlobalExistence} and  \ref{GlobalProperties} for pure power nonlinearities only for simplicity. 
One can replace the term $ \frac{1}{\si + 1} \int_{\R^N} | u |^{ 2 \si + 2} \, dx $ in the expression of $E(u)$ by 
$\int_{\R^N} F(u) \, dx $. 
Under appropriate assumptions on $F$, similar results can be obtained with minor changes in the proofs.
In order to avoid technicalities and to keep this paper reasonably long, we have chosen not to pursue in this direction.

\medskip

In the case  $ \si  > \frac 4N$ we have $ E_{min}(m) =  - \infty $ for any $ m > 0 $ and  the minimization problem $(\Po_m)$ does not make sense.  
In this case we  investigate the existence of {\it local } minimizers of $E$ when the $L^2-$norm is kept fixed. 
By {\it local minimizer } we mean a function $ u \in H^2( \R^N)$ such that there exists an open set $ \Uo \subset H^2( \R^N)$ 
having the property  that $ u \in \Uo $ and $ E(u ) = \inf \{ E(v) \; | \; v \in \Uo \mbox{ and } \| v \|_{L^2} = \| u \|_{L^2}\}$.

We find an open set $ \Oo \subset H^2( \R^N)$ (described in (\ref{O})) such that any possible local minimizer of $E$ at fixed $L^2-$norm must belong to $ \Oo$. 
The set $ \Oo \cup \{ 0 \} $ is star-shaped and is an open neighbourhood of the origin in 
$H^2 ( \R^N)$,  and $ \Oo $ is unbounded in $ H^2( \R^N)$. 
We denote 
$$
\tilde{E}_{min} ( m ) = \inf \{ E(u) \; | \; u \in \Oo \; \mbox{ and }  \| u \|_{L^2}^2 = m \}. 
$$
The problem of finding minimizers for $ \tilde{E}_{min} ( m ) $ (which are "best possible"
local minimizers of $E$ when the $L^2-$norm is fixed)
in the mass-supercritical case $ \si > \frac 4N$  is  much 
 harder than finding global minimizers for $ E_{min}$ in the subcritical and critical cases. 
To the best of our knowledge this problem has not been addressed in the literature. 
Understanding the behaviour of $E$ with respect to the $L^2-$norm is also an  important step in 
understanding  the dynamics associated to (\ref{BNLSg}).
Our main results in the case $ \si > \frac 4N$ are given below. 

\begin{Theorem}
\label{Localmin}
Suppose that $ \si > \frac 4N$ and $ \si < \infty $ if $ N \les 4$, respectively $ \si < \frac{4}{N-4}$ if $ N \ges 5 $. The following assertions are true. 

\medskip

(i) Assume that  $ N \ges 5 $ and $ \frac 4N < \si < 1$. 
Then  we have $  \tilde{E}_{min}(m) < -m$ for any $ m > 0 $.

\medskip

(ii) If $ \frac 4N < \si $ and $ \si \ges 1$, there exists $ m_0 > 0 $ such that $  \tilde{E}_{min}(m) = -m$ for any $ m \in ( 0, m_0] $ and the infimum $  \tilde{E}_{min}(m)  $ is not achieved 
for any $ m \in (0, m_0)$.

\medskip

(iii) Assume that $ 0 < m < \mu_0 $, where $ \mu _0$ is given by (\ref{m0}),  and $ \tilde{E}_{min}(m) < -m$. 
Then $ \tilde{E}_{min}(m) $ is achieved and any minimizing sequence for  $ \tilde{E}_{min}(m)$ has a subsequence that converges strongly in $H^2(\R^N)$ modulo translations. 

\medskip

(iv) Any minimizer for $ \tilde{E}_{min}(m)$ solves (\ref{Swc}) for some $ c = c(u) $ satisfying 
$$ 0 < c < -1 + \frac{(N \si - 2)^2}{N \si ( N \si - 4)} +  \frac{8( N \si - 2)}{N ( N \si - 4 )^2}.$$
Moreover, if $ N \ges 5 $ and $ \frac 4N < \si < 1$
we have $ c(u ) \lra 0 $ as $ m \lra 0 $.

\end{Theorem}

Any solution $u$  of { (\ref{SWc})} provided by Theorem \ref{Localmin} must satisfy (\ref{ordered}), 
and consequently there is some explicit constant $ C > 0 $ such that 
$ \| u \|_{H^2} \les C \| u \|_{L^2}$. 
Therefore if  $ N \ges 5 $ and  $ \frac 4N < \si < 1$ equation (\ref{BNLSg}) admits  standing waves 
with small $H^2-$norm and this rules out a scattering theory for small solutions  of  (\ref{BNLSg}). 
It is an open question whether small solutions of  (\ref{BNLSg}) scatter or not in the remaining cases. 

In the case $ \si > \frac 4N$, the least energy solutions of (\ref{SWc}) given by Proposition \ref{groundstates} have small $L^2-$norm 
as $ c \lra \infty$, but they have large $H^2-$norm  and do not belong to the set $ \Oo$ (see Remark \ref{double} and Corollary \ref{estimates}).  
Thus we have two types of interesting standing waves with small $L^2-$norm: 
the minimum action solutions for $ c \lra \infty$, 
and the local minimizers provided by Theorem \ref{Localmin}.

\medskip

Let us compare our results to  similar results in the cases $ \epsilon = -1 $ and $ \epsilon = 0$.
Let us  consider the problem $(\Po_m)$ with $E$ replaced by 
\beq
\label{sign-energy}
E(u ) =  \int_{\R^N}|\Delta u|^2\, dx - {2}\epsilon \int_{\R^N}|\nabla u|^2\, dx - \frac{1}{\si + 1}  \int_{\R^N}|u|^{2\sigma+2}\, dx .
\eeq
We define $E_{min}$ as in (\ref{EMIN}). 
Then { Theorem 1.1 p. 5030} in \cite{BoCaMoNa} and  Theorem 1.2 p. 2170 in \cite{BCGJ} give the following result:

\medskip

\noindent
{\bf Theorem. } (\cite{BoCaMoNa, BCGJ})  {\it Assume that $ \epsilon =- 1$. Then:

\medskip

(i) If $ 0 < \si < \frac 2N$, we have $ -\infty < E_{min}(m) < 0 $ for any $ m > 0 $. 

\medskip

(ii) If $  \frac{2}{N} \les \si < \frac 4N$, there exists $ m_{cr} > 0 $, depending on $ \si $ and $N$,  such that 
$ E_{min}(m) = 0 $ for any $ m\les m_{cr}   $ and  $- \infty<  E_{min}(m) < 0 $ for any $ m > m_{cr} $. 

\medskip

(iii) If $ \si = \frac 4N$, there exists $ m_{cr} > 0 $ such that $ E_{min}(m) = 0 $ for any $ m\les m_{cr}$ and  $  E_{min}(m) = - \infty$ for  $ m > m_{cr} $, and  $ E_{min}(m)$ is never achieved.

\medskip

(iv) If $ \si > \frac 4N$ we have $ E_{min}(m)  = - \infty$ for all $ m > 0 $.

\medskip

The problem $(\Po _m)$ admits solutions whenever $ -\infty < E_{min}(m) < 0 $. 
Moreover, all minimizing sequences have  subsequences that converge strongly in $H^2(\R^N)$ (modulo translations).  }

\medskip

Quite remarkably, Proposition 2.8 (ii) p. 5038 in \cite{BoCaMoNa} shows  that for $ \si \in\left( \frac 2N, \frac 4N \right)$ and for $ m = m_{cr}$, 
 problem $(\Po_{m_{cr}})$ admits solutions  
despite the fact that  there exist     minimizing sequences that do not have any convergent subsequence (modulo translations). 
In fact, proceeding as in (\ref{factorize}) it is easily seen that in this case one has
$$
\frac{\si + 1}{m _{cr} ^{\si }} = \sup \left\{ 
\frac{ \| u \|_{L^{2 \si +2}}^{2 \si + 2}}{\| u \|_{L^2}^{2 \si} \left( \| \Delta  u \|_{L^2} ^2 + 2\| \nabla  u \|_{L^2} ^2 \right) } 
\; \Big| \; u \in H^2( \R^N) \setminus \{ 0 \} \; \right\}
$$
and that $u$ is a minimizer for $(\Po_{m_{cr}})$ if and only if $ \| u \|_{L^2}^2 = m_{cr}$ and 
$u$ is an optimal function for the non-homogeneous Gagliardo-Nirenberg inequality
$$ \| v \|_{L^{ 2 \si + 2}} \les C \| v \|_{L^2}^{\frac{\si}{\si + 1}} \left(  \| \Delta  v \|_{L^2} ^2 + 2\| \nabla  v \|_{L^2} ^2 \right) ^{\frac{1}{2 \si + 2}}. $$

\medskip

We stress that local minimizers of the energy at fixed mass are specific to the case $ \epsilon = 1$. 
Such solutions do {\it  not } exist if $ \epsilon \les 0 $, see Remark \ref{no-local-min}. 
If $ \epsilon = -1$  it is shown in \cite{BCGJ}, Theorem 1.3 p. 2171  that one can minimize $E$ 
in the set $ \{ u \in H^2( \R^N) \; | \; \| u \|_{L^2}^2 = m \mbox{ and } P_1( u ) = 0 \}$ for some values of $ m > 0$, where $P_1$  is a Pohozaev-type functional  given in  (\ref{P1}). 
The minimizers found in \cite{BCGJ} are minimum action solutions of (\ref{SW}), 
but are not minimizers of $E$ at fixed $L^2-$ norm.
They correspond to solutions given by Theorem \ref{Tc} and Proposition \ref{groundstates} below. 
The  instability by blow-up of such  minimum action solutions has been proven in \cite{BCGJ-2}, 
Theorem 1.1 provided that they are radial and $ \frac 4N \les \si \les 4$ (and $ \si < \frac{4}{N-4}$ if $ n \ges 5$), and that instability result is an indication that those solutions cannot be local minimizers of the energy at fixed $L^2-$norm.

\medskip

The case $ \epsilon = 0 $ is much simpler. Proceeding as in (\ref{factorize}) we find 
\beq
\label{epsilon=0}
 E( u ) 
 =  \| \Delta u  \|_{L^2} ^2\left( 1 - \frac{\| u \|_{L^2}^{2 \si }}{\si + 1} \frac{\| u \|_{L^{2 \si +2}}^{2 \si + 2}}{\| u \|_{L^2}^{2 \si}  \| \Delta u \|_{L^2} ^2} \right) 
 =   \| \Delta u \|_{L^2} ^2\left( 1 - \frac{m^{\si }}{\si + 1}  \Qo_0(u)^{ 2 \si + 2} \right)
\eeq
for any $ u \in H^2( \R^N) $ such that $ \| u \|_{L^2}^2 = m $, where $ \Qo _0 ( u )  
= \frac{ \| u \|_{L^{ 2 \si + 2}}}{\| u \|_{L^2}{\frac{\si }{\si + 1}} \| \Delta u \|_{L^2}{\frac{1 }{\si + 1}} }$. 
A simple scaling argument shows that the quotient $ \Qo _0$ is unbounded on $ H^2( \R^N) \setminus \{0 \}$ if $ \si \neq \frac 4N$. With the above notation we get: 
\begin{itemize} 
\item[(i)] $ - \infty < E_{min}( m ) < 0 $ for any $ m > 0 $ if $ 0 < \si < \frac 4N$. 
\item[(ii)]  If $ \si = \frac 4N$, there exists $ m_{cr} > 0 $ such that $ E_{min}(m) = 0 $ for any $ m\les m_{cr}$ and  $  E_{min}(m) = - \infty$ for  $ m > m_{cr} $. 
$  E_{min}(m_{cr})$ is achieved by some optimal function for the    {Sobolev-Gagliardo-Nirenberg} inequality 
(\ref{GNS}), and $  E_{min}(m)$ is never achieved if $ m \neq m_{cr}$. 
\item[(iii)] $ E_{min}( m ) = - \infty $ if $ \si > \frac 4N$. 
\end{itemize}
Assertions (i) and (iii) are proven exactly as statements (v) and (i), respectively, in Proposition \ref{Emin} below, and (ii) follows from (\ref{epsilon=0}). 
The existence of minimizers for any $ m > 0 $ in case (i) is standard (one may use a simplified version of the proof of Theorem \ref{Global}).
 
 If $ \epsilon \les 0 $, equation (\ref{SW}) has infinitely many solutions that can be obtained by using topological methods  (see, e.g., Theorem 1.4 p. 2172 in \cite{BCGJ}). 
This is presumably true  for $ \epsilon = 1$, too. In this paper we focus on standing waves that minimize the energy at fixed mass, which are the most important for the dynamical study of (\ref{BNLSg}).

\medskip

Let us mention that the Cauchy problem for (\ref{BNLSg}) 
has been considered in several articles; see  \cite{Paus07} and references therein. 
In the energy-subcritical case ({that is, $ N \les 4 $} and $ \si \in (0, \infty)$, or $ N \ges 5 $ and $ 0 < \si < \frac{4}{N-4}$), B. Pausader  proved local existence in $ H^2( \R^N) $ as well as the conservation of mass and energy in all cases (see Proposition 4.1 p. 204 in \cite{Paus07}).


{In the defocusing case ($\vartheta = 1$) B. Pausader also proved global existence  for any $ \epsilon \in \{ -1, 0 , 1\} $ and all initial data (Corollary 4.1 (a) p. 205 in \cite{Paus07}), and scattering provided that $ \epsilon \in \{-1, 0 \}$, $ N \ges 5$  and $ \frac 4N < \si < \frac{4}{N-4}.$ In low dimensions $ 1 \les N \les 4$, scattering  has been proved in \cite{Paus-Xia} (see Theorem 1.1 p. 2177 in \cite{Paus-Xia}) provided that $\epsilon \in \{-1,0\}$ and $ \si > \frac 4N$. The latter condition can be weakened to $ \si > \frac 2N$ if $ \epsilon = -1$; this is due to the fact that Strichartz estimates are better for $ \epsilon = -1$. }

{In the focusing case ($\vartheta = -1)$, global existence holds provided that 
$ \si $ is energy-subcritical and the initial data is sufficiently small in $H^2(\R^N)$, 
or $ \si < \frac 4N$ and the initial data is arbitrary, 
or $ \si = \frac 4N$ and the initial data is sufficiently small in $ L^2(\R^N)$ (Corollary 4.1 (b)-(d) p. 205 in \cite{Paus07}). Global existence in the energy-critical case $ N \ges 5$ and $ \si = \frac{4}{N-4}$ in also  shown for any radial initial data (Theorem 1.1 p. 198 in \cite{Paus07}) and for non-radial but small initial data, as well as scattering for radial data if $ \epsilon \in \{ -1, 0 \}$. 
Notice that in the case $ \epsilon = 1$, $ N \ges 5 $ and $ \frac 4N < \si < 1$,  Theorem \ref{Localmin} above and the estimate (\ref{ordered}) below imply the existence of  standing waves with small $H^2-$norm, and these solutions do not scatter. 
}

\medskip

This paper is organized as follows.
In the next section we develop a  method to deal with non-homogeneous Gagliardo-Nirenberg inequalities and we prove Theorem \ref{main-ineq}.
We separate the cases  $ N = 1$ (when a simple argument based on the Hausdorff-Young inequality is sufficient, see Theorem  \ref{T-1D}) 
and $N\ges 2$ (when a more involved argument relying on the Tomas-Stein inequality is needed, see Theorem \ref{T-multiD}).
Examples \ref{Ex-dim1} and \ref{knapp} show that the results we obtain are optimal for the operator $ |D|^s - 1$. 
Extensions to more general operators are indicated in Remarks \ref{R-generalization} and \ref{qualitative-generalization}. 

In Section \ref{Globalmin} we consider the problem $(\Po_m)$ and we prove Theorems \ref{GlobalExistence}  and  \ref{GlobalProperties}.
Statements (i)-(iv) in Theorem \ref{GlobalExistence}   follow from Propositions \ref{Emin} and \ref{P-interesting}. The existence of minimizers and the pre-compactness of minimizing sequences are given by 
Theorem \ref{Global}.
Theorem \ref{GlobalProperties} (i)-(iii) follows from Propositions \ref{Lagrange} and \ref{PomToc}
(see also Remark \ref{PropEmin} and Proposition \ref{P-interesting}), 
and  Theorem \ref{GlobalProperties} (iv) is Corollary \ref{asympt-zero}.
Some asymptotic properties of minimum action solutions  as $ c \lra 0 $ and as $ c \lra \infty $ are given in Propositions \ref{assympto-zero} and  \ref{convergence}, respectively,  and in Corollary \ref{estimates}.

In Section \ref{Local} we consider the more delicate problem of minimizing  the energy at fixed $L^2-$norm in the set $ \Oo $ when  $ \si > \frac 4N$ and we prove Theorem \ref{Localmin}.
Statement (i) and the first part of (ii) in Theorem \ref{Localmin} follow from Lemma \ref{merdique}, and the second assertion in  (ii) follows from Remark \ref{nolocmin}.
Part (iii) is Theorem \ref{T2.6}. 
For (iv), see Remark \ref{loc-Lagrange}.

\section{A class of non-homogeneous Gagliardo-Nirenberg inequalities}
\label{inequalities}

Let $(X, \Ao, \mu) $ be an arbitrary measure space and let $f$ be a complex-valued measurable function on $X$.
We say that $ f \not\equiv 0 $ if the set $ \{ x \in X \;| \; f( x )\neq 0 \}$ has positive measure. 
Obviously, if $ f \not\equiv 0 $ then $ \int_X |f|^{\al } \, d \mu > 0$ for any $ \al > 0 $. 
If $ f_1, f_2 $ are measurable and $ f_1 f_2 \not\equiv 0$, then  $ f_1 \not\equiv 0 $ and $ f_2 \not\equiv 0 $.
We use the convention $ \frac{ f_1(x) }{f_2(x)} = 0 $ if $ f_1(x) = f_2(x) =  0$.

The following elementary lemma will be very useful in the sequel. 

\begin{Lemma}
\label{measure}
Let $(X, \Ao, \mu)$ be an arbitrary measure space, let   $ q \in [1, 2 ) $ and $ \ka \in (0, 1)$. 
Consider three measurable functions $ w, w_1, \, w_2: X \lra \C $ 
 such that $ w w_1 \not\equiv 0$,  $ w w_2 \not\equiv 0 $ and 
 $ w = 0 $ on the set $ \{ x \in X \; | \; w_1 (x ) = 0 \mbox{ and } w_2 ( x ) = 0 \}$.
Let 
$$
M_1 = \sup \left\{ 
\frac{\| \ph w \|_{L^q}}{ \| \ph w_1\|_{L^2}^{\ka } \cdot \| \ph w_2 \|_{L^2}^{ 1 - \ka }} \; \Big| \;
\begin{array}{l}
\ph \mbox{ is measurable, }
\ph w_i \in L^2 (X),  \\   \ph w_i \not\equiv 0 \mbox{ for } i =1, \, 2
\end{array}
 \right\}, 
 \mbox{ and }
$$
$$
M_2 = \ds \sup_{t > 0 } \left( t ^{\frac{ 1 - \ka}{2}} 
\Big\| \frac{ w}{\left(w_1 ^2 + t | w_2 |^2 \right)^{\frac 12}} \Big\|_{L^{ \frac{2q}{2 -q }}} \right).
$$
Then $ M_1 \les    ( 1 - \ka )^{ \frac{ \ka-1 }{2} } \ka ^{- \frac{\ka}{2} } M_2$.
Moreover, if $( X, \Ao, \mu)$ is $ \si-$finite or 
if there exists $ t_* > 0 $ such that 
 $ w \left( w_1 ^2 + t_* | w_2|^2 \right)^{- \frac 12} \in L^{ \frac{2q}{2 - q}}(X)$, then
 $ M_1 =    ( 1 - \ka )^{ \frac{ \ka-1 }{2} } \ka ^{- \frac{\ka}{2} } M_2$.

\end{Lemma}

\begin{remark} 
\label{R1.2}
\rm $\; $ 

(i)  Let  $ A_i = \{ x \in X \; | \; w_i (x) = 0 \}$ for $ i = 1, 2$.
We have 
$ \Big\| \frac{ w}{\left(w_1^2 + t | w_2 |^2 \right)^{\frac 12}} \Big\|_{L^{ \frac{2q}{2 -q }}(X)} 
\ges \Big\| \frac{w}{| w_1|} \Big\|_{L^{ \frac{2q}{2 -q }}(A_2)}$, 
and
$ \Big\| \frac{ w}{\left(w_1^2 + t | w_2 |^2 \right)^{\frac 12}} \Big\|_{L^{ \frac{2q}{2 -q }}(X)} 
\ges t^{ - \frac 12}\Big\| \frac{w}{| w_2|} \Big\|_{L^{ \frac{2q}{2 -q }}(A_1)}$.
We infer that if $M_2< \infty$, then necessarily 
$ w = 0 $ a.e. on $A_1 \cup A_2$ (here we use the assumption that $ w= 0 $ on $ A_1 \cap A_2$ and the convention  $\frac 00 = 0 $).

(ii) For any fixed $ a , b> 0$, the function 
$ g(t) = \frac{t^{ 1 - \ka}}{a + tb } = \frac{1}{ a t^{ \ka - 1}+ b t ^{\ka} }$ 
achieves its maximum on 
$(0, \infty ) $ at $t_{max } = \frac{(1 - \ka)a }{\ka b }$ and 
$ g( t_{max}) = \ka ^{\ka}( 1 - \ka )^{1 - \ka } a^{-\ka } b^{\ka - 1}$. 
Hence for any $ x \in X $ such that $ w_1 (x ) w_2 (x ) \neq 0 $ we have 
$\ds  \max _{t > 0 }  \left[ { t^{\frac{1 - \ka}{2}}}{\left(|w_1(x)|^2 + t | w_2 (x)|^2 \right)^{-\frac 12}} \right]
 = \ka ^{\frac{\ka}{2}}( 1 - \ka )^{\frac{1 - \ka}{2} } | w_1(x)|^{-\ka } | w_2(x) |^{  \ka -1 }$. 
We have thus a sufficient condition for the finiteness of $M_2$, namely  $M_2 < \infty $  if $ w = 0 $ a.e. on the set  $A_1 \cup A_2$ and 
$ w |w_1|^{-\ka}  | w_2|^{\ka -1 } \in L^{ \frac{2q}{2 -q }}(X).$

(iii) Assume that there is $ t_* > 0 $ such that
$ w \left( | w_1|^2 + t_* | w_2|^2 \right)^{- \frac 12} \in L^{ \frac{2q}{2 - q}}(X)$.
We have 
 \beq
 \label{parameter}
  \min \left(1, \frac{t}{t_*} \right) \les \frac{| w_1|^2 + t | w_2 |^2}{|w_1|^2 + t_* | w_2 |^2} \les \max \left(1, \frac{t}{t_*} \right) \qquad \mbox {  whenever  }  ( w_1, w_2 ) \neq (0,0).
  \eeq
Since  $ w = 0 $ if $( w_1, w_2 ) = (0,0)$,  we infer that 
 $ w \left( |w_1|^2+ t | w_2|^2 \right)^{- \frac 12} \in L^{ \frac{2q}{2 - q}}(X)$ for any $ t > 0 $. 
{Now let }
$$
F(t ) =
\Big\| \frac{  t ^{\frac{ 1 - \ka}{2}}  w}{\left(w_1 ^2 + t | w_2 |^2 \right)^{\frac 12}} \Big\|_{L^{ \frac{2q}{2 -q }}}^{ \frac{2q}{2 -q }}
= \int_{ X} \, \frac{| w(x)|^{\frac{2q}{2 -q }} }
{\left( t^{ \ka - 1} | w_1(x) |^2 + t^{\ka } | w_2(x)|^2 \right) ^{\frac{q}{2 -q }} } \, d \mu. 
$$
Clearly, $M_2$ is finite if and only if $F$ is bounded from above. 
For any $ x \in X$, the mapping 
$$
t \longmapsto \frac{| w(x)|^{\frac{2q}{2 -q }} }{\left( t^{ \ka - 1} | w_1(x) |^2 + t^{\ka } | w_2(x)|^2 \right) ^{\frac{q}{2 -q }} }
$$
is continuous on $(0, \infty )$. 
Then the estimates (\ref{parameter}) and the dominated convergence theorem imply that $ F$ is continuous on $(a, b ) $ for any $ 0 < a < t_* < b $, hence $F$ is continuous on $ (0, \infty)$. 
In order  to show that $M_2$ is finite, we only have  to prove that $F$ is bounded in a neighbourhood of zero and in a neighbourhood of infinity. 
 
(iv)  If  $( X, \Ao, \mu)$ is $ \si-$finite  and $ M_1 < \infty$, the second part of Lemma \ref{measure} implies that  we must have
$ w ( |w_1|^2 + t | w _2|^2)^{- \frac 12 } \in  L^{ \frac{2q}{2 - q}}(X)$ for all $ t > 0 $.

 \end{remark}

{\it Proof of Lemma \ref{measure}. } 
First notice that for any given $A, B > 0 $, the function $ f(t ) = At^{\ka - 1} + B t^{\ka } $ achieves its minimum on $(0, \infty )$ at $ t_{min} = \frac{( 1 - \ka)A}{\ka B} $ and $ f( t_{min} ) = ( 1 - \ka ) ^{\ka - 1} \ka^{-\ka } A^{\ka } B ^{1 - \ka }. $

Let $ \ph $ be a measurable function such that 
 $ \ph w_i \in L^2(X)$  and $ \ph w_i \not\equiv 0 $ for $ i = 1, \, 2 $.
Using the previous observation with $ A = \| \ph w_1 \|_{L^2}^2 $,  $ B = \| \ph w_2 \|_{L^2}^2$
and $ t_{min} = \frac{( 1 - \ka) \| \ph w_1 \|_{L^2}^2 }{\ka  \| \ph w_2 \|_{L^2}^2} $
we get
\beq
\label{1.1}
\begin{array}{l}
\| \ph w_1 \|_{L^2}^{2 \ka } \cdot \| \ph w_2 \|_{L^2}^{ 2 ( 1 - \ka ) }
= ( 1 - \ka )^{ 1 - \ka } \ka ^{\ka } 
\left( t_{min} ^{\ka -1}  \| \ph w_1 \|_{L^2}^2 + t_{min}^{\ka }  \| \ph w_2 \|_{L^2}^2 \right)
\\
\\
=   ( 1  - \ka )^{ 1 - \ka } \ka ^{\ka } t_{min}^{ \ka - 1} 
 \int_{ X } |\ph |^2 \left( | w_1|^2  +  t_{min} | w_2|^2 \right)  d \mu 
 \\
 \\
=  ( 1  -  \ka )^{ 1 - \ka } \ka ^{\ka } t_{min}^{ \ka - 1}  \big\| \ph \left( |w_1|^2  +  t_{min} | w_2|^2 \right)^{\frac 12} \big\|_{L^2}^2.
\end{array}
\eeq
H\"older's inequality implies that for any two measurable functions $ f, g $ defined on $ X$ there~holds 
\beq
\label{holder}
\| fg \|_{L^q } \les \| f \|_{L^2 } \| g \|_{L^{ \frac{2q}{2 -q }}}.
\eeq
Using (\ref{holder}) with $ f = \ph \left( | w_1|^2 + t_{min} | w_2|^2 \right)^{\frac 12} $ and 
$ g = w \left( |w_1|^2 + t_{min} | w_2|^2 \right)^{- \frac 12} $ we obtain 
\beq
\label{1.3}
\| \ph w \|_{L^q} \les \| \ph \left( | w_1|^2 + t_{min} | w_2|^2 \right)^{\frac 12}  \|_{L^2}
\cdot
\big\|  w \left( |w_1|^2 + t_{min} | w_2|^2 \right)^{- \frac 12}  \big\|_{L^{ \frac{2q}{2 -q }}}.
\eeq
From (\ref{1.1}) and (\ref{1.3}) we get 
\begin{align*}
\frac{\| \ph w \|_{L^q}}{ \| \ph w_1 \|_{L^2}^{\ka } \cdot \| \ph w_2 \|_{L^2}^{ 1 - \ka }} 
& \les  ( 1 - \ka )^{ \frac{ \ka-1 }{2} } \ka ^{- \frac{\ka}{2} } 
 t_{min} ^{\frac{1 - \ka}{2} }  
\big\|  w \left( | w_1|^2 + t_{min} | w_2|^2 \right)^{- \frac 12}  \big\|_{L^{ \frac{2q}{2 -q }}} \\
&\les  ( 1 - \ka )^{ \frac{ \ka-1 }{2} } \ka ^{- \frac{\ka}{2} }  M_2.
\end{align*}

Since the  above chain of inequalities holds 
for any measurable function $ \ph $ such that 
$ \ph w_i \in L^2(X)$  and $ \ph w_i \not\equiv 0 $ for $ i = 1, \, 2$, taking the supremum we get 
$ M_1 \les  ( 1 - \ka )^{ \frac{ \ka-1 }{2} } \ka ^{- \frac{\ka}{2} }  M_2.$

\medskip

Next, assume that there is $ t_* > 0 $ such that
$ w \left( | w_1|^2 + t_* | w_2|^2 \right)^{- \frac 12} \in L^{ \frac{2q}{2 - q}}(X)$.
By Remark \ref{R1.2} (iii) we have
 $ w \left( |w_1|^2+ t | w_2|^2 \right)^{- \frac 12} \in L^{ \frac{2q}{2 - q}}(X)$ for any $ t > 0 $. 

If  a measurable function $ \ph $ satisfies $ \ph w_i\in L^2(X)$  and $ \ph w_i \not\equiv 0 $ for $ i = 1, \, 2$, 
it is obvious that (\ref{1.1}) holds   if we replace  $ t_{min}$ by any $ t > 0 $ 
and the first "$=$"
by "$\les$." In other words, we have 
$$
\| \ph w_1\|_{L^2}^{ \ka } \cdot \| \ph w_2 \|_{L^2}^{ 1 - \ka  }
\les 
 (1  -  \ka )^{ \frac{1 - \ka}{2} } \ka ^{\frac{\ka}{2} } t^{ \frac{\ka - 1}{2}}  \| \ph \left( | w_1|^2  +  t | w_2|^2 \right)^{\frac 12} \|_{L^2}
 \qquad \mbox{ for any } t > 0 
$$
and consequently
\beq
\label{1.4}
\frac{\| \ph w \|_{L^q}}{ \| \ph w_1 \|_{L^2}^{\ka } \cdot \| \ph w_2 \|_{L^2}^{ 1 - \ka }}
\ges  ( 1 - \ka )^{ \frac{ \ka-1 }{2} } \ka ^{- \frac{\ka}{2} } t^{\frac{1 - \ka}{2}} 
\frac{\| \ph w \|_{L^q}}{  \big\| \ph \left( |w_1|^2  +  t | w_2|^2 \right)^{\frac 12} \big\|_{L^2} }
 \qquad \mbox{ for all } t > 0 .
\eeq
For any $ z \in \C$ we denote $ \mbox{sgn}(z) = 0 $ if $ z = 0 $ and $  \mbox{sgn}(z) = \frac{z}{|z |}$ if $ z \neq 0 $. 
Let $ \psi = \ov{  \mbox{sgn}(w)} | w|^{\frac{q}{2 - q}} \left( |w_1|^2 + t | w_2|^2 \right)^{- \frac{1}{2 - q}}. $
Then $ \psi w_i \not\equiv 0 $ because $ w w_i \not\equiv 0$,  and 
$$ 
| \psi w_1 | ^2 
\les  |w|^{\frac{2q}{2 - q}} \left( |w_1|^2 + t | w_2|^2 \right)^{- \frac{q}{2 - q}},
\qquad 
 | \psi w_2|^2 \les \frac 1t 
 |w|^{\frac{2q}{2 - q}} \left( |w_1|^2 + t | w_2|^2 \right)^{- \frac{q}{2 - q}}
$$
 hence $ \psi w_i \in L^2(X)$ for $ i = 1, 2$.
 It is easily seen that 
 $$
 \| \psi w \|_{L^q} \! = \! \big\|  w \left( |w_1|^2 + t | w_2|^2 \right)^{- \frac 12} \big\|_{L^{\frac{2q}{2 - q}}}^{\frac{2}{2 - q}},
$$
and
$$ 
\big\| \psi ( |w_1|^2+ t |w_2|^2) ^{\frac 12} \big\| _{L^2} \! = \!
\big\|  w \left( |w_1|^2 + t | w_2|^2 \right)^{- \frac 12} \big\|_{L^{\frac{2q}{2 - q}}}^{\frac{q}{2 - q}}.
$$
Using (\ref{1.4}) with $\varphi = \psi$,  we discover 
\beq
\label{1.5}
M_1 \ges \frac{\| \psi w \|_{L^q}}{ \| \psi w_1\|_{L^2}^{\ka } \cdot \| \psi w_2 \|_{L^2}^{ 1 - \ka }}
\ges  ( 1 - \ka )^{ \frac{ \ka-1 }{2} } \ka ^{- \frac{\ka}{2} } t^{\frac{1 - \ka}{2}} 
\big\|  w \left( |w_1|^2 + t | w_2|^2 \right)^{- \frac 12} \big\|_{L^{\frac{2q}{2 - q}}}.
\eeq
Since (\ref{1.5}) holds for any $ t > 0 $ we infer that $ M_1 \ges ( 1 - \ka )^{ \frac{ \ka-1 }{2} } \ka ^{- \frac{\ka}{2} } M_2$.

\medskip

{Finally} assume that $(X, \Ao, \mu)$ is $\si-$finite. 
Consider a collection of sets $(X_n)_{n \ges 1} \subset \Ao $ such that  
$ \mu ( X_n)  < \infty $, 
 $ X_n \subset X_{n+1} $ for all $ n $ and $\ds \cup_{n \ges 1 } X_n = X$.
 Fix $ t > 0 $ and denote
 $ A_n = \left\{ x \in X \; | \; \big| w (x ) ( |w_1(x)|^2+ t | w_2(x) |^2 )^{- \frac 12} \big| \les n \right\} \cap X_n$. 
Then $ A_n \subset A_{n + 1}$,  $ \mu ( A_n ) \les \mu ( X_n ) < \infty $ for any $n$ and 
$ \ds \cup_{n \ges 1 } A_n = X$.
Let $ \psi _n = \ov{  \mbox{sgn}(w)} | w|^{\frac{q}{2 - q}} \left( |w_1|^2 + t | w_2|^2 \right)^{- \frac{1}{2 - q}} \1_{A_n}$. 
For all $n $ sufficiently large we have $ w w_i \1_{A_n} \not\equiv 0 $, and consequently 
$ \psi _n  w_i \not\equiv 0 $.

As above we see that 
$$ 
| \psi _n w_1|^2 \les |w|^{\frac{2q}{2 - q}} \left( |w_1|^2 + t | w_2|^2 \right)^{- \frac{q}{2 - q}} \1_{A_n}
\les n ^{\frac{2q}{2 - q}} \1_{A_n} \quad \mbox{ and } 
$$
$$ | \psi _n w_2| ^2\les \frac 1t 
 |w|^{\frac{2q}{2 - q}} \left( |w_1|^2 + t | w_2|^2 \right)^{- \frac{q}{2 - q}}  \1_{A_n}
 \les  \frac 1t n ^{\frac{2q}{2 - q}}\1_{A_n}, 
$$
hence $ \psi_n w_i \in L^2 ( X)$. 
Proceeding as in the proof of  (\ref{1.5}) with $ \psi _n $ instead of $ \psi $ we get 
$$
 M_1 \ges  ( 1 - \ka )^{ \frac{ \ka-1 }{2} } \ka ^{- \frac{\ka}{2} } t^{\frac{1 - \ka}{2}} 
\big\|  w \left( |w_1|^2+ t | w_2|^2 \right)^{- \frac 12} \|_{L^{\frac{2q}{2 - q}}(A_n)} 
\quad \mbox{ for all $n $ sufficiently large. }
$$
Then letting $ n \lra \infty $ and using the monotone convergence theorem we 
find 
$$ M_1 \ges  ( 1 - \ka )^{ \frac{ \ka-1 }{2} } \ka ^{- \frac{\ka}{2} } t^{\frac{1 - \ka}{2}} 
\big\|  w_1 \left( 1 + t | w_2|^2 \right)^{- \frac 12} \|_{L^{\frac{2q}{2 - q}}}.$$
Since this is true for any $ t > 0$, we have
$ M_1 \ges  ( 1 - \ka )^{ \frac{ \ka-1 }{2} } \ka ^{- \frac{\ka}{2} } M_2$ 
 and Lemma \ref{measure} is proven. 
\hfill
$\Box$

\medskip

We consider the Fourier transform defined by 
$ \Fo ( u ) ( \xi ) = \wh{u}( \xi ) = \int_{\R^N} e^{- i x.\xi} u (x ) \, dx $ if $ u \in L^1( \R^N)$, and extended as usually to tempered distributions. 
We consider the Fourier integral operator $ |D|^s - 1 $ defined by 
$$
(|D|^s - 1 ) u = \Fo ^{-1} \left( ( |\cdot |^s - 1 ) \wh{u}\right).
$$
The space $H^s( \R^N)$ is defined by $ H^s( \R^N) = \{ u \in \So '( \R^N) \; | \; ( 1 + |\cdot |^2)^{\frac s2} \wh{u} \in L^2 ( \R^N) \} . $
Given a tempered distribution $ u $, we have $ u \in H^s ( \R^N)$ if and only if $ \wh{u} \in L^2 ( \R^N)$ and $  ( |\cdot |^s - 1 ) \wh{u} \in L^2 ( \R^N)$. 
Moreover, by Plancherel's identity we have 
\begin{equation}\label{plancherel}
 \| u \|_{L^2( \R^N) } = \frac{1}{( 2 \pi)^{\frac N2} } \| \wh{u} \|_{L^2( \R^N) } \, \mbox{ and } \,
 \| ( |D|^s - 1 ) u \|_{L^2( \R^N) } = \frac{1}{( 2 \pi)^{\frac N2} } \| ( |\cdot |^s - 1 )  \wh{u} \|_{L^2( \R^N) }.
\end{equation}

Let $ p \in (2, \infty)$ and let $ \ka \in (0, 1)$. Define
\beq
\label{Qka}
Q_{\ka }(u ) = \frac{ \| u \|_{L^p}}{\| u \|_{L^2}^{\ka } \| ( |D|^s - 1 ) u \|_{L^2}^{1 - \ka }}
\qquad \mbox{ for all } u \in H^s ( \R^N) \setminus \{ 0 \}, 
\eeq
and
\beq
\label{supr}
M := \ds \sup_{ u \in H^s ( \R^N) \setminus \{ 0 \} } Q_{\ka } ( u ) .
\eeq
We will investigate whether 
$M$  is finite. 
In the one-dimensional case we have the following:

\begin{Theorem}
\label{T-1D}
Assume that $ N = 1$, $ s \in (0, \infty)$,  $ p \in (2, \infty)$,  and $ \ka \in (0, 1)$.
The supremum $M$ in (\ref{supr}) is finite if and only if 
$ \frac{1}{2 s} \les  {\frac{( 1 - \ka )p}{p-2 }} \les \frac 12$. 
\end{Theorem}
If $N=1$, the condition in Theorem \ref{T-1D} is equivalent to  condition (\ref{condition}) in Theorem \ref{main-ineq}.

\medskip

{\it Proof.  } 
Let $N \in \N^*$.     Since $ p > 2$, 
by the Hausdorff-Young Theorem (see, e.g., Theorem 1.2.1 p. 6 in \cite{BerLo}) we have 
\beq
\label{Hausdorff-Young}
\| u \|_{L^p (\R^N)} \les ( 2 \pi)^{ \frac N p - N } \| \wh{u} \|_{L^{ p'}(\R^N)}, 
\qquad \mbox{ where } p'= \frac{p}{p-1}. 
\eeq
Taking into account Plancherel's identity and the fact that $ u \in H^s ( \R^N)$ if and only if 
$ \wh{u} \in L^2 ( \R^N)$ and $  ( |\cdot |^s - 1 ) \wh{u} \in L^2 ( \R^N)$, we infer that 
$ M \les ( 2 \pi)^{ \frac N p - \frac N 2} M_3 $ ,  where 
\beq
\label{HY}
M_3 := \sup \left\{  \frac{ \| \ph \|_{L^{ p'}}}{ \| \ph \|_{L^2}^{ \ka } 
\big\| \left( |\cdot |^s - 1  \right) \ph \big\|_{L^2}^{1 - \ka}} \; \Big| \; 
 \ph  , \; 
 \left( |\cdot |^s - 1  \right) \ph \in L^2( \R^N) \setminus \{ 0 \} \; \right\}. 
\eeq
Hence $ M $ is finite if $ M_3 $ is finite. 
To prove the finiteness of $ M_3$  
we may use Lemma \ref{measure} in $ \R^N$ endowed with the Lebesgue measure, with 
$ w = w_1 = 1 $ and $ w_2 ( \xi ) = | \xi |^s - 1 $. 
By Lemma \ref{measure}, it suffices to show that 
$ M_4 :=  {\ds \sup_{t > 0 }  }\left( t ^{\frac{ 1 - \ka}{2}} 
\Big\| \frac{ 1}{\left(1 + t \left(  \, | \cdot |^s -1 \right)^2 \right)^{\frac 12}} \Big\|_{L^{ \frac{2p'}{2 - p' }}} \right) $ is finite. 
Notice that $  \frac{2p'}{2 - p' } = \frac{2p}{p-2 }$.
Given any $ p > 2$ and any $ t > 0 $, the function $ \xi \longmapsto  \frac{ 1}{\left(1 + t \left(  \, | \cdot |^s -1 \right)^2 \right)^{\frac 12}} $ belongs to $L^{\frac{2p}{p-2 } }(\R^N)$ if and only if 
$ \frac{ 2 s p }{p -2} > N  $ (which is equivalent to  $ p ( N - 2s) < 2N$). 
Let 
\beq
\label{Fdet}
F( t ) = 
\Bigg\| \frac{  t ^{\frac{ 1 - \ka}{2}}  }{\left(1 + t \left(  \, | \cdot |^s -1 \right)^2 \right)^{\frac 12}} \Bigg\|_{L^{ \frac{2p}{p-2 }}} ^{ \frac{2p}{p-2 }}
= \int_{ \R^N} \frac{1}{\left( t^{\ka - 1} + t^{\ka } \left( |\xi |^s -1 \right)^2 \right)^{\frac{p}{p-2 }} }\, d \xi . 
\eeq
In view of Remark \ref{R1.2} (iii), $F$ is continuous on $(0, \infty )$ provided that  $ \frac{2 s p}{p-2 } > N$.

Assume now that $ N = 1$ and $ \frac{1}{2 s} \les  {\frac{( 1 - \ka )p}{p-2 }} \les \frac 12$. 
Then $ \frac{2 s p}{p-2 } \ges \frac{1}{1 - \ka } > 1$, hence  $F$ is continuous on $(0, \infty )$ and 
we need  only to check that $F$ is bounded in a neighbourhood of zero and of infinity. 
We have 
$$ 
F(t) = 2 \int_0 ^{\infty } f( r, t ) \, dr, \quad \textup{ where } \quad 
 f( r, t ) =  \frac{1}{\left( t^{\ka - 1} + t^{\ka } \left( r^s -1 \right)^2 \right)^{\frac{p}{p-2 }} }.
$$
For any fixed $ A > 0 $ we have 
$$
\int_0^A | f( r, t ) | \, dr \les \int_0^A t ^{\frac{( 1 - \ka )p}{p-2 } } \, dr = A t ^{\frac{( 1 - \ka )p}{p-2 } }\lra 0 \quad \mbox{ as } t \lra 0 . 
$$
Choose $ A \ges 2^{\frac 1s}$, so that $ r^s >  r^s -1 > \frac 12 r^s $ for $ r > A$. We have then 
$$
  \frac{1}{\left( t^{\ka - 1} + t^{\ka }  r^{2s}  \right)^{\frac{p}{p-2 }} }
< f ( r, t ) <  \frac{1}{\left( t^{\ka - 1} +  \frac 14 
t^{\ka }   r^{2s}  \right)^{\frac{p}{p-2 }} } \qquad \mbox{ for any } r > A. 
$$ 
Using the change of variable $ r = t^{-  \frac{1}{2s} } y $ we get 
$$ 
0 < \int_{ A} ^{\infty} f( r, t ) \, dr 
<  \int_{A}^{\infty}   \frac{1}{\left( t^{\ka - 1} +  \frac 14 
t^{\ka }   r^{2s}  \right)^{\frac{p}{p-2 }} }  \, dr 
= t^{\frac{ ( 1 - \ka)p}{p-2} - \frac{1}{2s} }\int_{ t^{  \frac{1}{2s} } \! A}^{\infty } 
\; \frac{1}{\left( 1 + \frac 14 y^{2s} \right)^{\frac{p}{p-2}} } \, dy. 
$$
We conclude that $F $ is bounded in a neighbourhood of zero if 
$ \frac{ ( 1 - \ka)p}{p-2} \ges \frac{1}{2s} $.

Let us study the behaviour of $F$ as $ t \lra \infty$. We have
$$
 0 <  \int_{ 2} ^{\infty} f( r, t ) \, dr < \int_{2} ^{\infty}   \frac{1}{\left( t^{\ka } \left( r^s -1 \right)^2 \right)^{\frac{p}{p-2 }} }  \, dr 
 = t^{- \frac{ \ka p}{p-2} }  \int_{2} ^{\infty}  \frac{1}{(r^s - 1 ) ^{\frac{2p }{p-2}} } \, dr
 \lra 0 
  \mbox{ as } t \lra \infty. 
 $$
There exist two positive constants $ c_2,\,  c_2$ such that 
$ c_1 |y| \les \big| \, | 1 + y |^s -1 \big| \les c_2 |y | $ 
for all $ y \in [-1, 1]$. 
Using the change of variable $ r = 1 + y$, the above estimate and then the change of variable $ y = t^{- \frac 12} z $ we get 
$$
\begin{array}{l}
\ds \int_0^2 f( r, t ) \, dr 
= \int_{ -1}^1  \frac{1}{\left( t^{\ka - 1} + t^{\ka } \big| \, | 1 + y |^s -1 \big|^2 \right)^{\frac{p}{p-2 }} } \, dy 
\les 
\int_{-1}^1  \frac{1}{\left( t^{\ka - 1} + t^{\ka } c_1 ^2 y^2 \right)^{\frac{p}{p-2 }} } \, dy 
\\
\\ 
\ds = t ^{ - \frac 12 + {\frac{( 1 - \ka )p}{p-2 }} } 
\int_{ - t^{\frac 12} }^{  t^{\frac 12} } \frac{1}{\left( 1 + c_1 ^2 z^2 \right) ^{\frac{p}{p-2 }} } \, dz
\end{array}
$$
and we infer that $F$ is bounded in a neighbourhood of infinity if 
$  {\frac{( 1 - \ka )p}{p-2 }} \les \frac 12$. 

So far we have proved that $M_3$ is finite (and consequently $M$ is finite)  if 
$ \frac{1}{2 s} \les  {\frac{( 1 - \ka )p}{p-2 }} \les \frac 12$. 
The fact that $ M = \infty $ if one of these two inequalities is not satisfied follows from the next example. 
\hfill
$\Box$

\begin{example}  \rm 
\label{Ex-dim1}
Given  $\tau > 0$, we define 
$u_\tau(x):=e^{ix-\frac{x^2}{2\tau^2}}.$
It is clear that $u \in \mathcal{S}(\R)$ and direct computations give
\beq
\label{ExD1-1}
\wh{u}_{\tau}(\xi) = \sqrt{2\pi} \tau e^{-\frac{\tau^2 (\xi-1)^2}{2}}, \quad 
\|u_{\tau}\|_{L^2}^2 = \sqrt{\pi} \tau \quad \textup{ and } 
\quad \|u_{\tau}\|_{L^p}^{p} =  \sqrt{\frac{ 2\pi}{p} } \tau.
\eeq
Using  Plancherel's formula we get 
\beq
\label{ExD1-2}
     \begin{aligned}
   \|(|D|^s  - 1) u_\tau\|_{L^2}^2  &  = \frac{1}{2\pi}\|(|\cdot|^s-1)\wh{u}_\tau\|_{L^2}^2  \\
   & =  \tau^2 \int_\R (|\xi|^s-1)^2 e^{- \tau^2(\xi-1)^2 }\,d\xi\\
   & =  \tau^{1-2s} \int_\R (|\tau+x|^s-\tau^s)^2 e^{-x^2} dx  
    \sim \begin{cases}
      \tau^{-1} &\text{as }\tau\to \infty, \\
      \tau^{1-2s} &\text{as }\tau\to 0. \\
    \end{cases}
    \end{aligned}
\eeq
Thus we obtain
  \begin{equation} \label{eq:utau}
    Q_{\ka }(u_\tau)
    \sim \frac{\tau^{\frac{1}{p}}}{\tau^{\frac{\kappa}{2}}  \tau^{\frac{\kappa-1}{2}}}
    = \tau^{\frac 1p + \frac 12 -\kappa}   \lra \infty \quad \text{as }\tau\to\infty 
    \quad \mbox{ if } \frac 1p + \frac 12 -\kappa >0, 
\end{equation}   
while in the limit $ \tau \to 0 $ we have 
\begin{equation}
    Q_{\ka}(u_\tau)
    \sim \frac{\tau^{\frac{1}{p }}}{\tau^{\frac{\kappa}{2}}  \tau^{\frac{(1-2s)(1-\kappa)}{2}}}
    = \tau^{\frac{1}{p } - \frac{\ka }{2 } -  \frac{(1 - 2s) ( 1 - \ka )}{2} } \lra \infty
    \quad \mbox{ if }    \frac{1}{p } - \frac{\ka }{2 } -  \frac{(1 - 2s) ( 1 - \ka )}{2} <0.
  \end{equation}
Notice that for $ p > 2 $ and $ \ka \in (0, 1)$, 
the inequality  $ \frac 1p + \frac 12 -\kappa >0 $ 
is equivalent to $  {\frac{( 1 - \ka )p}{p-2 }} > \frac 12 $ and the inequality 
$  \frac{1}{p } - \frac{\ka }{2 } -  \frac{(1 - 2s) ( 1 - \ka )}{2} <0 $ 
is equivalent to $    {\frac{( 1 - \ka )p}{p-2 }} < \frac{1}{2s}. $
\end{example}

\begin{remark}
\label{not-enough} 
\rm
Let $F$ be as in (\ref{Fdet}). 
Using polar coordinates in $ \R^N$ and proceeding as in the proof of Theorem \ref{T-1D}, 
one can prove that $F$ is bounded near zero if and only if $ \frac{( 1 - \ka )p}{p - 2} \ges \frac{N}{ 2s}$ 
and $F$ is bounded near infinity if and only if $ \frac{( 1 - \ka )p}{p - 2} \les \frac 12$. 
By Lemma \ref{measure}, $M_3$ is finite if and only if $F$ {is bounded, and a necessary} condition for the boundedness of $F$ would be 
$ \frac{N}{ 2s}\les \frac 12$, or equivalently  $ N \les s $. 
Thus any attempt to show that $ M< \infty $ by proving that $ M_3 < \infty $  will fall short from  providing the optimal range of parameters $ (\ka, s)$ for which the supremum in (\ref{supr}) is finite, given by Theorem \ref{T-multiD}.
The Hausdorff-Young inequality (that we have used successfully in dimension one) is not sufficiently accurate  in higher dimensions and a more subtle argument is needed.
See also the second part of Remark \ref{R-generalization}.
\end{remark}

The next theorem gives optimal conditions for the finiteness of $M$ in any space dimension $N \ges 2$. 
Its proof is based on   the Tomas-Stein Theorem, which asserts that for $ p \ges \frac{2N + 2}{N-1}$, there exists a positive constant
$C_{TS}$ depending only on $ p $ and on $N$ such that for any $ \ph \in L^2( \Sp^{N-1} , d \si )$ there holds
\beq
\label{TS}
\| \wh{f \, d \si} \|_{L^p ( \R^N)} \les C_{TS} \| f \|_{L^2( \Sp^{N-1} ) } 
\eeq
(see, e.g., Theorem 7.1 p. 45 in \cite{Wolff}). 
Here $ \Sp^{N-1} = \{ \omega \in \R^N \; | \; | \omega | = 1 \} $ is the unit sphere in $ \R^N$, 
$ \si $ is the usual surface measure on the unit sphere, 
$L^2( \Sp^{N-1} , d \si )$ is the space of measurable functions defined on the unit sphere which are square integrable with respect to the surface measure, and $ \wh{f \, d \si} $ is the Fourier transform of the measure  $ f\,d\sigma$,  given by 
$$
\wh{f \, d \si} ( \xi) = \int_{\Sp^{N-1}} f( \omega ) e^{- i \omega . \xi } \, d \si ( \omega) \qquad \mbox{ for any }  \xi \in \R^N \mbox{ and any }  f \in L^1( \Sp^{N-1}, d \si).
$$

\begin{Theorem}
\label{T-multiD}
Let $ N \in \N$, $ N \ges 2 $,  
$p \in (2, \infty)$, $ \ka \in (0, 1)$,  and $ s > 0 $. 
Then $ Q_{\ka}$ is bounded on $ H^s( \R^N) \setminus \{ 0 \}$ (that is, $M$ in (\ref{supr}) is finite) 
if and only if 
\beq
\label{conditions}
 \ka \ges \frac12  \quad \mbox{  and }  \quad
 \frac{N}{s}\left(\frac 12 - \frac{1}{p} \right)  \les 1 - \ka \les    \frac{N+1}{2}\left(\frac 12 - \frac{1}{p} \right).
\eeq

\end{Theorem}

Observe that the last condition in (\ref{conditions}) implies that $s \ges \frac{2N}{N + 1}$
and it  is equivalent to $   \frac{2(N+1)}{N + 4\ka -3 } \les p \les \frac{ 2N}{N - 2( 1 - \ka) s}$ if $  2( 1 - \ka) s < N$, respectively to
$   \frac{2(N+1)}{N + 4\ka -3 } \les p $ if $  2( 1 - \ka) s \ges N$.

\medskip

{\it Proof. } 
Assume that (\ref{conditions}) hold and 
assume also in a first stage
 that $ p \ges \frac{2(N+1)}{N-1}$, so that we may apply the  Tomas-Stein Theorem.
Let $ u \in \So( \R^N)$. 
Using the Fourier inversion formula and passing to polar coordinates in $ \R^N$ we get
$$
u(x) 
=  \frac{1}{(2 \pi)^N} \int_0 ^{\infty} r^{N-1} \int_{\Sp^{N-1} } \wh{u}( r \omega ) e^{i (rx).\omega } d \si ( \omega) \, dr
=  \frac{1}{(2 \pi)^N} \int_0 ^{\infty} r^{N-1} \wh{\wh{u} ( r \cdot ) \, d \si } ( - r x) \, dr. 
$$
Using Minkowski's inequality in integral form (see, e.g., Theorem 2.4 p. 47 in \cite{lieb-loss}), then the Tomas-Stein inequality we get 
$$
\begin{array}{l}
\ds \| u \|_{L^p ( \R^N )} \les   \frac{1}{(2 \pi)^N}  \int_0^{\infty} r^{N-1} \big\|  \wh{\wh{u} ( r \cdot ) \, d \si } ( - r \cdot )  \big\|_{L^p ( \R^N )} \, dr
\\
\\
=   \ds \frac{1}{(2 \pi)^N}  \int_0^{\infty} r^{N-1 - \frac Np} \big\|  \wh{\wh{u} ( r \cdot ) \, d \si }   \big\|_{L^p ( \R^N )} \, dr
\\
\\
\ds \les \frac{C_{TS}}{(2 \pi)^N}  \int_0^{\infty }  r^{N-1 - \frac Np} \big\| \wh{u} ( r \cdot ) \big\|_{L^2 ( \Sp^{N-1})}                      \, dr .
\end{array}
$$
Denoting $ z_u ( r ) = r^{\frac{N-1}{2}} \left( \int_{\Sp ^{N-1}} | \wh{u}( r \omega ) |^2 \, d \si ( \omega) \right)^{\frac 12}, $
we have proved that there exists $ C > 0 $ depending only on $ p $ and on $ N$ such that 
\beq
\label{Tomas-Stein}
\| u \|_{L^p ( \R^N)} \les C \int_0^{\infty } r ^{ \frac{N-1}{2} - \frac N p } z_u ( r ) \, dr 
\qquad \mbox{ for all }  u \in \So ( \R^N) \mbox{ and for all } p \ges \frac{ 2N + 2}{N-1}. 
\eeq
On the other hand, using  Fourier's inversion formula and polar coordinates in $ \R^N$ we have 
\beq
\label{L2-polar}
\| u \|_{L^2 } ^2= \frac{1}{( 2 \pi)^N} \int_0^{\infty } \int_{\Sp ^{N-1}} | \wh{u} ( r \omega ) |^2 \, d \si ( \omega ) \, dr =  \frac{1}{( 2 \pi)^N} \int_0 ^{\infty } z_u ^2 ( r ) \, dr
\eeq
and 
\beq
\label{OpDs}
\begin{array}{l}
\ds \| ( | D |^s - 1 ) u \|_{L^2}^2 =  \frac{1}{( 2 \pi)^N}  \big\| ( |\cdot |^s - 1 ) \wh{u} \big\|_{L^2}^2 
\\
\\
\ds =  \frac{1}{( 2 \pi)^N}  \!  \int_0^{\infty } r^{N-1} ( r^s - 1 )^2  \int_{\Sp^{N-1}}  | \wh{u} ( r \omega ) |^2 \, d \si ( \omega ) \, dr
=  \frac{1}{( 2 \pi)^N}  \!  \int_0^{\infty } ( r^s - 1 )^2 z_u ^2 ( r ) \, dr.
\end{array}
\eeq
From (\ref{Tomas-Stein}) - (\ref{OpDs}) it follows that there is $C > 0 $ such that for all $u \in \So ( \R^N)$ we have 
\beq
\label{Qk-bound}
Q_{\ka }(u) \les C
\frac{ \int_0^{\infty } r ^{ \frac{N-1}{2} - \frac N p } z_u ( r ) \, dr  }{ \| z_u \| _{L^2(0, \infty) }^{\ka } 
\left( \int_0^{\infty } ( r^s - 1 )^2 z_u ^2 ( r ) \, dr \right)^{\frac{ 1 - \ka}{2}} }.
\eeq
Notice that $ z_u   \in L^2(0, \infty) $ and $ \left( | \cdot |^s - 1 \right) z_u   \in L^2(0, \infty)$ by (\ref{L2-polar}) and (\ref{OpDs}).
We use Lemma \ref{measure} in $ (0, \infty) $ endowed with the usual Lebesgue measure and we take
 $ w( r ) =  r ^{ \frac{N-1}{2} - \frac N p } $, $ w_1 ( r ) = 1 $ and $ w_2 ( r ) = r^s -1 $. We get
\beq
\label{ineq-sup}
\sup_{\ph \in L^2(0, \infty) \setminus \{ 0 \}   }
 \frac{\|  w \ph \|_{L^1(0, \infty)}}{\| \ph \|_{L^2(0, \infty)}^{\ka } \| w_2 \ph \|_{L^2(0, \infty)}^{ 1-\ka }  }
\les  C 
 \sup_{t > 0 } \left( t ^{\frac{ 1 - \ka}{2}} 
\Big\| \frac{ w}{\left(1 + t | w_2 |^2 \right)^{\frac 12} } \Big\|_{L^2(0, \infty)} \right).
\eeq
Let 
\beq
\label{Gdet}
G(t ) := \Big\| \frac{  t ^{\frac{ 1 - \ka}{2}}  w}{\left(1 + t | w_2 |^2 \right)^{\frac 12} } \Big\|_{L^2(0, \infty)} ^2
= \int_0^{\infty} \frac{ r^{N-1 - \frac{2N}{p}} }{ t^{\ka -1 } + t^{\ka } (r^s -1 )^2 } \, dr
=: \int_0^{\infty} g(r,t) \, dr.
\eeq
Since  $ N - 2s - \frac{2N}{p} < 0$ (because $ \frac{N}{s}\left(\frac 12 - \frac{1}{p} \right)  \les 1 - \ka< 1$), we have $ g ( \cdot, t) \in L^1(0, \infty)$ and then 
Remark \ref{R1.2} (iii) implies that $G$ is continuous on $(0, \infty)$. 

For any fixed $ A > 0 $ we have 
$$
0 < \int_0^A g( r, t ) \, dr < t^{ 1 - \ka } \int_0^A  r^{N-1 - \frac{2N}{p}} \, dr
=  t^{ 1 - \ka } \frac{A^{N - \frac{2N}{p}} }{N - \frac{2N}{p}} \lra 0 \qquad \mbox{ as } t \lra 0 .
$$
We have $ r^s - 1 > \frac 12 r^s $ if $ r > 2^{\frac 1s}$. 
Taking $ A \ges  2^{\frac 1s}$ and using the change of variable $ r = t^{- \frac{1}{2s}} y $ we find
$$
0 < \int_A^{\infty} g(r, t ) \, dr 
< \int_A^{\infty}  \frac{ r^{N-1 - \frac{2N}{p}} }{ t^{\ka -1 } + \frac 14 t^{\ka } r^{2s} } \, dr 
= t^{ 1 - \ka - \frac{1}{2s}(N - \frac{2N}{p})} \int_{t^{\frac{1}{2s}} A }^{\infty} 
\frac{y^{N-1 - \frac{2N}{p}} }{1 + \frac 14 y^{2s}} \, dy.
$$
We infer that $G$ is bounded 
as $ t \lra 0 $ if $ 1 - \ka \ges \frac{N}{s}\left(\frac 12 - \frac{1}{p}\right)$.

It is clear that 
$$
0 < \int_2^{\infty } g( r, t ) \, dr < t^{-\ka }  \int_2^{\infty } \frac{ r^{N-1 - \frac{2N}{p}} }{ (r^s -1 )^2 }\, dr \lra 0 \qquad \mbox{ as } t \lra \infty.
$$
There is $ C > 0 $ such that $  r^{N-1 - \frac{2N}{p}} \les C$ for $ r \in [0, 2 ]$ because $ N-1 - \frac{2N}{p} > 0 $ (recall that $p \ges \frac{2N + 2}{N-1}$).
There is  $ c_1 > 0 $ such that $ \left( | 1 + y |^s -1 \right)^2 \ges c_1 y^2 $ for all $ y \in [-1, 1]$.
Using the change of variable $ r = 1 + y$, the previous observations, then the change of variable $ z =  t^{\frac 12} y $ we get 
$$
0 < \! \int_0^2 g(r, t) \, dr \les \! \int_{-1}^1 \frac{C}{t^{\ka - 1} + t^{\ka } ( | 1 + y |^s -1)^2 } \, dy 
\les  \! \int_{-1}^1 \frac{C}{t^{\ka - 1} + t^{\ka } c_1^2 y^2 } \, dy
=  \! \int_{-\sqrt{t}}^{\sqrt{t}} \frac{  C t ^{ \frac 12 - \ka } }{1 + c_1^2 z^2 }\, dz. 
$$
We conclude that $ G $ is bounded as $ t \lra \infty $ if $ \ka \ges \frac 12$.

We have thus proved that 
if (\ref{conditions}) and the additional assumption $ p \ges \frac{2(N+1)}{N-1}$ hold, 
 the function 
$G$ is bounded on $(0, \infty)$ and 
therefore the supremum on the right hand side of (\ref{ineq-sup}) is finite. 
Then (\ref{Qk-bound}) and (\ref{ineq-sup}) imply that there exists $ C(\ka )> 0 $ such that 
$ Q_{\ka } ( u ) \les C(\ka )$ for any $ u \in \So ( \R^N) \setminus \{ 0 \}$. 
Since $ u \longmapsto  Q_{\ka } ( u ) $ is continuous on $ H^s( \R^N) \setminus \{ 0 \}$ and 
$  \So ( \R^N) \setminus \{ 0 \}$ is dense in $ H^s( \R^N)$, 
we infer that $ Q_{\ka }(u) \les C(\ka )$ for any  $ u \in H^s( \R^N) \setminus \{ 0 \}$.

\medskip

It remains to consider the case $    \frac{2(N+1)}{N + 4\ka -3 } \les p < \frac{2(N+1)}{N-1}$.
We proceed by interpolation.
Denote $ q : = \frac{2(N+1)}{N-1}$. 
We see that $ \frac{N+1}{2}\left(\frac 12 - \frac{1}{q} \right) = \frac 12$.
Since $ 2 < p < q$, there is some $ \theta \in (0, 1 )$ such that 
$ \frac 1p = \frac{\theta}{2} + \frac{1 - \theta}{q}$.  We have
$$ 
\theta = \left( \frac 1p - \frac 1q \right) \left( \frac 12 - \frac 1q \right)^{-1}
= (N+1) \left( \frac 1p - \frac 1q \right) 
\les  (N+1) \left( \frac{N + 4\ka -3 }{2(N+1)} - \frac{N-1}{2(N+1)} \right) = 2 \ka -1. 
$$
Since $ \ka <1$, the above inequality gives $ \theta < \ka $. By H\"older's inequality we have 
$$
\| u \|_{L^p} \les \| u \|_{L^2}^{\theta}\cdot  \| u \|_{L^q}^{1 - \theta}
$$
and then for any $ u \in H^s( \R^N) \setminus \{ 0 \}$ we find
\beq
\label{sub-TS}
Q_{\ka }(u ) \les \frac{\| u \|_{L^2}^{\theta}\cdot  \| u \|_{L^q}^{1 - \theta}}{\| u \|_{L^2}^{ \ka } 
\| (|D|^s - 1) u \|_{L^2}^{1 - \ka } }
= \left( \tilde{Q}_{\tilde{\ka }}( u ) \right)^{ 1 - \theta}, 
\eeq
where  $ \tilde{\ka} = \frac{ \ka - \theta}{1 - \theta}$, 
so that $ 1 -  \tilde{\ka} = \frac{ 1 - \ka }{1 - \theta}$, and 
$ \tilde{Q}_{\tilde{\ka }}( u )  = \frac{ \| u \|_{L^q} }{\| u \|_{L^2}^{ \tilde{\ka }} 
\| (|D|^s - 1) u \|_{L^2}^{1 - \tilde{\ka }} }.$
Notice that the inequality $ \theta < 2 \ka -1 $ implies that $ \tilde{\ka } \ges \frac 12$ and then we get 
$ 1 - \tilde{\ka} \les \frac 12 = \frac{N+1}{2}\left(\frac 12 - \frac{1}{q} \right) $.
Using (\ref{conditions})  and  the fact that
$ \frac 12 - \frac 1p = ( 1 - \theta ) \left(\frac 12 - \frac 1q \right)$, we get
$$
 1 -  \tilde{\ka} = \frac{ 1 - \ka }{1 - \theta} \ges (1 - \theta)^{-1} \frac Ns \left( \frac 12 - \frac 1p \right) = \frac Ns \left( \frac 12 - \frac 1q \right).
 $$
Thus we see that (\ref{conditions}) is satisfied with $q$ and $ \tilde{\ka }$ 
  instead of $ p $ and $ \ka$, respectively. 
From the first part of the proof we infer that    $ \tilde{Q}_{\tilde{\ka }} $ is bounded from above on 
$H^s( \R^N) \setminus \{ 0 \}$, and then (\ref{sub-TS}) implies that $ Q_{\ka } $ is also bounded.

\medskip

So far we have proved that $Q_{\ka}$ is bounded on $ H^s( \R^N) \setminus \{ 0 \}$ if (\ref{conditions}) hold. 
Now let us show that (\ref{conditions}) is necessary for the boundedness of $Q_{\ka}$. 
Let $ u \in \So ( \R^N)$, $ u \neq 0 $. 
For $ \tau > 0 $ let  $ u _{\tau }(x) = u \left( \frac{x}{\tau} \right)$. A simple 
 computation gives $ \| u _\tau \|_{L^q } = \tau^{\frac{N}{q}} \| u \|_{L^q} $ 
 for any $ q \in [1, \infty )$ and
$$
\| ( |D|^s - 1 ) u \|_{L^2}^2 = \frac{\tau^{N - 2s}}{(2 \pi)^N} \int_{\R^N} \left( |\xi |^{2 s} - 2 \tau ^s |\xi |^s + \tau^{ 2 s} \right) |\wh{u}( \xi ) |^2 \, d \xi .
$$
Thus we find 
$$
Q_{\ka } ( u_{\tau} ) = ( 2 \pi)^{\frac{(1 - \ka)N}{2}} \tau^{\frac Np - \frac{\ka N}{2} -  (1 - \ka)(\frac N2 -s ) } \frac{ \| u \|_{L^p}}{\| u \|_{L^2}^{\ka } 
\left( \int_{\R^N} \left( |\xi |^{2 s} - 2 \tau ^s |\xi |^s + \tau^{ 2 s} \right) |\wh{u}( \xi ) |^2 \, d \xi \right)^{\frac{ 1 - \ka}{2} }}.
$$
If $ Q_{\ka }( u _{\tau}) $ remains bounded as $ \tau \lra 0 $ we must have
$ \frac Np - \frac{\ka N}{2} -  (1 - \ka)(\frac N2 -s )  \ges 0 $ 
and this is equivalent to $ 1 - \ka \ges \frac{N}{s}\left(\frac 12 - \frac{2}{p} \right)$.

The next example shows that 
 $ Q_{\ka }$ is not bounded if 
$ \ka < \frac 12$ or if  $ \frac{N+1}{2} \left( \frac 12 - \frac 1p \right) < 1 - \ka $.
\hfill
$\Box$

\begin{example}
\label{knapp}
\rm
We consider a variant of Knapp's example related to the Tomas-Stein inequality (see, e.g., \cite{Wolff} p. 46).
For small $ \de > 0 $, let $ S_{\de } = \{ \omega = ( \omega_1, \dots, \omega _N) \in \Sp^{N-1} \; | \; \omega _N > 1 - \de ^2 \}$. 
It is easily seen that 
there exist positive constants $ C_1, \, C_2$ such that  for all $ \de \in(0,  \frac{1}{10})$, say, we have
$$
 C_1 \de ^{N-1} \les \si( S_{\de}) \les  C_2 \de ^{N-1},
$$
where $ \si $ is the surface measure on $ \Sp^{N-1}$.
For small $ \e >0$, $ \de > 0 $ we define $ v_{\e, \de } : \R^N \lra \R$ by 
$$
v_{\e, \de } ( \xi ) = \left\{ 
\begin{array}{l} 1 \; \mbox{ if } 1 - \e < | \xi | < 1 + \e \mbox{ and } \frac{ \xi }{| \xi | } \in S_{\de}, 
\\ 0 \; \mbox{ otherwise.}
\end{array}
\right.
$$
Let 
\beq
\label{u-epsilon}
 u_{\e, \de } = \Fo ^{-1} ( v_{\e, \de}), 
 \qquad \mbox{ that is }
 u_{\e, \de } (x) = \frac{1}{(2 \pi)^N} \int_{1 - \e}^{ 1 + \e} r^{N-1} 
 \int_{S_{\de} }e^{ i x.(r\omega)} \, d \si ( \omega ) \, dr.
\eeq
Since $ \wh{u_{\e, \de }} = v_{\e, \de }$ is bounded and compactly supported, we have $ u_{\e, \de } \in H^s( \R^N)$ for all $s$. 
By Plancherel's identity we get 
\beq
\label{ex-1}
\| u_{\e, \de } \|_{L^2}^2 = \frac{1}{(2 \pi)^N } \int_{\R^N} |v_{\e, \de }( \xi )|^2 \, d \xi 
=  \frac{1}{(2 \pi)^N } \int_{1-\e }^{1 + \e} r^{N-1} \si ( S_{\de } ) \, dr
\sim \e \de^{N-1} 
\eeq
and
\beq
\label{ex-2}
\begin{array}{l}
\ds \|\left( |D|^s - 1 \right) u_{\e, \de } \|_{L^2}^2 
= \frac{1}{(2 \pi)^N } \int_{\R^N} (|\xi |^s - 1 )^2 |v_{\e, \de }( \xi )|^2 \, d \xi 
\\
\\
\ds =  \frac{1}{(2 \pi)^N } \int_{1-\e }^{1 + \e} r^{N-1} ( r^s - 1 )^2 \si ( S_{\de } ) \, dr
\sim \e ^3 \de^{N-1}.
\end{array}
\eeq
Let $ e_N = ( 0, \dots, 0, 1) \in \R^N$. It is clear that 
$$
 | u_{\e, \de} (x) | = |e^{-i x.e_N} u_{\e, \de} (x) | 
=  \frac{1}{(2 \pi)^N }  \Big|  \int_{1-\e }^{1 + \e} r^{N-1} 
 \int_{S_{\de} }e^{ i x.(r\omega - e_N)} \, d \si ( \omega ) \, dr \Big|.
$$
Let $ A_{\e, \de } = \left\{ x \in \R^N \;\big| \; |  x.(r\omega - e_N) | \les \frac{\pi}{3} \mbox{ for all } 
r \in (1- \e, 1 + \e) \mbox{ and all } \omega \in S_{\de} \right\}$. 
For any $ x \in A_{\e, \de}$, $r \in (1-\e, 1 + \e)$, and $ \omega \in S_{\de}$ we have 
$\text{Re} ( e^{ i x.(r\omega - e_N)} ) \ges \frac 12$, hence
\begin{align*}
 | u_{\e, \de} (x) | & \ges  \frac{1}{(2 \pi)^N }  \int_{1-\e }^{1 + \e} r^{N-1} 
 \int_{S_{\de} } \text{Re}  \left( e^{ i x.(r\omega - e_N)} \right)\, d \si ( \omega ) \, dr  \\
& \ges C  \int_{1-\e }^{1 + \e} r^{N-1} \si( S_{\de}) \, dr \ges C \e \de^{N-1} 
\end{align*}
for some $ C > 0 $ independent of $ \e $ and $ \de$. We infer that 
\beq
\label{ex-3}
\| u_{\e, \de} \|_{L^p ( \R^N ) } \ges \| u_{\e, \de} \|_{L^p ( A_{\e, \de} ) }
\ges C \e \de^{N-1}  |  A_{\e, \de} |^{\frac 1p}.
\eeq
We will find a lower bound for $|  A_{\e, \de} |$. 
Denote $ x = ( x', x_N)$, $ \omega = ( \omega ', \omega _N)$ where $ x', \omega ' \in \R^{N-1}$,
and assume that $ \e \les 1$.  We have
$$
 |  x.(r\omega - e_N) | = | x.( r \omega', r \omega_N - 1) | \les r | x'.\omega '| + |x_N (r \omega _N -1 )|.
$$
For $ \omega = ( \omega', \omega _N) \in S_{\de }$ and $ r \in [1-\e, 1 + \e ]$ 
we have $ |\omega '| < \sqrt{2} \de $ 
and  $|r \omega _N -1 | \les 2 \de ^2 + \e $, hence   
$ r | x'.\omega '| \les \frac{\pi}{6} $ if  $ | x'| \les \frac{\sqrt{2} \pi}{24 }\de^{-1}$
and 
$|x_N( r \omega _N -1) | \les \frac{\pi}{6} $   if $ | x_N | \les \frac{ \pi}{6}(2\de^2 + \e)^{-1}$.
We conclude that
$$
\left\{ ( x', x_N) \in \R^{N-1} \times \R \; \big| \; |x'| \les \frac{\sqrt{2} \pi}{24 }\de^{-1}, \; | x_N | \les \frac{ \pi}{6}(2\de^2 + \e)^{-1} \right\} 
\subset A_{\e, \de}.
$$
Hence there exists $C > 0 $ independent of $ \e $ and $ \de$ such that $| A_{\e, \de }| \ges  \frac{C}{ \de^{N-1} ( \de^2 + \e )}$ and  (\ref{ex-3}) gives
\beq
\label{ex-4}
\| u_{\e, \de} \|_{L^p ( \R^N ) } \ges C \e \de^{(N-1)( 1 - \frac 1p)} (\de^2 + \e )^{- \frac 1p}.
\eeq
From (\ref{ex-1}), (\ref{ex-2}), and (\ref{ex-4}) we obtain 
\beq
\label{ex-5}
Q_{\ka } ( u _{\e, \de} ) \ges C \de^{(N-1)(\frac 12 - \frac 1p )} \e^{\ka - \frac 12} ( \de^2 + \e )^{- \frac 1p }.
\eeq

Fix $ \de _0 \in (0, \frac{1}{10})$ and let $ \e \lra 0 $. 
If $ Q_{\ka }( u_{\e, \de_0})$ remains bounded, (\ref{ex-5}) implies that $ \ka \ges \frac 12$. 

\medskip

Putting $ \e = \de ^2 $ in (\ref{ex-5}) we get 
$ Q_{\ka }(u_{\de^2, \de} ) \ges C \de ^{(N-1)(\frac 12 - \frac 1p ) + 2 \ka - 1 - \frac 2p}. $
If $Q_{\ka }(u_{\de^2, \de} ) $ remains bounded as $ \de \lra 0 $ we must have 
$(N-1)(\frac 12 - \frac 1p ) + 2 \ka - 1 - \frac 2p \ges 0 $, 
and this is equivalent to $ \frac{N+1}{2} \left( \frac 12 - \frac 1p \right) \ges 1 - \ka. $

\medskip

Obviously, if $ Q_{\ka}$ is bounded on $ H^s( \R^N, \C) \setminus \{  0 \} $ then it is also bounded on $ H^s ( \R^N, \R) \setminus \{  0 \}$. 
The converse is also true. 
Indeed, it is easily seen that $ | D|^s u $ is a real-valued function whenever $ u \in H^s ( \R^N)$ is real-valued. 
If $ u \in H^s ( \R^N)$ is such that $ \Re ( u ) \not\equiv 0 $ and $ \Im ( u ) \not\equiv 0$, 
where $ \Re ( u ) $ and $ \Im ( u ) $ are its real and imaginary parts, respectively, we have
$$
 Q_{\ka} (u ) \les 
 \frac{\| \Re (u) \|_{L^{p}} + \| \Im (u) \|_{L^{p}}}{  \| u \|_{L^2}^{\ka}  \| \left( |D|^s - 1\right) u  \|_{L^2} ^{1 - \ka} }
 \les Q_{\ka} ( \Re ( u ) ) +  Q_{\ka} ( \Im ( u ) ).
$$

\end{example}

\begin{remark} 
\label{R-generalization}
\rm 
The method used in the proof of Theorem \ref{T-multiD} is very flexible and can be used  to prove non-homogeneous Gagliardo-Nirenberg inequalities of the form 
\beq
\label{general}
\| u \|_{L^p} \les C \| P_1(D) u \|_{L^2} ^{\ka } \| P_2(D) u \|_{L^2} ^{1 - \ka } , 
\eeq
where $ N \ges 2 $,  $ p \ges \frac{2N + 2}{N-1}$ and $ P_1(D)$,  $P_2(D)$ are Fourier integral operators defined by 
$$
P_i (D) (u ) = \Fo ^{-1} \left( P_i ( \cdot ) \wh{u} \right), \quad i = 1, 2.
$$
Assuming that there exist non-negative functions $ p_1, p_2 : [0, \infty ) \lra \R_+$ such that 
$ | P_i ( \xi ) | \ges p_i ( | \xi | )$ for all $ \xi \in \R^N$, $ i = 1,\, 2 $
and proceeding as in (\ref{L2-polar}) and (\ref{OpDs}) we get for all $ u \in \So ( \R^N)$
\beq
\label{ineq-gen}
\| P_i (D) u \|_{L^2}^2 
\ges  \frac{1}{( 2 \pi)^N} \int_0^{\infty } p_i^2(r)  z_u ^2 ( r ) \, dr, 
\quad 
\mbox{ where } z_u ( r ) = r^{\frac{N-1}{2}} \| \wh{u}( r \cdot )\|_{L^2( \Sp^{N-1})}.
\eeq
In order to prove the inequality (\ref{general}) in some function space $ \Xo $ (typically $ \Xo = H^s( \R^N)$, but other spaces might be considered), one needs to show the continuity 
of the $ L^p- $norm and of the operators $ P_1(D)$ and $ P_2(D)$ on $ \Xo $, as well as the density of 
$ \So ( \R^N)$ in $ \Xo $. 
Then, taking into account (\ref{Tomas-Stein}) and (\ref{ineq-gen}) and denoting 
$ w(r ) =  r ^{ \frac{N-1}{2} - \frac N p } $, it suffices  to show that 
$$
\sup  
\left\{  \frac{\|  w \ph \|_{L^1(0, \infty)}}{\| p_1 \ph \|_{L^2(0, \infty)}^{\ka } \| p_2 \ph \|_{L^2(0, \infty)}^{ 1-\ka }  } \; \Big| \; p_i \ph \in L^2(0, \infty) \setminus \{ 0 \} , \, i = 1,2 \right\}
$$
is finite. 
To do this, by Lemma \ref{measure} and Remark \ref{R1.2} (iii) it suffices to prove that the function 
\beq
\label{Hdet}
H(t ) := \Big\| \frac{  t ^{\frac{ 1 - \ka}{2}}  w}{\left(p_1^2 + t | p_2 |^2 \right)^{\frac 12} } \Big\|_{L^2(0, \infty)} ^2
= \int_0^{\infty} \frac{ r^{N-1 - \frac{2N}{p}} }{ t^{\ka -1 } p_1^2 (r) + t^{\ka } p_2^2(r) } \, dr
\eeq
is bounded on $(0, \infty)$.

If $ N = 1$ or if the Tomas-Stein inequality is not available (for instance, if $2 < p < \frac{2N + 2}{N-1}$), one may try to use the Hausdorff-Young Theorem  to prove the inequality (\ref{general}), as in the proof of Theorem \ref{T-1D}.
Indeed, to establish (\ref{general}) it suffices to show that the supremum 
$$
\sup \left\{  \frac{ \| \ph \|_{L^{ p'}}}{ \| P_1\ph \|_{L^2}^{ \ka } 
\big\| P_2 \ph \big\|_{L^2}} \; \Big| \; 
P_i \ph \in L^2( \R^N) \setminus \{ 0 \} , \;  i = 1,\, 2 \right \}
$$
is finite, and by Lemma \ref{measure} this amounts to proving that the function
\beq
\label{Kdet}
K(t) := \Bigg\| \frac{  t ^{\frac{ 1 - \ka}{2}}  }{\left(P_1^2 + t P_2^2 \right)^{\frac 12}} \Bigg\|_{L^{ \frac{2p}{p-2 }}} ^{ \frac{2p}{p-2 }}
= \int_{ \R^N} \frac{1}{\left( t^{\ka - 1} | P_1( \xi )|^2+ t^{\ka }|P_2(\xi)|^2 \right)^{\frac{p}{p-2 }} }\, d \xi 
\eeq
is bounded on $(0, \infty)$. 
However, we expect the approach based on the Hausdorff-Young inequality to give weaker results than the approach based on the Tomas-Stein inequality. 
For instance, if $P_1$ and $P_2$ are radial and non-negative 
 (that is, if $ P_i( \xi ) = p_i ( |\xi | ) \ges 0 $), 
it is easily seen that the boundedness of the function $K$ implies the boundedness of the function  $H$, but the converse might not be true. See Remark \ref{not-enough}.

As a matter of fact, the Gagliardo-Nirenberg-Sobolev inequality 
$$
\| u \|_{L^p } \les C \| u \|_{L^2}^{1 - \frac{N( p-2)}{4 sp} } \|( - \Delta )^s u \|_{L^2}^{ \frac{N( p-2)}{4 sp} } 
$$
(with  $ 2 < p < \frac{2N}{N - 4s} $ if $ s < \frac N4$) can be proven by using our method and the 
Hausdorff-Young inequality;  in this case $ P_1 = 1$ and $ P_2( \xi ) = |\xi |^{ 2s}$, 
and the integral in $K(t)$ is easily evaluated  using polar coordinates and the change of variables 
$ r = t^{- \frac{1}{4s}}y$.
For $ s = 1 $  and $ p = 2 \si + 2$, this gives (\ref{GNS}).
\end{remark}

\begin{remark} 
\label{qualitative-generalization}
\rm
More quantitative  variants of (\ref{general}) can be proved, too. 
For instance, in some applications it is useful to dispose of inequalities of the form (\ref{general}) under  the additional constraint $ \| P_2( D ) u \|_{L^2} \les R  \| P_1( D ) u \|_{L^2}$, where $ R>  0 $ is given.
To obtain such inequalities we may use a slight modification of Lemma \ref{measure}. 

With the notation and the assumptions in Lemma \ref{measure}, let 
$$
M_1^R = \sup \left\{ 
\frac{\| \ph w \|_{L^q}}{ \| \ph w_1\|_{L^2}^{\ka } \cdot \| \ph w_2 \|_{L^2}^{ 1 - \ka }} \; \Bigg| \;
\begin{array}{l} \ph \mbox{ is measurable, }
\ph w_i \in L^2 (X)  \mbox{ for } i =1, \, 2 , 
\\
\mbox{and } 0<\| \ph w_2 \|_{L^2} \les R \| \ph w_1 \|_{L^2}  
\end{array}
\right\}, 
$$
$$
M_2 ^a= \ds \sup_{t > a } \left( t ^{\frac{ 1 - \ka}{2}} 
\Big\| \frac{ w}{\left(w_1 ^2 + t | w_2 |^2 \right)^{\frac 12}} \Big\|_{L^{ \frac{2q}{2 -q }}} \right).
$$
Then we have $ M_1^R \les    ( 1 - \ka )^{ \frac{ \ka-1 }{2} } \ka ^{- \frac{\ka}{2} } M_2^{\frac{( 1 - \ka)  }{\ka R^2}} $.

To prove the above statement  we use again the observation that 
for any  $A, B > 0 $, the function $ f(t ) = At^{\ka - 1} + B t^{\ka } $ achieves its minimum on $(0, \infty )$ at $ t_{min} = \frac{( 1 - \ka)A}{\ka B} $ and $ f( t_{min} ) = ( 1 - \ka ) ^{\ka - 1} \ka^{-\ka } A^{\ka } B ^{1 - \ka }. $
If  $ \ph $ is a measurable function satisfying 
 $ \ph w_i \in L^2(X)$ for $ i = 1, \, 2 $  and $ 0 < \| \ph w_2 \|_{L^2} \les R \| \ph w _1 \|_{L^2}  $, 
taking  $ A = \| \ph w_1 \|_{L^2}^2 $,  $ B = \| \ph w_2 \|_{L^2}^2$
and $ t_{min} = \frac{( 1 - \ka) \| \ph w_1 \|_{L^2}^2 }{\ka  \| \ph w_2 \|_{L^2}^2} $, 
we see that 
 (\ref{1.1}) holds and, moreover, $ t_{min} \ges \frac{( 1 - \ka)  }{\ka R^2}$.
Then we use (\ref{1.1}) and (\ref{1.3}) and we proceed exactly as in the proof of Lemma \ref{measure}.

\medskip

Assume that $ P_1$ and $P_2$ are radial, that is $ P_i ( \xi ) = p_i (|\xi|) $ for $i = 1, 2$. 
Then we have  equality in (\ref{ineq-gen}) and the condition 
$ \| P_2( D ) u \|_{L^2} \les R  \| P_1( D ) u \|_{L^2}$ 
is equivalent to $ \| p_2 (|\cdot |)  \wh{u} \|_{L^2(\R^N)} \les R  \| p_1  (|\cdot |)\wh{u} \|_{L^2(\R^N)}$ 
and to 
$ \| p_2 z_ u \|_{L^2(0, \infty)} \les R  \| p_1 z_ u \|_{L^2(0, \infty)}$.
We infer that 
\beq
\label{qualitat-sup}
\sup \left\{ \frac{ \| u \|_{L^p}}{ \| P_1(D) u \|_{L^2} ^{\ka } \| P_2(D) u \|_{L^2} ^{1 - \ka } }
\; \Big| \; u \in \So ( \R^N), 0 < \| P_2(D) u \|_{L^2} \les R \| P_2(D) u \|_{L^2} \right\}
\eeq
is finite if one of the functions $H$ or $K$ defined in (\ref{Hdet}) and  in (\ref{Kdet}) is bounded on 
$ [  \frac{( 1 - \ka)  }{\ka R^2}, \infty )$. 
If $H(t)$ (respectively $K(t)$) is finite for some $ t > 0 $, it suffices to verify the boundedness of $H$ 
(respectively of $K$) in a neighbourhood of infinity. 
Of course, having explicit bounds on $H$ or on $K$ on the interval $ [  \frac{( 1 - \ka)  }{\ka R^2}, \infty )$
would provide explicit bounds on the supremum in (\ref{qualitat-sup}).

\end{remark}

Remark  \ref{qualitative-generalization}  enables us to state the following quantitative variant of Theorems \ref{T-1D} and \ref{T-multiD}.

\begin{Corollary}
\label{C-quant}
Let $ Q_{\ka}$ be as in (\ref{Qka}). The supremum    
\beq
\label{sup-qualitat}
\sup \left\{ Q_{\ka }( u ) \; \big| \; u \in H^s( \R^N) \setminus \{ 0 \} \mbox{ and } 
\| ( |D|^s - 1 ) u \|_{L^2} \les R \| u \|_{L^2} \right\} 
\eeq
 is finite for any fixed $R>0$ if
\beq
\label{conditions-qualitat}
 \frac 1p > \frac 12 - \frac sN, \qquad  \ka \ges \frac 12  \qquad \mbox{ and } \qquad 
  1 - \ka \les    \frac{N+1}{2}\left(\frac 12 - \frac{1}{p} \right).
\eeq
\end{Corollary}

{\it Proof. } 
We have already seen in the proofs of Theorems \ref{T-1D} and \ref{T-multiD} that the functions $F$ given by (\ref{Fdet}) and $G$ given by (\ref{Gdet}) 
are well-defined and continuous on $(0, \infty)$ if $ \frac 1p > \frac 12 - \frac sN$.

In the proof of Theorem \ref{T-1D} it is shown  that $F$ is bounded in a neighbourhood of infinity if 
$ 1 - \ka \les \frac 12 - \frac 1p$, and Remark \ref{qualitative-generalization} above  implies Corollary \ref{C-quant}  in dimension $N=1$.

Assume that $N \ges 2$. In the case $ p \ges \frac{2(N+1)}{N -1}$, it is shown in the proof of Theorem \ref{T-multiD} that the function $G$ is bounded near infinity if $ \ka \ges \frac 12$ and this proves Corollary \ref{C-quant}. 
In the case $    \frac{2(N+1)}{N + 4\ka -3 } \les p < \frac{2(N+1)}{N-1}$, 
the conclusion  follows from the case $ p = \frac{2(N+1)}{N -1}$ by interpolation, using
  (\ref{sub-TS}). 
\hfill
$\Box$

\section{Global minimisation of the energy at fixed $L^2-$norm}
\label{Globalmin}

In this section we study the minimisation problem $(\Po_m)$. Recall that $E_{min}$ has been introduced in (\ref{EMIN}).
Scaling properties of various terms appearing in $E$ will be important. 
It is easily seen that  for any function $ u \in H^2 ( \R^N)$ and for any $ a, b > 0 $, letting $ u_{a, b }( x) = a u \left( \frac xb \right)$
we have 
\beq
\label{scaling}
\begin{array}{c}
\ds \int_{\R^N} | \Delta u_{a, b }|^2 \, dx  = a^2 b^{ N-4 } \! \! \int_{\R^N} | \Delta u|^2 \, dx,  \; \; 
\int_{\R^N} | \nabla u_{a, b }|^2 \, dx  = a^2 b^{ N-2 } \! \! \int_{\R^N} | \nabla u|^2 \, dx,  
\\ \\
\ds \int_{\R^N} |  u_{a, b }|^{2 \sigma + 2} \, dx  = a^{2\sigma + 2} b^{ N} \! \! \int_{\R^N} |  u|^2 \, dx,  \qquad
\int_{\R^N} |  u_{a, b }|^2 \, dx  = a^2 b^{ N } \int_{\R^N} |  u|^2 \, dx.
\end{array}
\eeq

\medskip

Using the Plancherel Theorem we have for all $ u \in H^2 ( \R^N)$
$$
 \| \Delta u \|_{L^2 }^2 = \frac{ 1}{ ( 2 \pi )^N} \| | \xi |^2 \wh{u} \|_{L^2}^2 = \frac{ 1}{ ( 2 \pi )^N} \int_{\R^N} | \xi  |^4 | \wh{u} ( \xi )|^2 \, d \xi 
\; \; 
\mbox{ and }
\; \; 
\Big\| \frac{ \p^2 u}{\p x_j \p x_k } \Big\|_{L^2 }^2 = \frac{ 1}{ ( 2 \pi )^N} \Big\|  \xi _j \xi _k  \wh{u} \|_{L^2}^2. 
$$
It is then obvious that $ \big\| \frac{ \p^2 u}{\p x_j \p x_k } \big\|_{L^2 } \les \| \Delta u \|_{L^2 }$.
We have also the interpolation inequality 
\beq
\label{interpol}
\| \nabla u \|_{L^2}^2 = \frac{1}{(2 \pi )^N } \| \, | \cdot | \wh{u} \|_{L^2}^2 
\les  \frac{1}{(2 \pi )^N } \| \, | \cdot | ^2 \wh{u} \|_{L^2} \|  \wh{u} \|_{L^2}
= \| \Delta u \|_{L^2} \| u \|_{L^2}. 
\eeq
Notice that we have {\it strict} inequality in (\ref{interpol}), except for $ u = 0 $.

We denote $ 2^{**} = \infty $ if $ N \les 4$ and $ 2^{**} = \frac{ 2N}{N-4}$    if $ N \ges 5$. 
It is well-known (see, e.g., \cite{brezis} section 9.3) that $ H^2 ( \R^N) \subset  L^{\infty }( \R^N)$ if $ N
\les 3 $, $ H^2 ( \R^N) \subset  L^{p }( \R^N)$ for any $ p \in [2, \infty )$ if $ N = 4 $ and  $ H^2 ( \R^N) \subset L^{2^{**} }( \R^N)$ if $ N \ges 5$. 
Moreover, in the latter case we have the Sobolev inequality $ \| u \|_{L^{ 2^{**}} } \les C_S \| \Delta u \|_{L^2}$ for any $ u \in H^2 ( \R^N)$. 

For any $ \sigma \in [0, \frac{2^{**}}{2} - 1) $ we have the Gagliardo-Nirenberg-Sobolev inequality 
\beq
\label{GNS}
\| u \|_{ L^{ 2 \si + 2 } }^{ 2 \si + 2 } \les B \| \Delta u \|_{L^2 }^{ \frac{ \si N}{2}}  \|  u \|_{L^2 }^{  2 + 2 \si -\frac{\si N}{2}} 
\qquad \mbox{ for all } u \in H^2 ( \R^N), 
\eeq
where $B$ is independent of $u$ (see e.g. \cite{Nir} or the end of Remark \ref{R-generalization}).
We denote by $ B( N, \si )$ the best possible value of the constant $B$ in (\ref{GNS}), namely 
\beq
\label{maxGNS}
B( N, \si ) = \sup_{ u \in H^2 ( \R^N), \; u \neq 0 } \; \; 
\frac{ \| u \|_{ L^{ 2 \si + 2 } }^{ 2 \si + 2 } }{ \| \Delta u \|_{L^2 }^{ \frac{ \si N}{2}}  \|  u \|_{L^2 }^{  2 + 2 \si -\frac{ \si N}{2}} }.
\eeq
It is also well known that there exist optimal functions for (\ref{maxGNS}); that is, the supremum in (\ref{maxGNS}) is, in fact, a maximum
(see, e.g., Lemma \ref{L3.13} (ii) below). 

\medskip

Let $E_{min}$ be as in (\ref{EMIN}).
The properties of the function $E_{min}$ will be crucial in the sequel. They are summarized in the next Proposition and in the  remark following it.

\begin{Proposition}
\label{Emin}
The function $ m \longmapsto E_{min}(m )$ has the following properties: 

\medskip

(i) If $ \si N > 4 $ we have $ E_{min}(m ) = - \infty $ for all $ m > 0 $. 

\medskip

\noindent
For the following statements we assume that $ 0 < \si N \les 4$. We have: 

\medskip

(ii) The function $E_{min} $ is concave on $(0, \infty )$. 

\medskip

(iii) For any $ m > 0 $ there holds $ E_{min}( m ) \les - m. $

\medskip

(iv) 
 $ \ds \lim _{ m \downarrow 0 } E_{min}( m ) = 0 $ and $ \ds \lim _{ m \downarrow 0 } \frac{E_{min}(m )}{m } 
= -1 $.

\medskip

(v) If $  0 < \si N <4$ we have $ E_{min}(m ) > - \infty $ for all $ m > 0 $ and there exist $ A\in \R$, $ B > 0 $ such that 
$  E_{min}( m ) <  A m - B m ^{ \si + 1}$ (thus, in particular, 
$\ds \frac{ E_{min}( m ) }{m } \lra - \infty $ as {$ m \lra \infty$}). 
Moreover, for any $ k _1 , \; k_2 >0$ the set $ \{ u \in H^2 ( \R^N) \; | \; \| u \|_{L^2 } \les k _1 \; \mbox{ and } \; E( u ) \les k_2 \}$ 
is bounded in $ H^2 ( \R^N)$. 

\medskip

(vi) Assume that $ \si  N = 4$. 
Let $B(N, \si )$ be as in (\ref{maxGNS}) and let $ k_* = ( \si + 1 ) ^{\frac{1}{\si}} B(N, \si )^{ - \frac{1}{\si}}. $

Then $ E_{min}( m  ) $ is finite for any $ m \in (0, k_* ) $ and $ E_{min}(m ) = - \infty $ if $ m \ges k_* .$

In addition, for any $ k_1  < k_*  $ and any $ k_2 > 0 $,   the set 
$ \big\{ u \in H^2 ( \R^N) \; | \; \; \| u \|_{L^2 }^2 \les k_1  \; \; \mbox{ and } \; \; E( u ) \les k_2 \big\}$ 
is bounded in $ H^2 ( \R^N)$.

\end{Proposition}

\begin{remark}
\label{PropEmin}
\rm 
The function $E_{min}$ is finite and concave on $ (0, \infty )$ if $ \si N < 4 $,  respectively  on $( 0, k_* ) $ if $ \si N = 4$, 
hence it is continuous  and admits left and right derivatives at any point of these intervals. 
We denote by $E_{min, \ell}'(m) $ and $E_{min, r}'(m)$, respectively, the left and right derivatives of $ E_{min} $ at $ m $. 
The functions  $E_{min, \ell}' $ and $E_{min, r}' $ are nonincreasing,  we have $ E_{min, \ell}' (m) \ges E_{min, r}' (m)$ for all $m$ and
equality must occur at all but countably many $m$'s. 
Proposition \ref{Emin} (iv) implies that 
$$
 \lim_{ m \downarrow 0 } E_{min, \ell}'(m) = \sup_{m > 0 } E_{min, \ell}'(m) =  \lim_{ m \downarrow 0 } E_{min, r}'(m) 
 = \sup_{m > 0 } E_{min, r}'(m) = -1 .
$$
Let 
\beq
\label{m0}
m _0 := \sup \{ m > 0 \; | \; E_{min}(m ) = - m \}. 
\eeq
It is clear that $  E_{min}(m ) = - m $ on $ (0, m_0)$ and $  E_{min}(m ) < - m $ on $ (m_0, \infty)$. 
If $ m > m_0 $ and $ E_{min}(m) > -\infty  $ we must have $ E_{min, \ell}' ( m ) < -1$.
If $ \si N < 4$  we have $ m _0 < \infty $ and 
$ \ds \lim_{m \ra \infty }  E_{min, \ell}'(m)  =  \lim_{m \ra \infty }  E_{min, r}'(m) = - \infty $ because 
$  {\ds \lim_{m \ra \infty } } \frac{ E_{min}(m ) }{m} = - \infty$ by 
Proposition \ref{Emin} (v).

\end{remark}

\medskip 

{\it Proof of Proposition \ref{Emin}. } 
(i) Let $ m > 0 $.  Choose  $ u \in H^2 ( \R^N)$ satisfying $ \| u \|_{L^2 } ^2= m $. 
We use (\ref{scaling}) with $ u_{a, b } = a u \left( \frac{\cdot}{b}\right)$ and $ a =  t ^{\frac N4}$, $ b = t ^{- \frac 12 } $.
It is obvious that $ \| u _{t ^{ N/ 4}, \,  t ^{-  1/2 } } \| _{L^2 } ^2 =  \| u \|_{L^2 } ^2= m$ for any $ t > 0 $, and consequently 
$ E_{min}(m ) \les E(u _{t ^{ N/ 4}, \,  t ^{-  1/2 } }) $ for all $t$. 
From (\ref{scaling}) we have 
\beq
\label{scalm}
E(  u _{t ^{ N/ 4}, \,  t ^{-  1/2 } }) = t^2 \int_{\R^N} |\Delta u |^2 \, dx - 2 t \int_{\R^N} |\nabla u |^2 \, dx 
- t ^{ \frac{N \si}{2}} \int_{\R^N} | u |^{2 \si + 2} \, dx . 
\eeq
If $ N \si > 4$, letting $ t \lra \infty $ we discover $  E_{min}(m ) \les \ds \lim_{t \ra \infty } E(u _{t ^{ N/ 4}, \,  t ^{-  1/2 } }) = -\infty .$

\medskip

(ii)  It is obvious that $ u \in S( m ) $ if and only if there exists $ v \in S(1)$ such that $ u = \sqrt{ m } v$. 
Hence for any $ m > 0 $ we have 
$$
\begin{array}{l}
E_{min}(m)   =  \inf \{ E( \sqrt{m } v ) \; | \; v \in S(1) \} \\ \\
 =  \ds \inf \left\{ m \left(   \int_{\R^N}|\Delta v|^2\, dx - {2}\int_{\R^N}|\nabla v|^2\, dx \right) - \frac{m^{ \si + 1 }}{\si + 1} \int_{\R^N}|v|^{2\sigma+2} \, dx \; \; \Big| \; v \in S(1)  \right\}.
\end{array}
$$
For any $ A \in \R$ and any $ B \ges 0 $ the function $ m \longmapsto A m - B m^{\si + 1 } $ is concave on $ (0,\infty)$. 
The infimum of a family of concave functions is also a concave function and statement (ii) follows. 

\medskip

(iii)  
Let $ m > 0 $ and $ \e > 0 $. 
Choose a function $ \eta \in C_c ^{\infty }( \R^N)$ such that 
{$ \| \eta \|_{L^2 }^2 = (2 \pi )^N m$} and the support of $ \eta $ is contained in the annulus 
$ B( 0, 1  ) \setminus B(0, 1 - \e)$.
Let $ u = \Fo ^{-1} ( \eta)$. Then  $ u \in \So ( \R^N)$ and $ \| u \|_{L^2}^2 = \frac{1}{(2 \pi )^N} \| \eta \|_{L^2}^2 = m $.
Using the basic properties of the Fourier transform,  Plancherel's formula and the fact that 
$ 0 \les  \left( |\xi |^2 - 1 \right)^2 \les 4 \e ^2 $ 
on the support of $ \eta$  we get 
$$
\begin{array}{l}
 \ds \int_{\R^N}|\Delta u|^2\, dx - {2}\int_{\R^N}|\nabla u|^2\, dx + \int_{\R^N}| u|^2\, dx
= \frac{1}{(2 \pi )^N} \int_{\R^N} \left( |\xi |^4 - 2 |\xi |^2 + 1 \right) |\wh{u}(\xi )|^2 \, d\xi 
\\ \\
= \ds \frac{1}{(2 \pi )^N}  \int_{B( 0, 1 ) \setminus B( 0, 1 - \e ) } \left( 1 - |\xi |^2  \right)^2  |\eta(\xi )|^2 \, d\xi 
\les \frac{ 4 \e^2 }{(2 \pi )^N}  \int_{\R^N}  |\eta(\xi )|^2 \, d\xi 
=  4 \e^2  m. 
\end{array}
$$
We infer that 
$$
E_{min}(m ) + m \les E(u) + \| u \|_{L^2}^2 \les 
\int_{\R^N}|\Delta u|^2\, dx - {2}\int_{\R^N}|\nabla u|^2\, dx + \int_{\R^N}| u|^2\, dx 
\les  4 \e^2  m, 
$$
that is $ E_{min}(m ) \les -m + 4 \e^2 m. $ Since $ \e > 0 $ is arbitrary, (iii) follows. 

\medskip

(iv) Consider first the case $ 0 < \si N < 4$. 
Let $  0 < \e <1 $. 
Using the Gagliardo-Nirenberg-Sobolev inequality (\ref{GNS}) and Young's inequality 
($|ab | \les \frac{ |a|^p }{p } + \frac{|b |^q }{q} $ if $ \frac 1p + \frac 1q = 1$)  with exponents 
$ p = \frac{4}{\si N } $ and $ q = \frac{4}{4 - \si N }$, we get for any $v \in H^2 ( \R^N)$
$$
\frac{1}{\si + 1 } \| v\|_{L^{ 2 \si + 2}}^{2 \si + 2 } 
\les B \frac{1}{\si + 1 }  \| \Delta v \|_{L^2 }^{ \frac{ \si N}{2}}  \|  v \|_{L^2 }^{2 + 2 \si - \frac{ \si N}{2}} 
\les \e \| \Delta v \|_{L^2 }^2 + C_1 ( \e ) \|  v \|_{L^2 }^{ \frac{8 ( \si + 1) - 2 \si N }{4 - \si N}}, 
$$
where $ C_1 ( \e )$ is independent of $v$. 
It follows that 
\beq
\label{Eninf}
E(v) \ges ( 1 - \e ) \| \Delta v \|_{L^2 }^2 - 2 \| \nabla v \|_{L^2 }^2 - C_1 ( \e ) \|  v \|_{L^2 }^{ \frac{8 ( \si + 1) - 2 \si N }{4 - \si N}}.
\eeq
Using Plancherel's formula we get 
\beq
\label{Planch}
( 1 - \e ) \| \Delta v \|_{L^2 }^2 - 2 \| \nabla v \|_{L^2 }^2 + \frac{1}{1  - \e } \| v \|_{L^2}^2 
= \frac{1 - \e }{(2 \pi)^N} \int_{\R^N} \left( |\xi |^2 - \frac{1}{ 1 - \e } \right)^2 |\wh{v} ( \xi )|^2 \, d \xi \ges 0 . 
\eeq
Notice that  the inequality in (\ref{Planch}) is strict  if $ u \neq 0$.
From (\ref{Eninf}) and (\ref{Planch}) we get 
\beq
\label{lim0}
E(v) \ges - \frac{1}{1 - \e} \| v \|_{L^2 }^2 - C_1 ( \e ) \|  v \|_{L^2 }^{ \frac{8 ( \si + 1) - 2 \si N }{4 - \si N}} \qquad 
\mbox{ for all } v \in H^2 ( \R^N). 
\eeq
Taking the infimum in (\ref{lim0}) over all $ v \in H^2( \R^N)$ satisfying $ \| v \|_{L^2}^2 = m $ we discover
\beq
\label{lim1}
E_{min}(m ) \ges - \frac{m}{1 - \e} - C_ 1 ( \e ) m^ {\frac{4 ( \si + 1) -  \si N }{4 - \si N}} \qquad \mbox{ for any } m > 0 . 
\eeq
From (iii) and (\ref{lim1}) it follows that $ E_{min}( m ) \lra 0 $ as $ m \lra 0 $ and 
$ \ds - \frac{1}{ 1- \e} \les  \liminf_{m \downarrow 0 } \frac{E_{min}(m)}{m} \les \limsup_{m \downarrow 0 } \frac{E_{min}(m)}{m} \les -1. $
Since $ \e $ is arbitrary, (iv) is proven. 

Next consider the case $ \si  N = 4$. 
The Gagliardo-Nirenberg-Sobolev inequality (\ref{GNS}) becomes 
\beq
\label{GNS2}
\| v \|_{L^{2 \si + 2 }}^{ 2\si + 2 } \les B \| \Delta v \|_{L^2 }^2 \| v \|_{L^2 }^{ 2 \si } . 
\eeq
Let $ 0 < \e < 1 $. Using (\ref{GNS2}), for any $ v \in H^2( \R^N)$ satisfying $ \frac{B}{\si + 1 } \| v \|_{L^2}^{2 \si } \les \e $ we get 
$$
E(v ) \ges ( 1 - \e ) \| \Delta v \|_{L^2 }^2 - 2 \| \nabla v \|_{L^2 }^2 
$$
and then using (\ref{Planch}) we obtain $E(v ) \ges - \frac{1}{ 1- \e } \|v \|_{L^2 }^2 . $
This gives $ E_{min}(m) \ges -  \frac{m}{ 1- \e } $ for all $ m > 0 $ satisfying  $ \frac{B}{\si + 1 } m^{ 2 \si } \les \e$. 
Taking into account (iii), statement (iv)  is now  obvious.

\medskip

(v) Using (\ref{lim1}) with $\e = \frac 12$, say, it is clear that $E_{min}(m) > -\infty $ for any $ m > 0 $. 

Fix $u \in H^2( \R^N)$ such that $ \| u \|_{L^2} = 1$. Let $ A = \ds \int_{\R^N } | \Delta u | ^2 \, dx - 2 \int_{\R^N } | \Delta u | ^2 \, dx $ and 
$ B = \ds \frac{1}{\si + 1} \int_{\R^N } | u | ^{2 \si + 2}  \, dx$. 
For all $ m > 0 $ we have $ \| m^{\frac 12} u \|_{L^2}^2 = m $, 
hence $ E_{min}(m) \les E(u ) = A m - B m^{\si + 1}$.

If $ \| u \| _{L^2 } \les k_1 $ and $ E(u ) \les k_2$, using (\ref{interpol}) and (\ref{GNS}) we get 
$$
k_ 1 \ges \| \Delta u \|_{L^2}^2 - 2 k_1  \| \Delta u \|_{L^2} - \frac{ B}{\si + 1} k_1 ^{2 + 2 \si - \frac{ \si N }{2} } \| \Delta u \|_{L^2}^{\frac{\si N}{2}}.
$$
Since $\frac{\si N}{2} < 2$, the above inequality implies that $ \| \Delta u \|_{L^2} $ is bounded. 
Then (\ref{interpol}) and the inequality $ \Big\| \frac{\p^2 u }{\p x_i \p x _j } \Big\|  _{L^2} \les \| \Delta u \|_{L^2}$ imply that 
$ \| u \|_{H^2( \R^N)}$ is bounded. 

\medskip

(vi) Assume $ \si  N = 4$. 
By (\ref{interpol}) and (\ref{GNS}) we have for all $  u \in H^2(\R^N)$
\beq
\label{Ecrit}
\begin{array}{rcl}
E(u ) 
& \ges & \| \Delta u \|_{L^2}^2 - 2 \| u \|_{L^2} \| \Delta u \|_{L^2} - \frac{B(N, \si )}{\si + 1} \| \Delta u \|_{L^2}^2 \| u \|_{L^2}^{ 2 \si}
\\
\\
 & = & \left( 1 - \frac{B(N, \si )}{\si + 1} \| u \|_{L^2}^{ 2 \si} \right) \| \Delta u \|_{L^2}^2 - 2 \|
 u \|_{L^2} \| \Delta u \|_{L^2} .
 \end{array}
\eeq

Let $ k_1 < k_* = ( \si + 1 ) ^{\frac{1}{\si}} B(N, \si )^{- \frac{1}{\si}}$. Let $ \tau ( k_1 ) = 1 -  \frac{B(N, \si )}{\si + 1} k_1 ^{ \si } > 0 $. 
For any $ u \in H^2 ( \R^N)$ such that $ \| u \|_{L^2} ^2 \les k_1 $,  by (\ref{Ecrit}) 
we get 
\beq
\label{Ecrit1}
E(u ) \ges \tau( k_1) \| \Delta u \|_{L^2}^2 - 2 k_1  \| \Delta u \|_{L^2} \ges  { \ds \min _{s \in \R}  ( \tau( k_1) s ^2 - 2 k_1 s )  } = - \frac{k_1 ^2}{ \tau ( k_1 )}. 
\eeq
We infer that $ E_{min}( m ) \ges - \frac{k_1 ^2}{ \tau ( k_1 ) }> - \infty $ for all $ m \in (0, k_1]$. 
Since $ k_1 < k_* $ was arbitrary, we see that $E_{min} $ is finite on $(0, k_*)$. 
Moreover, if $ \| u \|_{L^2} ^2 \les k_1 < k_* $ and $ E(u) \les k_2 $ then (\ref{Ecrit1}) implies that $ \| \Delta u \|_{L^2}$ is bounded and arguing as in part (v)  we see that 
$ \| u \|_{H^2( \R^N)}$ is bounded.

Let $Q$ be an optimal function for the Gagliardo-Nirenberg-Sobolev inequality (\ref{GNS}) with $ \si = \frac 4N$ such that 
$ \| Q \|_{L^2 } = k_*^{\frac 12} = ( \si + 1 ) ^{\frac{1}{2\si}} B(N, \si )^{- \frac{1}{2\si}}.$ 
Such a function $Q$ exists because whenever $u$ is an optimal function for (\ref{GNS}), the rescaled functions 
$ u_{a, b }(\cdot ) = a u \left( \frac{\cdot}{b} \right)$
are   optimal functions, too. 
We have 
$$ \frac{1}{ \si + 1 } \| Q \| _{L^{2 \si + 2}}^{2 \si + 2 }
 = \frac{1}{ \si + 1 }  B(N, \si ) \| \Delta Q \|_{L^2 }^{ 2 }  \|  Q \|_{L^2 }^{   2 \si } = \| \Delta Q \|_{L^2 }^{ 2 } .$$
For $ t > 0 $ let $ u_t ( x ) = t^{\frac N4} Q( t ^{\frac 12 } x)$. 
From (\ref{scaling}) and (\ref{scalm}) it follows that $ \| u_t \|_{L^2}^2 = \| Q \|_{L^2 } ^2 = k_* $ and 
$E( u_t ) = - 2 t \| \nabla Q \|_{L^2} ^2. $ 
Letting $ t \lra \infty $ we discover $ E_{min}( k_*) = - \infty$. 
If $ m > k_*$, using the test function $  m ^{\frac 12}  k_* ^{- \frac 12 }   u_t$ and letting $ t \lra \infty $ we find 
$ E_{min}( m ) = -\infty$.  
\hfill
$\Box$

\medskip

As we will see later, problem $(\Po _m)$ admits solutions if and only if $ -\infty < E_{min}(m) < -m$. 
We have already seen that $  E_{min}(m) = - \infty $ for all $m$ if $ \si > \frac 4N$. 
Proposition \ref{P-interesting} gives necessary and sufficient conditions to have $E_{min}(m) < -m$
whenever  $ \si \in (0, \frac 4N]$. 
Its proof relies on the functional inequalities proved in Section \ref{inequalities} and on the test functions constructed there.

\begin{Proposition}
\label{P-interesting}

Let $ N \in \N^*$. We have: 

\medskip

(i) If $ 0 < \si < \max \left( 1, \frac{4}{N+1} \right), $ then $E_{min}(m) < - m $ for all $ m > 0 $. 

\medskip

(ii) If $ \max \left( 1, \frac{4}{N+1} \right) \les \si \les \frac{4}{N}$, there exists $ m_0 > 0 $ 
such that $ E_{min}(m ) = - m $ for any $ m \in (0, m_0]$ and  $E_{min}(m) < - m $ for any $ m >m_ 0 $.
Moreover, $m_0$ is given by (\ref{explicitm0}) below.

\medskip

(iii) Assume that $ \si = \frac 4N$ and let $ k_*$ be as in Proposition \ref{Emin} (vi). 
Let $m_0$ be as in  (\ref{explicitm0}). Then   $ m_0 < k_*$ and   we have
$ E_{min}(m) = -m $ if $ m \in (0, m_0]$, respectively $ -\infty < E_{min}(m) < -m $ if $ m \in (m_0, k_*)$.

\end{Proposition}

{\it Proof. }
(i) Assume that $N = 1$ and  $ 0 < \si < 2$. 
Fix  $ m > 0 $ and let $ u _{\tau }$ be as in Example \ref{Ex-dim1}.
 Denote $ v_ {\tau } = \pi^{- \frac 14} \tau ^{-\frac 12 } m ^{\frac 12} u_{\tau }$.
By (\ref{ExD1-1}) we have 
$\|v_{\tau } \|_{L^2}^2 = m $, 
$\| v_{\tau }  \|_{L^{2\si + 2}}^{2\si + 2} = 
( \si + 1)^{- \frac 12} \pi^{- \frac{ \si}{2} } m^{\si + 1} \tau ^{ - \si } $
and (\ref{ExD1-2}) gives 
$\| ( \Delta + 1 )v_{\tau }  \|_{L^2}^2 
\sim Cm \tau^{-2}$ as $ \tau \lra \infty$. 
For $ \tau $   sufficiently large we have 
$ E(v_{\tau }   )   + \| v_{\tau }  \|_{L^2}^2  = \| ( \Delta + 1 )v_{\tau }  \|_{L^2}^2  - 
\frac{1}{\si + 1}\| v_{\tau }  \|_{L^{2\si + 2}}^{2\si + 2}
< 0 $, 
and this implies $ E_{min}(m ) < - m $.

If $ N \ges 2$, for small $ \e, \de > 0 $ let $ u_{\e, \de}$ be as in (\ref{u-epsilon}). 
Denote $ w_{\e, \de}= \frac{\sqrt{m}}{\| u _{\e, \de}\|_{L^2}}  u_{\e, \de}$, so that $ \|  w_{\e, \de}\|_{L^2}^2 = m $.
By (\ref{ex-1}) and (\ref{ex-2}) we have 
$ \| ( \Delta + 1 )w_{\e, \de }  \|_{L^2}^2 \sim m \e ^2$,  
while (\ref{ex-1}) and (\ref{ex-4}) give
$\| w_{\e, \de }  \|_{L^{2\si + 2}}^{2\si + 2 } \ges C m^{\si + 1} \e^{\si+1} \de^{\si ( N -1)} ( \de ^2 + \e )^{-1} $ for some $C>0$.

If $ \si \in (0, 1)$, fix a small $ \de_0 > 0 $   and observe that
\beq
\label{wedee}
  E(w_{\e, \de_0 })  +  \| w_{\e, \de_0 }  \|_{L^{ 2}}^{2} =  \| ( \Delta + 1 )w_{\e, \de _0 }  \|_{L^2}^2 -\frac{1}{\si + 1}  \| w_{\e, \de_0 }  \|_{L^{2\si + 2}}^{2\si + 2} 
 < 0 
\eeq  if  $ \e $ is sufficiently small, 
 hence $  E_{min}(m ) + m < 0 $. 

If $ \si \in \left(0, \frac{4}{N + 1}\right)$, taking $ \e = \de^2$ it follows from the above estimates that  
$ \| ( \Delta + 1 )w_{\de ^2, \de }  \|_{L^2}^2 \sim m \de ^4$ and 
$\| w_{\de^2, \de }  \|_{L^{2\si + 2}}^{2\si + 2 } \ges C m^{\si + 1}  \de^{\si ( N +1)}$ for some $C> 0 $ and any small $ \de > 0 $. 
As in (\ref{wedee}), this implies $ E( w_{\de^2, \de }) + \| w_{\de^2, \de } \|_{L^2}^2 < 0 $ for sufficiently small $ \de$, 
and consequently $  E_{min}(m ) + m < 0 $.

\medskip

(ii)  
It is easy to see that for any $ u \in H^2( \R^N)$ we have
$$
E(u ) + \| u \|_{L^2}^2 
= \| ( \Delta + 1 ) u \|_{L^2}^2 - \frac{1}{\si + 1} \| u \|_{L^{2 \si + 2}}^{2 \si + 2}
=  \ds \| ( \Delta + 1 ) u \|_{L^2}^2 \left( 1 - \frac{\| u \|_{L^2}^{ 2 \si }}{\si + 1}
 \frac{ \| u \|_{L^{2 \si + 2}}^{2 \si + 2}}{ \| u \|_{L^2}^{ 2 \si } \| ( \Delta + 1 ) u \|_{L^2}^2} \right).
 $$
Let $  \ka \! = \! \frac{\si}{\si + 1}$ and $  \ds Q_{\ka}(u ) \! =  \! \frac{ \| u \|_{L^{2 \si + 2}} }
{ \| u \|_{L^2}^{ \ka } \| ( \Delta \! +\!  1  ) u \|_{L^2}^{1 - \ka}}$ (see (\ref{Qka})).
The above equality can be written~as 
\beq
\label{energy-rewritten}
E(u ) + \| u \|_{L^2}^2 = 
\| ( \Delta + 1 ) u \|_{L^2}^2 \left( 1 - \frac{\| u \|_{L^2}^{ 2 \si }}{\si + 1} Q_{\ka}(u )^{2 \si + 2} \right)
\quad \mbox{ for all } u \in H^2( \R^N) \setminus \{ 0 \}. 
\eeq
We use the results in Section \ref{inequalities} with $ s = 2$, $ p = 2 \si + 2 $ and $ \ka = \frac{\si}{\si + 1}$.

If $ N=1$, 
 condition $ \frac{1}{2 s} \les  {\frac{( 1 - \ka )p}{p-2 }} \les \frac 12$ in Theorem \ref{T-1D} is equivalent to $ 4 \ges \si \ges 2$. Hence $ Q_{\ka}$ is bounded from above if $ \si \in [2, 4]$.

If $ N \ges 2$, the condition 
$ \ka \ges \frac12  $ in (\ref{conditions})
is equivalent to $\si \ges 1$, 
$
 \frac{N}{s}\left(\frac 12 - \frac{1}{p} \right)  \les 1 - \ka
$
is equivalent to $  \si \les \frac 4N$, 
and 
 $ 1 - \ka \les    \frac{N+1}{2}\left(\frac 12 - \frac{1}{p} \right)$ is equivalent to 
 $ \si \ges \frac{4}{N+1}$.
 By Theorem \ref{T-multiD}, $ Q_{\ka}$ is bounded from above if and only if 
 $ \max \left( 1, \frac{4}{N+1} \right) \les \si \les \frac 4N$.

Whenever $Q_{\ka}$ is bounded from above, let $M = \ds \sup_{ u \in H^s ( \R^N) \setminus \{ 0 \} } Q_{\ka } ( u ) $, as in  (\ref{supr}), and let 
\beq
\label{explicitm0}
m_0 = ( \si + 1 )^{\frac{1}{\si}} M^{- \frac{ 2 \si + 2}{\si}}. 
\eeq
If $m \in (0, m_0]$, using (\ref{energy-rewritten}) we infer that $ E(u) + \| u \|_{L^2}^2 \ges 0 $ for any $ u \in H^2( \R^N)$ satisfying $ \| u \|_{L^2}^2 = m $, hence $E_{min}(m) + m \ges 0 $. 
Then Proposition \ref{Emin} (iii) implies  $E_{min}(m) = - m$.

If $ m > m_0$,  we have $ ( \si + 1)^{\frac{1}{2 \si + 2}} m ^{-\frac{\si}{2 \si + 2}} < M$. 
Choose $ u \in H^2( \R^N)$, $ u \neq 0 $, such that $ Q_{\ka }(u) >  ( \si + 1)^{\frac{1}{2 \si + 2}} m ^{-\frac{\si}{2 \si + 2}} $. 
Let $ v = \frac{\sqrt{m}}{\| u \|_{L^2} } u $, so that $ \| v \|_{L^2}^2 = m $ and $ Q_{\ka }(v ) = Q_{\ka }(u)$.
From (\ref{energy-rewritten}) we get $ E(v) + \| v \|_{L^2}^2 < 0 $, hence 
 $E_{min}(m) < - m$.

\medskip

(iii) Taking into account (\ref{explicitm0}) and the expression of $ k_*$ in Proposition \ref{Emin} (vi), 
the inequality $ m_0 < k_*$ is equivalent to 
 $ B(N, \si ) <  M^{ 2 \si + 2 }$, where $ B(N, \si )$ is given by (\ref{maxGNS}).
Denote by $ \Qo (u)$ the quotient appearing in (\ref{maxGNS}).
Let $ u_* $ be an optimal function for (\ref{maxGNS}).
Then $ u_{\tau } = u_*\left(\frac{ \cdot}{\sqrt{\tau}} \right) $ is also an optimal function for (\ref{maxGNS}), that is 
$ \Qo ( u_{\tau}) = B(N, \si)$ for all $ \tau > 0 $. 
The conclusion follows if we find $ \tau > 0 $ such that $  \Qo ( u_{\tau}) < Q_{\ka }( u _{\tau})^{2 \si + 2}$, and this is equivalent to $ \|  ( \Delta + 1) u_{\tau} \|_{L^2}^2 < \| \Delta u _{\tau } \|_{L^2}^2$, 
or using Plancherel's theorem, 
$\ds \int_{\R^N} \left( |\xi |^2 - \tau \right) ^2 |\wh{u }_* |^2 ( \xi ) \, d \xi < 
 \int_{\R^N} |\xi |^4  |\wh{u }_* |^2 ( \xi ) \, d \xi. $
The last inequality can be written as 
$$ \ds - 2 \tau \int_{\R^N} |\xi |^2  |\wh{u }_* |^2 ( \xi ) \, d \xi + 
\tau ^2  \int_{\R^N}  |\wh{u }_* |^2 ( \xi ) \, d \xi < 0 $$
 and  holds  true if 
 $ 0 < \tau < \frac{2 \| \, | \cdot | \wh{u} _* \|_{L^2}^2 }{      \| \wh{u} _* \|_{L^2}^2} $.
We have thus shown that $ m_0 < k_*$. 
The rest follows from part (ii) and Proposition \ref{Emin} (vi).
 \hfill
 $\Box$

\medskip

The next Theorem establishes the existence of minimizers for the problem $(\Po _m )$ as well as the pre-compactness modulo translations of 
all minimizing sequences.

\begin{Theorem}
\label{Global}
Assume that $ N \si < 4 $ and $m > 0 $ is such that $ E_{min}(m) < -m. $

Then for any sequence $ (u_n)_{n \ges 1} \subset H^2( \R^N)$ satisfying $ M( u_n) \lra m $ and $ E( u_n) \lra E_{min}(m)$ 
there exist a subsequence, still denoted $( u_n)_{n \ges 1}$, a sequence of points $(x_n )_{n \ges 1 } \subset \R^N$ and a function 
$ u \in H^2 ( \R^N)$ such that 
$ u_n ( \cdot + x_n) \lra u $ strongly in $ H^2( \R^N)$. 

In particular,  there exists a solution $ u \in H^2( \R^N)$ to the minimization problem $(\Po_m)$.

The same conclusion holds if $ N \si = 4$, $ 0 < m < k_*$ (where $ k_*$ is as in Proposition \ref{Emin} (vi)),  and $ E_{min}(m ) < - m $. 
\end{Theorem}

{\it Proof. } 
Let $( u_n)_{n \ges 1 }$ be a minimizing sequence.
It follows from Proposition \ref{Emin} (v) or (vi) that $(u_n)_{n \ges 1}$ is bounded in $H^2(\R^N)$. 

Using (\ref{Planch}) with $ \e = 0 $ we infer that for any $ u \in H^2( \R^N)$ there holds
\beq
\label{nonlin}
\begin{array}{l}
\ds \int_{\R^N} |u |^{ 2 \si + 2} \, dx 
 =  ( \si + 1 ) \int_{\R^N} |\Delta  u |^2 - 2 |\nabla u |^2 + |u |^2 \, dx - ( \si + 1 ) ( E(u ) + \| u \|_{L^2}^2 )
\\
\\
 \ges   - ( \si + 1 ) ( E(u ) + \| u \|_{L^2}^2 ). 
\end{array}
\eeq
Choose $ \si '> \si $ such that $ 2 \si '+ 2 < 2^{**}$. We denote by $ \Lo ^N$ the Lebesgue measure in $ \R^N$. 
Using H\"older's inequality and the Sobolev embedding we get for any $ u \in H^2( \R^N)$  and for any $ t > 0 ,$
\beq
\label{nontriv}
\begin{array}{l}
\int_{\R^N} |u |^{ 2 \si + 2} \, dx  
 =  \int_{\{ | u | < t \} } |u |^{ 2 \si + 2} \, dx +  \int_{\{ | u | \ges t \} }    |u |^{ 2 \si + 2} \, dx   
\\
\\
  \les  t^{ 2 \si } \int_{\{ | u | < t \} } |u |^{ 2} \, dx 
+   \left( \int_{\{ | u | \ges t \} }    |u |^{ 2 \si '+ 2} \, dx   \right)^{\frac{\si + 1}{\si ' + 1} }
\Lo ^N \left( \{ |u | \ges t \} \right ) ^{1 - \frac{\si + 1}{\si ' + 1} }
\\
\\
 \les 
t^{ 2 \si } \| u \|_{L^2 }^2 + \left(C_S  \| u \|_{H^2} \right)^{ 2 \si + 2 } \Lo ^N \left( \{ |u | \ges t \} \right ) ^{1 - \frac{\si + 1}{\si ' + 1} } .
\end{array}
\eeq
Choose $ \de > 0 $ such that $ 2 \de < - ( \si + 1 ) ( E_{min}(m) + m ) $ (this is possible because $ E_{min}(m) < - m $). 
Since $ \| u _n \|_{L^2 }^2 \lra m $ and $ E( u_n ) \lra E_{min}(m)$, (\ref{nonlin}) implies that 
$\ds \int_{\R^N} | u_n |^{ 2 \si + 2 } \, dx > 2 \de $ for all $ n $ sufficiently large. 
Choose $ t_0 > 0 $ such that $ t_0 ^{ 2 \si } ( m + 1 ) = \de$. 
Using (\ref{nontriv}) for $ u_n $ and the boundedness of $ ( u_n )_{ n \ges 1 } $ in $ H^2( \R^N)$, 
we infer that there exists a constant $ a > 0 $, independent of $n$, such that 
$ \Lo ^N \left( \{ | u _n | \ges t_0 \} \right) \ges a $ for all sufficiently large $n$.
Using Lieb's Lemma (see Lemma 6 p. 447 in \cite{lieb} or Appendix 2 in \cite{MM}) 
we infer that there exists a constant $ b > 0 $, independent of $n$,  and for each $n$ large there exists $ x_n \in \R^N$ such that 
$$
 \Lo ^N \left( \left\{ x \in B( x_n, 1) \; \big| \; | u _n | \ges \frac{t_0}{2} \right\} \right) \ges b .
$$ 
From now on we replace $ u_n $ by $ u_n ( \cdot +  x_n)$, which is still a minimizing sequence and satisfies 
$  \Lo ^N \left( \{ x \in B(0, 1) \; \big| \; | u _n | \ges \frac{t_0}{2} \} \right) \ges b .$
Since $ ( u_n )_{n \ges 1}$ is bounded in $ H^2( \R^N)$ there exists a subsequence, still denoted $ ( u_n )_{n \ges 1}$, 
and there is $ u \in H^2( \R^N)$ such that 
\beq
\label{conv}
\begin{array}{l}
u _n \rightharpoonup u \qquad \mbox{ weakly  in }  \; H^2( \R^N), 
\\
u_n \lra u \qquad \mbox{ in } L_{loc}^p ( \R^N ) \mbox{ for } 1 \les p < 2^{**} \mbox{  and a.e. } 
\end{array}
\eeq

It is clear that $ \int_{B( 0, 1 ) } | u_n |^p \, dx \ges b \left( \frac{t_0}{2} \right)^p $ for all $ n $ sufficiently large. 
Take any $ p \in [1, 2^{**})$ and pass to the limit 
to get 
$ \int_{B( 0, 1 ) } | u |^p \, dx \ges b \left( \frac{t_0}{2} \right)^p $. In particular, we infer that $ u \neq 0 $. 

Let $ m_1 = \| u \|_{L^2}^2$. It is clear that $ 0 < m_1 \les \liminf_{n \ra \infty } \| u _n \|_{L^2}^2 = m . $
We will show that $ m_1 = m $. 
We argue by contradiction and we assume that $ m_1 < m $. 
The weak convergence in a Hilbert space gives as $ n \lra \infty $ 
\beq
\label{conv1}
\begin{array}{c}
\| u _n \|_{L^2 }^2 = \| u \|_{L^2}^2 + \| u_n - u \|_{L^2}^2 + o(1), 
\\
\\
\| \nabla u _n \|_{L^2 }^2 = \| \nabla u \|_{L^2}^2 + \|\nabla (u_n - u ) \|_{L^2}^2 + o(1), 
\\
\\
\| \Delta u _n \|_{L^2 }^2 = \| \Delta u \|_{L^2}^2 + \| \Delta ( u_n - u ) \|_{L^2}^2 + o(1) .
\end{array}
\eeq
Using the Brezis-Lieb Lemma (see e.g. Lemma 4.6 p. 10 in \cite{kavian}) we get 
\beq
\label{conv2}
\int_{\R^N} | u _n |^{ 2 \si + 2} \, dx = \int_{\R^N} | u  |^{ 2 \si + 2} \, dx  + \int_{\R^N} | u _n - u |^{ 2 \si + 2} \, dx + o(1). 
\eeq

From (\ref{conv1}) and (\ref{conv2}) we get 
\beq
\label{conv3}
E( u_n ) = E(u) + E( u_n - u ) +   o (1) \qquad \mbox{ as } n \lra \infty . 
\eeq
It is obvious that $ E(u) \ges E_{min} ( m_1) $ and  $ E( u_n - u ) \ges E_{min} ( \| u_n - u \|_{L^2} ^2) $. 
By (\ref{conv1}) we have $ \| u_n - u \|_{L^2} ^2 \lra m - m_1$. 
The function $E_{min}$ is continuous on $(0, \infty)$ if $ 0 < \si N < 4 $, respectively   on $ (0, k_*)$ if $ \si N = 4$, 
and passing to the limit in (\ref{conv3}) we get 
\beq
\label{conv4}
E_{min}(m) \ges E_{min} ( m_1) + E_{min} ( m- m_1). 
\eeq
Since $E_{min} $ is concave and $ E_{min} ( \eta ) \lra 0 $ as $ \eta \lra 0 $ we have 
$ E_{min} ( m_1 ) \ges \frac{m_1}{m} E_{min} ( m ) $ and $ E_{min} ( m - m_1 ) \ges \frac{m - m_1}{m} E_{min} ( m ) $. 
Moreover, equality may occur in one of these inequalities if and only if $E_{min} $ is linear on $(0, m]$. 
Summing up and comparing to (\ref{conv4}) we infer that necessarily we have $ E_{min} ( m_1 ) = \frac{m_1}{m} E_{min} ( m ) $ and 
$E_{min}$ must be linear on $(0, m]$. 
Then Proposition \ref{Emin} (iv) implies that $ E_{min}(\eta) = - \eta $ for any $ \eta \in (0, m]$, 
and in particular $ E_{min}(m) = -m $, contradicting the fact that $ E_{min}(m) < -m$. 
This contradiction shows that necessarily $ m_1 = m $.

Since $ u_n \rightharpoonup u $ weakly in $ L^2 ( \R^N)$ and $ \| u_n \|_{L^2} ^2 \lra m  = \| u \|_{L^2}^2$ we infer that $ u_n \lra u $ strongly in $L^2( \R^N)$. 
Using (\ref{interpol}) and (\ref{GNS}) for $ u_n - u $ we infer that 
$ \nabla u _n \lra \nabla u $ strongly in $ L^2 ( \R^N)$ and 
$ u_n \lra u $ strongly in $ L^{ 2 \si + 2 } ( \R^N)$. 
The weak convergence $ u_n  \rightharpoonup u $  in $ H^2 ( \R^N)$ gives 
$ \| \Delta u \|_{L^2}^2 \les \ds \liminf_{n \ra \infty } \| \Delta u _n \|_{L^2 }^2$,
and consequently we get $ E( u ) \les \ds \liminf_{n \ra \infty } E( u_n ) = E_{min}(m)$. 
On the other hand we have $ E( u ) \ges E_{min}(m)$ because  $ \| u \|_{L^2 }^2 = m $. 
Therefore $ E( u ) = E_{min}(m)$ and $ u $ solves the problem $(\Po _m )$. 
Moreover, we have $ \| \Delta u _n \|_{L^2 }^2 \lra \| \Delta u  \|_{L^2 }^2$. 
Since $ \Delta u _n \rightharpoonup \Delta u $ weakly in $ L^2 ( \R^N)$, we infer that 
$ \Delta u _n \lra\Delta u $ strongly in $ L^2 ( \R^N)$.
The inequality $ \Big\| \frac{ \p ^2 v}{\p x_i \p x_j } \Big\|_{L^2 } \les \| \Delta v \|_{L^2} $ for any $ v \in H^2 ( \R^N)$ 
implies that $ u_n \lra u $ strongly in $H^2 ( \R^N)$ and Theorem \ref{Global} is proven. 
\hfill
$\Box$

\begin{Proposition}
\label{Lagrange}
Assume that $ \si N \les 4$, $ m > 0 $ and $ u \in H^2( \R^N)$ is a solution of the minimization problem $( \Po _m )$. 
Then there exists $ c = c(u ) > 0 $ such that $ u $ satisfies the equation 
\beq
\label{Eq} 
\Delta ^2 u + 2 \Delta u + ( 1 + c ) u - |u |^{ 2 \si } u = 0 \qquad \mbox{ in } H^{-2} ( \R^N). 
\eeq
Moreover, we have: 

\medskip

(i) $ 1 + c \in [ - E_{min, \ell} '( m),    - E_{min, r} '( m) ] $. 

\medskip

(ii) If $ m_0 = 0 $ (where $ m_0 $ is given by (\ref{m0})), then $ c( u ) \lra 0 $ as $ m \lra 0 $. 

\medskip

(iii) If $ \si N < 4 $ we have $ c( u ) \lra \infty $ as $ m \lra \infty$.  

\medskip

(iv) If $ m > m_0 $ 
and
$ E_{min, \ell} '( m) > E_{min, r} '( m) $, there exist at least  two solutions $ u_1 $ and $u_2$ for the problem $(\Po_m)$ such that 
$ 1 + c( u_1) =  - E_{min, \ell} '( m) $ and $ 1 + c( u_2) =  - E_{min, r} '( m) $.

\medskip

(v) If $ m _1 < m_2$, the function $ u_ 1 $ solves $ (\Po _{m_1})$ and $ u_ 2 $ solves $ (\Po _{m_2})$,   then $ c( u_1 ) < c( u_2)$. 

\medskip

(vi) If $ m_0 > 0 $, problem $(\Po_{m})$ does not admit solutions for any $ m\in (0, m _0)$.

\end{Proposition}

\medskip

\noindent
{\it Proof.} 
Since $ E $ and $ M(u ) : = \|  u \|_{L^2}^2 $ are $C^1$ functionals on $ H^2 ( \R^N)$,  the existence of  a Lagrange multiplier 
$ \la _u \in \R$ such that $ E' ( u ) = \la _u M'( u ) $ in $ H^{-2}( \R^N)$ is standard. Then we have
\beq
\label{Eqla}
 \Delta ^2 u + 2 \Delta u - \la _u u - | u |^{ 2 \si } u = 0 \qquad \mbox{ in } H^{-2}( \R^N).
\eeq
We claim that $ \la \in [E_{min, r} '( m), E_{min, \ell} '( m)]$. 
We have $ \| ( 1 \pm t ) u \|_{L^2}^2 = ( 1 \pm t)^2 m$, hence $E(u \pm tu) \ges E_{min}\left( (1 \pm t )^2 m \right)$ and therefore
$$
\begin{array}{l}
2 \la _u m = 2 \la _u \| u \|_{L^2}^2 = \la _u M'(u) .u = E'(u) . u 
= \ds \lim_{t \downarrow 0 } \frac{E( u + tu ) - E(u)}{t} 
\\
\\
\ges \ds \lim_{t \downarrow 0 }  \frac{ E_{min} \left( ( 1 + t )^2 m \right) - E_{min}(m) }{t } = 2 m  E_{min, r} '( m). 
\end{array}
$$
We conclude that  $ \la _u \ges E_{min, r} '( m)$. 
Proceeding similarly with $ 1 - t $ instead of $ 1 + t $ we get $ -\la _u \ges - E_{min, \ell} '( m)$ and the claim is proven.
Denoting $ c(u ) = - \la _u - 1 $, statement (i) follows. 

Taking the $ H^{-2} - H^2 $ duality product of (\ref{Eqla}) and of $u $ we get 
\beq
\label{idla}
\int_{\R^N} |\Delta  u |^2 \, dx - 2 \int_{\R^N} |\nabla u|^2 \, dx - \la _u \int_{\R^N} |u |^2 \, dx - \int_{\R^N} |u|^{ 2 \si + 2} \, dx = 0. 
\eeq
Using (\ref{idla}) and the identities $ \| u \|_{L^2}^2 = m $, $ E(u) = E_{min}(m)$ we get 
\beq
\label{idsimple} 
\| u \|_{L^{2 \si + 2}}^{ 2 \si + 2} 
=  \! \frac{\si + 1}{\si } \left( E_{min}(m) \! - \! \la _u m \right) 
\,   \mbox{ and } 
\int_{\! \R^N} \! |\Delta  u |^2  -2 |\nabla u|^2  dx  = \! \frac{ \si + 1}{\si } E_{min}(m) \! - \!  \frac{\la _u }{\si } m.
\eeq
Since $ E_{min}(m) \les - m $ for all $m > 0 $, the first part in (\ref{idsimple}) implies that for any $ m > 0 $ and for any solution $u$ of 
$(\Po _m)$ we must have $ \la _u < -1$, that is $ c(u) > 0 $. 
If $ m_0 > 0 $,  $ m \in (0, m_0)$ and $ u $ is a solution to the problem $(\Po_m)$ by  (i) we should have $ \la _u = -1$, a contradiction.
Thus (vi) is  proven. 

\medskip

Consider $ m > m_0 $ such that $E_{min}(m ) $ is finite. 
Take an increasing sequence $ (m_n)_{n \ges 2} $ in $(m_0, m )$ such that $ m_n \lra m $. 
For each $ n$, let $u_n$ be a solution of the minimization problem $ (\Po _{m_n} )$ 
(the existence of $u_n$ is guaranteed by Theorem \ref{Global}). 
Then $ (u_n)_{n \ges 2}$ 
satisfies $ \| u _n \|_{L^2} ^2 \lra m $ and $ E_{min}( u_n ) \lra E_{min}(m)$, and 
using Theorem \ref{Global} again we see that there exists a subsequence, still denoted $ (u_n)_{n \ges 2}$, 
and there exists a solution $u_1$ of the problem $ (\Po _{m} )$ such that $ u_n \lra u_1 $ strongly in $ H^2( \R^N)$. 
Identity (\ref{idla}) and the strong convergence in $H^2(\R^N)$ imply that $ \la _{u_1} = \ds \lim_{n \ra \infty} \la_{u_n}.$
On the other hand,  (i) and the basic properties of concave functions imply that 
$ \la_{u_n} \lra E_{min, \ell}'(m)$. Thus we have  $ \la _{u_1} =  E_{min, \ell}'(m)$.
Taking a decreasing sequence $ m_n \lra m $ and proceeding similarly we see that there exists a solution $ u_2$ of $ (\Po _{m} )$
such that $ \la _{u_2} =  E_{min, r}'(m)$.
This proves (iv). 

\medskip

To prove (v) we argue by contradiction and we assume that there are $ m_1 < m_2$ and there are solutions $ u_1$ and $u_2$ of 
$(\Po_{m_1})$ and of $(\Po_{m_2})$, respectively, such that $ c( u_1) = c( u_2)$. 
Since $ -1 - c( u_1) \ges E_{min, r }'( m_1) \ges E_{min, \ell}'(m_2) \ges -1 - c( u_2)$, 
we see that $E_{min, r }'( m_1) = E_{min, \ell}'(m_2)$, and this implies that $ E_{min}$ is affine on $[m_1, m_2]$. 
Hence there exist $ \la < -1 $ and $ B \in \R$ such that $ E_{min}(m) = \la m + B$ for any $ m \in [m_1, m_2]$. 
For any $ m \in (m_1, m_2), $ Theorem \ref{Global} gives the existence of a solution $u$ to the problem $(\Po_m)$ and 
statement (i)  above implies that $ \la _u = \la$, or equivalently $ c(u ) = c(u_1) = c( u_2)$. 
Fix  $ m_3 \in (m_1, m_2)$ and let $u$ be a minimizer for $(\Po_{m_3})$.
The first part of  (\ref{idsimple}) gives $ \|  u \|_{L^{2 \si + 2}}^{ 2 \si + 2} = \frac{ \si + 1}{\si} B$, hence $ B > 0 $.
Using $ \sqrt{t} u $ as test function and taking (\ref{idsimple}) into account we get for $t$ sufficiently close to $1$, 
\beq
\label{concav5}
\begin{array}{l}
\la t m_3 + B  = E_{min} ( tm _3 ) \les E( \sqrt{t} u ) = t \! \int_{\! \R^N} \! |\Delta  u |^2  -2 |\nabla u|^2  dx   
- \frac{t^{\si + 1}}{\si + 1} \! \int_{\R^N} |u |^{ 2 \si + 2} \, dx 
\\
\\
\ds = \la t m_3 + B + \frac{ \si + 1}{\si } B \left( t - \frac{ t^{\si + 1}}{\si + 1} - \frac{ \si }{\si + 1} \right).
\end{array}
\eeq
Since $ B > 0 $ and $  t - \frac{ t^{\si + 1}}{\si + 1} - \frac{ \si }{\si + 1}  < 0 $ for 
$ t \neq 1$, 
  (\ref{concav5}) gives a contradiction. 
This proves (v).

\medskip

All other statements in Proposition \ref{Lagrange} are obvious.
\hfill
$\Box$

\medskip

So far we have solved the global minimization problem $(\Po _m )$  in $H^2( \R^N)$ in the case $ \si N \les 4 $ 
and we have shown that any solution satisfies (\ref{Eq}) for some $ c> 0 $.
Obviously,  $ u \in  H^2( \R^N)$ solves (\ref{Eq}) if and only if $u$ is a critical point of the following functional, called {\it action}:
$$
S_c(u ) = \int_{\R^N} |\Delta  u |^2 \, dx - 2 \int_{\R^N} |\nabla u|^2 \, dx  + ( 1 + c) \int_{\R^N} |u |^2 \, dx 
- \frac{1}{\si + 1}  \int_{\R^N} |u|^{ 2 \si + 2} \, dx.
$$
At this stage it is not clear that given any $ c> 0 $, there exists $m > 0 $ and a solution $u $ of $(\Po _m)$ 
such that $ c(u ) = c$. 
We will show that for any $ c > 0 $   and for any $ \si > 0 $ satisfying $ 2 \si + 2 < 2^{**}$, equation (\ref{Eq}) has solutions and, moreover, 
it has solutions minimizing the action $S_c $ among all solutions (these are called {\it minimum action solutions} or {\it ground states}). 
Moreover, we will show that all minimizers of a problem $(\Po _m)$ are ground states. 
To this end we introduce another family of minimization problems.

Let $ c \ges 0 $. We consider the minimization problem 
$$
\label{Toc}
\begin{array}{c}
\mbox{ minimize } 
\ds {T_c} (u):= \int_{\R^N}|\Delta u|^2\, dx -2 \int_{\R^N}|\nabla u|^2\, dx
+ (1 + c ) \int_{\R^N}|u|^{2}\, dx \\ \\
\mbox{ in the set } 
\ds  U :=\left\{ u\in H^2(\R^N) \; \Big| \; \int_{\R^N}|u|^{2\sigma+2}\, dx=1\right\} .
\end{array}
\leqno{({\To}_c)}
$$
We denote $ t(c) = \ds \inf \{ T_c ( u ) \; | \; u \in U \}$. It is clear that 
\beq
\label{relation}
S_c (u) = T_c (u) - \frac{1}{\si + 1}  \int_{\R^N} |u|^{ 2 \si + 2} \, dx = E(u) + ( 1 + c) \int_{\R^N} | u|^2 \, dx. 
\eeq

\begin{Theorem} 
\label{Tc}
Assume that $ 0 < \si < \infty $ if $ N \les 4 $ and  $ 0 < \si <  \frac{4}{N-4}$ if $ N \ges 5$. 
Then for any $ c > 0 $ we have $ t(c ) > 0 $ and the minimization problem ($\To _c $) admits solutions. 
Moreover, for any sequence  $ ( u_n )_{ n \ges 1 } \subset H^2 ( \R^N)$ satisfying $ \int_{\R^N} |u_n |^{2\sigma+2}\, dx \lra 1 $ and 
$ T_c ( u_n ) \lra t(c)$ there exist a subsequence $( u_{n _k} )_{k \ges 1 }$, a sequence $ ( x_k )_{k \ges 1} \subset \R^N $ and a minimizer $u$ for ($\To _c $) 
such that $ u_{n _k }( \cdot + x_k ) \lra u $ strongly in $ H^2 ( \R^N)$. 

\end{Theorem} 

{\it Proof. } 
The proof is standard, so we only sketch it. 
Fix $ c > 0 $. Then fix $\e > 0 $ such that $ \frac{ 1}{1 - \e} + \e < 1 + c$. 
Using (\ref{Planch}) we get $ T_c ( v ) \ges \e \| \Delta v \|_{L^2}^2 + \e \|  v \|_{L^2}^2 $ for any $ v\in H^2( \R^N)$, and then it is clear that 
$T_c ^{\frac 12}$ is a norm on $ H^2( \R^N)$ and that it is equivalent to the usual norm.
By the Sobolev embedding there exists $ K_c > 0 $ such that $ \| v \|_{L^{ 2 \si + 2}} \les K_c T_c^{\frac 12}(v)$, thus $ t( c ) \ges K_c^{-2} > 0 $. 
For any $ v \in H^2( \R^N)$, $ v \neq 0 $ we have $ \frac{ v}{\| v \|_{L^{2 \si + 2}} } \in U$, hence 
$ T_c \left( \frac{ v}{\| v \|_{L^{2 \si + 2}} } \right) \ges t(c)$ and this gives
\beq
\label{jos}
T_c ( v) \ges t(c ) \| v \|_{L^{2 \si + 2}} ^2 \qquad \mbox{ for any } v \in H^2( \R^N).
\eeq

Let $ ( u_n)_{n \ges 1}$ be a sequence as in Theorem \ref{Tc}. 
It is obvious that $ ( u_n)_{n \ges 1}$ is bounded in $ H^2( \R^N)$. \
We choose $ \si '> \si  $  such that $ 2 \si '+ 2 < 2^{**}$ and we use 
(\ref{nontriv}) for $ u_n $  
to infer that there exists  constants $ t_0, \, a > 0 $, independent of $n$, such that 
$ \Lo ^N \left( \{ | u _n | \ges t_0 \} \right) \ges a $ for all sufficiently large $n$.
Then Lieb's Lemma 
implies that there exists a constant $ b > 0 $, independent of $n$,  and for each $n$ large there exists $ x_n \in \R^N$ such that 
$
 \Lo ^N \left( \left\{ x \in B( x_n, 1) \; \big| \; | u _n | \ges \frac{t_0}{2} \right\} \right) \ges b .
$ 
We replace $ u_n $ by $ u_n ( \cdot +  x_n)$, which is still a minimizing sequence and satisfies 
$ \int_{B( 0, 1 ) } | u _n |^p \, dx \ges b \left( \frac{t_0}{2} \right)^p $ for all $n$.
Since $ ( u_n )_{n \ges 1}$ is bounded in $ H^2( \R^N)$ there is a subsequence, still denoted $ ( u_n )_{n \ges 1}$, 
and there is $ u \in H^2( \R^N)$ such that (\ref{conv}) holds.
The convergence $ u_n \lra u $ in $ L_{loc}^p $ for $ 1 \les p < 2^{**}$ gives 
$ \int_{B( 0, 1 ) } | u |^p \, dx \ges b \left( \frac{t_0}{2} \right)^p $, and therefore $ u \neq 0 $. 
Denote $ \eta = \int_{\R^N} |u |^{ 2 \si + 2} \, dx> 0$. By Fatou's Lemma we get $\eta \les 1$. 
It is obvious that (\ref{conv1}) and (\ref{conv2}) hold.
By (\ref{conv2}) we have $ \int_{\R^N} |u_n - u |^{ 2 \si + 2} \, dx \lra 1 - \eta$ and then using (\ref{conv1}) and (\ref{jos}) we find
$$
T_c ( u_n ) = T_c ( u_n - u ) + T_c ( u) + o(1)
\ges t(c) \| u_n -u \|_{L^{2 \si + 2}} ^2 + t(c ) \| u \|_{L^{2 \si + 2}} ^2 + o(1). 
$$
Letting $ n \lra \infty $ in the above inequality we obtain $ 1 \ges ( 1 - \eta )^{\frac{1}{\si + 1}} +  \eta ^{\frac{1}{\si + 1}}$
and this implies that  $ \eta = 1$, that is $ u \in U$. 
Then we must have $ T_c ( u ) \ges t(c)$. 
On the other hand,  
$ T_c ( u ) \les \ds \liminf_{n \ra \infty } T_c ( u_n) = t(c)$ by weak convergence,  and therefore
$ T_c( u) = t(c) = \ds \lim_{n \ra \infty } T_c ( u_n)$. Since $ T_c$ is a norm on $H^2( \R^N)$ we infer that 
$ u_n \lra u $ strongly in $H^2( \R^N)$, as desired.
\hfill
$\Box$

\medskip

\begin{Proposition}
\label{tc-zero} 
The mapping $ c \longmapsto t(c)$ is strictly increasing on $(0, \infty)$ and 
there is $C>0 $ such that $ t(c ) \les C \sqrt{c}$ for all sufficiently small $c$.
In particular we have
$ t(c) \lra 0 $ as $ c \lra 0 $.
\end{Proposition}

{\it Proof. } Let $ 0 < c_1 < c_2 $. 
Let $ u $ be a minimizer for the problem $(\To _{c_2})$. We have $ u \in U$ and
$ t( c_2 ) = T_{c_2} ( u ) > T_{c_1} ( u ) \ges t( c_1)$. 
Hence the mapping $ c \longmapsto t(c)$ is strictly increasing. 

Let $ u_{\e, \de }$ be as in (\ref{u-epsilon}). 
Fix $ \de_0 = \frac{1}{20}$ and let $ v_c = u_{\sqrt{c}, \de_0}$. 
By (\ref{ex-1}), (\ref{ex-2}) and (\ref{ex-4}) we have
$$
\| v_c\|_{L^2}^2 \les C_1 c^{\frac 12}, \qquad
\| \Delta v_c\|_{L^2}^2 - 2 \|\nabla  v_c\|_{L^2}^2 + \| v_c\|_{L^2}^2 \les C_2 c^{\frac 32}
\quad \mbox{ and } \quad
\| v_c \|_{L^{2 \si + 2}} \ges C_3 c^{\frac 12}
$$
for some $ C_1, C_2, C_3 > 0 $, so that $ T_c( v_c) \les C_4  c^{\frac 32}$. 
Using (\ref{jos}) we see that 
$ t(c) \les \frac{T_c(v_c)}{\| v_c \|_{L^{2 \si +2}}^2} \les C \sqrt{c}$. 
\hfill
$\Box$

\begin{Proposition}
\label{groundstates}
Let $ u $ be any minimizer for the problem $(\To _{c})$.
Then $v:=  t(c)^{\frac{1}{2 \si}} u $ is a solution of (\ref{Eq}). Moreover, for any solution $ w \in H^2( \R^N)$ of (\ref{Eq})
we have $ S_c \left( w \right) \ges S_c( v) = \frac{\si }{\si + 1}  t(c)^{\frac{\si + 1}{ \si}}$. 
In other words, $v$ is a least action solution of (\ref{Eq}).
\\
Conversely, if $\tilde{v}$ is any least action solution of (\ref{Eq}) then $ \tilde{u}:=  t(c)^{- \frac{1}{2 \si}} \tilde{v} $ is a solution of $(\To _{c})$.
\\
Furthermore, if $u$ solves  $(\To _{c})$ there exist $ a \in \R$ and a real-valued function $ \tilde{u} $ such that 
$ u = e^{ia } \tilde{u}$. 

\end{Proposition}

{\it Proof. } 
Assume that $ u $ solves $(\To _c)$.
The functionals $ T_c$ and $ u \longmapsto \| u \|_{L^{2 \si + 2} }^{ 2 \si + 2}$  are $ C^1 $ on $H^2( \R^N)$, 
and consequently there exists a Lagrange multiplier $ \kappa \in \R$ such that 
\beq
\label{Ek}
\Delta ^2 u + 2 \Delta u + ( 1 + c) u = \kappa |u|^{ 2 \si } u \qquad \mbox{ in } H^{ -2}( \R^N).
\eeq
Taking the $H^{-2} - H^2 $ duality product of (\ref{Ek}) with $ u $ we get 
$T_c ( u ) = \kappa \int_{\R^N} | u|^{ 2 \si + 2 } \, dx $, which implies that $ t(c) = \kappa$.
Denoting $ v =  t(c)^{\frac{1}{2 \si}} u  = \kappa^{\frac{1}{2\si}} u$, it is clear  that $v$ solves (\ref{Eq}).
We have 
\beq
\label{value}
 S_c ( v ) = t(c)^{\frac{1}{ \si}} T_c ( u ) 
 - \frac{1 }{\si + 1} t(c) ^{\frac{\si + 1}{\si }}  \int_{\R^N} |u |^{ 2 \si + 2} \, dx 
 = \frac{\si  }{\si + 1} t(c) ^{\frac{\si + 1}{\si }}  .
\eeq

Let $ w \neq 0 $ be an arbitrary solution of (\ref{Eq}). 
The duality product of (\ref{Eq}) with $w$ gives
\beq
\label{facile}
 T_c ( w) = \int_{\R^N} |w |^{ 2\si + 2}\, dx  . 
\eeq
From (\ref{facile}) and (\ref{jos}) we obtain 
$ T_c ( w) \les  t(c )^{ -\si - 1} T_c ( w) ^{ \si + 1} $ and this implies 
$ T_c ( w)  \ges t(c)^{\frac{ \si + 1}{\si }}$.  
(By (\ref{facile}) we have also the lower bound 
$
 \int_{\R^N} |w |^{ 2\si + 2}\, dx =  T_c ( w)  \ges t(c)^{\frac{ \si + 1}{\si }}.) 
$
Using again (\ref{facile})  we find 
$$
S_c ( w) = T_c ( w) - \frac{1}{\si + 1}  \int_{\R^N} |w |^{ 2\si + 2}\, dx = 
\frac{ \si }{\si + 1} T_c(w)  \ges  \frac{ \si }{\si + 1} t(c)^{\frac{ \si + 1}{\si }} = S_c ( v). 
$$

Conversely, let $ \tilde{v}$ be a least action solution  of (\ref{Eq}). 
Then we have $ S_c ( \tilde{v}) =  S_c ({v}) =  \frac{ \si }{\si + 1} t(c)^{\frac{ \si + 1}{\si }}$. 
On the other hand, by (\ref{facile})  we  get $ S_c ( \tilde{v}) =   \frac{ \si }{\si + 1}  T_c( \tilde{v}) 
=  \frac{ \si }{\si + 1} \| \tilde{v} \|_{L^{ 2 \si + 2}}^{ 2\si + 2} $. 
We conclude that $  T_c( \tilde{v})  =  \| \tilde{v} \|_{L^{ 2 \si + 2}}^{ 2\si + 2} = t(c)^{\frac{ \si + 1}{\si }} $
and then one immediately checks that 
$ \tilde{u}=  t(c)^{- \frac{1}{2 \si}} \tilde{v} $ is a minimizer for the problem $(\To _{c})$.

\medskip

Let $ u $ be a solution of  $(\To _{c})$ and let $ u_1$ and $ u_2$ be the real and imaginary parts of $u$, respectively. 
If $ \| u _i \|_{L^{ 2 \si + 2 } } = 0 $ for $ i \in \{ 1, 2 \}$, the last statement in Proposition \ref{groundstates} is obvious. 
Otherwise let $ \al _ i =   \| u _i \|_{L^{ 2 \si + 2 } } $. 
Clearly, we have $ \al _ i \in (0, 1)$ because 
$ \| u  \|_{L^{ 2 \si + 2 } }  =1$,  
and 
$T_c ( u _i ) \geq \al _i ^2 t(c) $ by (\ref{jos}).
Then we get 
$$
t(c ) = T_c ( u ) = T_c ( u_1 ) + T_ c( u_2 ) \ges ( \al _ 1 ^2  + \al _2 ^2 ) t(c), 
$$
hence $ \al _1 ^2 + \al_2 ^2 \les 1 $. 
This gives
$$
\| u _1 ^2 \|_{L^{\si + 1}} + \| u _2 ^2 \|_{L^{\si + 1}} = \al_1 ^2 + \al _2 ^2 \les 1 = \| u \|_{L^{ 2 \si + 2 }}^2 = \| u_1 ^2 + u_2 ^2 \|_{L^{\si + 1}}.
$$
Then Minkowski's inequality implies that $ \al _1 ^2 + \al _2 ^2 = 1 $ and $ u_1 ^2 $ and $ u_2 ^2 $ are proportional. 
Since those functions do not vanish a.e., there exists $ b \in (0, \infty)$ such that $ u_2 = \pm b u_1 $ a.e. on $ \R^N$. 
If $ u_2 = b u_1 $ a.e. we have $ u  = ( 1 + i b ) u _1$, where $ u_1 $ is a real-valued function,  and the conclusion follows. 
Similarly, if $ u_2 = - b u_1 $ a.e. we have   $ u  = ( 1 - i b ) u _1$. 
We argue by contradiction to show that $ u_1 $ and $ u_2 $ are always proportional a.e. on $ \R^N$.
Assuming that this is false, the sets $ A_1 = \{ x \in \R^N \; | \; u_1 ( x ) \neq 0  \mbox{ and } u_2 (x) = b u_1 (x) \} $ and 
 $ A_1 = \{ x \in \R^N \; | \; u_1 ( x ) \neq 0  \mbox{ and } u_2 (x) = - b u_1 (x) \} $ have both positive Lebesgue measure. 
Let $ w = \frac{ 1 - ib}{| 1 - ib |} u $. 
Then $ w $ is also a solution of the minimization problem $( \To_c)$. 
Denoting by $ w_1 $ and $ w_2 $ its real and imaginary parts, respectively, we see as above that $ w_1 ^2 $ and $ w_2 ^2$ must be proportional  on $ \R^N$. 
However, we have 
$ w = \sqrt{ 1 + b^2} u _1 $ on $ A_1 $ and $ w =  \sqrt{ 1 + b^2} u _1 - i \cdot \frac{ 2 b }{\sqrt{ 1 + b^2}} u_1 $ on $ A_2$, a contradiction. 
\hfill
$\Box$

\medskip

We do not know whether minimum action solutions of (\ref{Eq}) are radially symmetric (after a translation in $ \R^N$), 
nor whether real-valued minimum action solutions have constant sign. 
Theorem 4 p. 33 in \cite{BoNa} asserts that in some cases (more precisely, if $N \les 3$ and $ \si = 1$),  if real-valued minimum action solutions are radial then they must  change sign.

\begin{Proposition}
\label{PomToc}
Assume that $u \in H^2( \R^N)$ is a solution of the minimization problem $(\Po _m )$ for some $ m > 0 $ and that $ u $ solves (\ref{Eq}). Then:

\medskip

(i) $u $ is a minimum action solution of (\ref{Eq}). 

\medskip

(ii) If $ v $ is any minimum action solution of (\ref{Eq}) we have $ \| v \|_{L^2}^2 = m $ and $ v $ is a minimizer for $(\Po _m)$.
\end{Proposition}

{\it Proof.}
Since $ u $ solves $(\Po _m )$ and (\ref{Eq}), using Proposition \ref{groundstates} we get 
\beq
\label{Eminu}
E_{min}(m) + ( 1 + c ) m = E(u ) + ( 1 + c ) \| u \|_{L^2}^2 = S_c ( u ) \ges  \frac{\si }{\si + 1}  t(c)^{\frac{\si + 1}{ \si}}.
\eeq
Let $ v $ be an arbitrary  minimum action solution for 
(\ref{Eq}) (the existence of such solutions follows from Theorem \ref{Tc}  and Proposition \ref{groundstates}). 
From the proof of Proposition \ref{groundstates} we know that $ T_c ( v ) = \int_{\R^N} | v |^{ 2 \si + 2} \, dx = t(c)^{ \frac{\si + 1}{\si }}. $
Denote $ m ' = \| v \|_{L^2}^2$. 
For any $ a > 0 $ we have $ \| {a}^{\frac 12}  v \|_{L^2}^2 = a m'$ and taking $ { a}^{\frac 12} v $ as test function we discover
\beq
\label{Eminv}
\begin{array}{l}
E_{min}(a m') + ( 1 + c ) a m' \les E( {a}^{\frac 12} v ) + ( 1 + c ) \|  {a}^{\frac 12}  v \|_{L^2}^2 
\\
\\
\ds = T_c ( {a}^{\frac 12 } v  ) - \frac{1}{\si + 1 } \int_{\R^N} |{a}^{\frac 12}  v |^{2 \si + 2} \, dx 
= \left( a - \frac{ a^{\si + 1}}{\si + 1} \right) t(c)^{ \frac{\si + 1}{\si }} . 
\end{array}
\eeq
The mapping $ a \longmapsto \ph ( a ) := a - \frac{a^{\si + 1}}{\si + 1} $ reaches its maximum value on $( 0 , \infty )$ only at $ a = 1 $ 
and the maximum  is $ \ph (1) = \frac{ \si}{\si + 1}$. 
Comparing (\ref{Eminu})  and   (\ref{Eminv})  we get 
$$
\frac{\si }{\si + 1}  t(c)^{\frac{\si + 1}{ \si}} \les E_{min}(m) + ( 1 + c ) m 
\\
\\
\les \ph \left(\frac{m}{m'} \right) t(c)^{\frac{\si + 1}{ \si}}
\les \frac{\si }{\si + 1} t(c)^{\frac{\si + 1}{ \si}} .
$$
We infer that we must have equality throughout in the above sequence of inequalities. 
Therefore $ S_c ( u ) = \frac{\si }{\si + 1} t(c)^{\frac{\si + 1}{ \si}} $ and $ u $ is a minimum action solution for (\ref{Eq}). 
Moreover, we must have $ m = m'$, that is any minimum action solution $v$ of (\ref{Eq}) satisfies $ \|  v \|_{L^2}^2 = m $. Then we find 
$ E( v ) = S_c ( v ) - ( 1 + c ) \| v \|_{L^2}^2 =  S_c ( u ) - ( 1 + c ) \| u \|_{L^2}^2 = E(u) = E_{min}(m)$, 
and consequently $v$ solves $ (\Po _m)$. 
\hfill
$\Box $

\begin{remark} 
\label{IntId}
\rm (Some integral identities) 
Taking the $H^{-2} - H^2$ duality product of (\ref{Eq}) with $u$ we see that any solution $ u \in H^2( R^N)$ of (\ref{Eq}) satisfies the identity $ N_c(u) = 0$, where 
\beq
\label{Nc}
N_c ( u ) 
= \int_{\R^N} |\Delta u |^2 \, dx - 2  \int_{\R^N} |\nabla u |^2 \, dx 
+ ( 1 + c )  \int_{\R^N} | u |^2 \, dx -  \int_{\R^N} | u |^{2\si + 2 } \, dx .
\eeq

Any solution $ u \in H^2( R^N)$ of (\ref{Eq}) satisfies the identity $ P_c(u) = 0$, where 
\beq
\label{Pc}
\begin{array}{l}
P_c ( u ) = \ds 
\frac{ N-4}{N}  \int_{ \R^N} |\Delta u |^2 \, dx - 2 \frac{N-2}{N} \int_{\R^N} |\nabla u |^2 \, dx 
+ ( 1 + c )  \int_{\R^N} | u |^2 \, dx 
\\
\\
\qquad \qquad \ds -  \frac{1}{\si + 1} \int_{ \R^N} | u |^{2\si + 2 } \, dx .
\end{array}
\eeq
The functionals $N_c $ and $ P_c$ are the Nehari and Pohozaev functionals, respectively, 
and $N_c( u ) =0 $ and $P_c( u ) = 0 $ are the Nehari and Pohozaev (or Derrick-Pohozaev) identities.
 The Pohozaev identity expresses the behaviour of the action functional $ S_c$ with respect to dilations: for any $ u \in H^2( \R^N)$ we have 
$ P_c ( u ) = \frac{ d}{dt }_{\mid _{t = 1}} S_c \left(u  \left( \frac{ \cdot}{t} \right) \right) $, and consequently
one expects  $ P_c ( u ) = 0 $ for any critical point of $S_c$.
To give a formal proof of this fact, one first uses a bootstrap argument to prove some regularity of solutions of (\ref{Eq}) 
($ u \in H^3 ( \R^N) $ is enough). 
Then consider a cut-off function $ \chi \in C_c ^{\infty}( \R^N)$ such that $ \chi = 1 $ on $B( 0, 1)$ and 
$ \mbox{supp} (\chi ) \subset B(0, 2)$, 
take the  $H^{-2} - H^2$ duality product of (\ref{Eq}) with
$ \chi \left( \frac{ \cdot} {n} \right) \sum_{j = 1}^N x_j  \frac{ \p u }{\p x_j } $ and integrate by parts, then let $ n \lra \infty $. 
See Lemma 2.1 in \cite{BCGJ-2} for details. 

Two other functionals are of interest: 
\beq
\label{P1}
P_1 ( u ) = \frac N4 \left( N_c ( u ) - P_c ( u) \right) 
= \int_{\R^N} |\Delta u |^2 \, dx - \int_{\R^N} |\nabla u |^2 \, dx - \frac{ N \si }{4( \si + 1)}  \int_{\R^N} | u |^{2\si + 2 } \, dx
\eeq
and
$$
\begin{array}{l}
P_2( u )   =  \ds \frac{1}{\si } N_c ( u ) - \frac{ \si + 1}{\si } P_c ( u  ) 
\\
\\
 = \ds  \left( \frac{ 4( \si + 1)}{N \si } -1 \right) \int_{\R^N} |\Delta u |^2 \, dx 
- 2 \left( \frac{ 2( \si + 1)}{N \si } -1 \right) \int_{\R^N} |\nabla u |^2 \, dx 
- ( 1 + c )  \int_{\R^N} | u |^2 \, dx .
\end{array}
$$
Obviously, any solution $ u \in H^2( \R^N)$ of (\ref{Eq}) satisfies $ P_1 ( u ) = P_2( u ) = 0 $.
If 
$ u \in H^2( \R^N)$ satisfies $N_c ( u ) = 0 $ and 
$ P_i( u ) = 0 $ for some $ i \in \{ 1, 2\}$, then $ P_c ( u ) = P_1 ( u ) = P_2 ( u ) = 0 $. 

Given $ u  \in H^2( \R^N)$ and $ t > 0 $, we denote
\beq
\label{fiber}
u_t ( x) = t^{\frac N4} u ( t^{\frac 12} x) \qquad \mbox{ and } \qquad u^t ( x ) = t^{\frac{N}{2 \si + 2}} u ( tx).
\eeq
By (\ref{scaling}) we have $ \| u _t \| _{L^2 } = \| u \|_{L^2} $ and $ \| u ^t \| _{L^{2\si + 2} } = \| u \|_{L^{2\si + 2}} $ for all $ t > 0 $. 

One has $  \frac{ d}{dt } \left( S_c \left(u _t \right) \right)= \frac{ d}{dt } \left( E \left(u _t \right) \right) = \frac 2t P_1 ( u_t)$. 
If the mapping $ t \longmapsto E( u_t)$ (or, equivalently, $t \longmapsto S_c( u_t)$) 
achieves a local minimum or a local maximum at $ t =1 $ we must have
 $\frac{ d}{dt }_{\mid _{t = 1}} E \left(u _t \right) = 0 $ and this gives $ P_1 ( u ) = 0 $. 

If  $ t \longmapsto T_c( u^t )$ (or, equivalently, $t \longmapsto S_c( u^t)$) 
achieves a local minimum  at $ t =1 $ we must have
 $\frac{ d}{dt }_{\mid _{t = 1}} T_c \left(u ^t \right) = 0 $ and this gives $ P_2 ( u ) = 0 $.

\end{remark}

It follows from Proposition \ref{tc-zero} that $ t(0) = 0 $ and that the minimization problem $(\To _0)$ has no solution
(any function $ u \in H^2( \R^N) $ satisfying $ T_0 ( u ) = 0 $ must be zero a.e.). 
Therefore one would expect that after renormalization, minimum action solutions of (\ref{Eq}) do not have a meaningful limit as $ c \lra 0 $. 
This is proven in  the next Proposition.

\begin{Proposition}
\label{assympto-zero} 
Assume that $ \si > 0 $ and $ ( N - 4) \si < 4 $.
Let $ u_c$ be any minimum action solution of (\ref{Eq}) and let $ v_c = \frac{ u_c}{\| u_c\|_{L^2}}$, so that $ \| v_c \|_{L^2} =1$. 
We have 
\beq
\label{ass0}
\| \Delta v_c \|_{L^2} \lra 1, \quad 
\| \nabla v_c \|_{L^2} \lra 1, \quad 
\| (\Delta + 1)v_c \|_{L^2} \lra 0 \qquad \mbox{ as } c \lra 0 , 
\eeq
and $ \| v_c \|_{L^p} \lra 0 $ for any $ p \in (2, \infty)$ if $ N \les 4$, respectively  for any $ p \in (2, 2^{**})$ if $ N \geq 5$.

\end{Proposition}

{\it Proof.} 
We claim that for any $ c_0 > 0 $ there exists $ C( N, \si, c _0) > 0 $ such that for   any solution
 $ u \in H^2( \R^N) \setminus \{ 0 \}$  of (\ref{Eq}) with $ 0 < c < c_0 $  we have 
 $ \| \Delta u \|_{L^2} \les  C( N, \si, c _0) \| u \|_{L^2}$. 

To prove the claim we use the identity  $ P_2 ( u ) = 0 $ (which is satisfied by all solutions of (\ref{Eq}) in $H^2( \R^N)$, see Remark \ref{IntId}). 
If $ ( N - 2) \si \ges 2$ we have $  \frac{ 2( \si + 1)}{N \si } -1 \les 0 $ 
and the identity $ P_2 ( u ) = 0  $ implies that 
$   \left( \frac{ 4( \si + 1)}{N \si } -1 \right) \int_{\R^N} |\Delta u |^2 \, dx 
\les ( 1 + c )  \int_{\R^N} | u |^2 \, dx $. 

If $ ( N - 2) \si > 2$, denote $ t =  \frac{\| \Delta u \|_{L^2}}{\| u \|_{L^2}}$.  
We divide $ P_2 ( u ) = 0 $ by $ \| u \|_{L^2}^2$ and we use (\ref{interpol}) to get 
$  \left( \frac{ 4( \si + 1)}{N \si } -1 \right) t^2 - 2 \left( \frac{ 2( \si + 1)}{N \si } -1 \right) t - ( 1 + c) \les 0 $.
This inequality implies that $ t $ remains bounded (the bound can be explicitly computed in terms of $ N$, $\si $ and $c$), and the claim follows.

If $ u $ is a solution of (\ref{Eq})  as above, using (\ref{GNS}) we see that there exists $ \tilde{C}( N, \si , c ) > 0 $ such that $ \| u \|_{L^{2 \si + 2}} \les  \tilde{C}( N, \si , c ) \| u \|_{L^2}$, and 
$ \tilde{C}( N, \si , c ) $ is bounded if $ c $ is bounded.

Now let  $ u_c $ be any minimum action solution of (\ref{Eq}). 
We have seen in the proof of Proposition \ref{groundstates} that 
$ T_c ( u _c ) = \| u _c \|_{L^{ 2 \si +2}}^{ 2 \si + 2} = t(c )^{ \frac{ \si + 1}{\si }}.$
Let $ v_c = \frac{ u _c}{\| u _c \|_{L^2}}$. Using Proposition \ref{tc-zero} we get 
\beq
\label{ass0-1}
T_c ( v_c ) = \frac{T_c ( u_c )}{\| u_c \|_{L^2}^2 } 
= \frac{ \| u _c \|_{L^{ 2 \si +2}}^{ 2 \si + 2}}{\| u_c \|_{L^2}^2 }
\les \tilde{C}( N ,\si, c )^2  \| u _c \|_{L^{ 2 \si +2}}^{ 2 \si}
= \tilde{C}( N, \si, c )^2  t(c ) \lra 0 
\eeq
as $ c \lra 0 $.
Dividing (\ref{Pc}) by $ \| u _c\|_{L^2}^2 $ and proceeding as above we get 
\beq
\label{ass0-2}
\frac{N-4}{N} \| \Delta  v_c \|_{L^2}^2 - 2 \frac{N-2}{N} \| \nabla  v_c \|_{L^2}^2 + ( 1 + c ) \|v_c \|_{L^2}^2 \lra 0 
 \mbox{ as } c \lra 0 .
\eeq
From  (\ref{ass0-1}),  (\ref{ass0-2})
and  the fact that $ \| v _c \|_{L^2} = 1$ we infer that (\ref{ass0})  holds.

We claim that for any sequence $ c_n \lra 0 $ and for any sequence $( x_n )_{n \ges 1} \subset \R^N$, the only possible weak limit  in $H^2( \R^N)$  of $ v_{c_n}(\cdot + x_n)$  is zero. 
Indeed, assume that 
  $v_{c_n}(\cdot + x_n) \rightharpoonup w $ weakly in $ H^2( \R^N)$. 
 Then by weak convergence we have
 $$ \| ( \Delta + 1) w \|_{L^2}\les \ds \liminf_{n \ra \infty } \| ( \Delta + 1) v_{c_n}(\cdot + x_n) \|_{L^2} =  \liminf_{n \ra \infty } \| ( \Delta + 1) v_{c_n}\|_{L^2} = 0, 
 $$
 and then Plancherel's identity  implies $ ( |\cdot | ^2 - 1 ) \wh{w} = 0 $ a.e. in $ \R^N$,
 thus $ w= 0$. 
 The claim is thus proven. 
It implies that $  \ds s(c):= \sup_{y \in \R^N} \int_{B(y, 1 )} |v_c |^2 \, dx \lra 0 $ as $ c \lra 0 $.
For otherwise,  there would exist a sequence of $( c_n)_{n \ges 1}$, a sequence $ (x_n )_{n \ges 1 } \subset \R^N$ and $ \de > 0 $ such that $ \| v_{c_n } \|_{L^2( B( x_n, 1))} \ges \de $ for all $ n $. 
But $(v_{c_n}( \cdot + x_n))_{n \ges 1}$ is bounded in $H^2( \R^N)$, hence it has a subsequence 
$(v_{c_{n_k}}( \cdot + x_{n_k)})_{k \ges 1}$ that converges weakly in $ H^2( \R^N)$ and 
we have seen that its weak limit must be zero.  Since $H^2 ( B(0,1))$  is compactly embedded into $ L^2 ( B(0, 1))$, we get 
$ \| v_{c_{n_k}} \|_{L^2( B( x_{n_k}, 1 ))} \lra 0 $, a contradiction.

Fix $ r_0 \in (2, \infty)$ if $ N \in\{ 1, 2 \}$ and $ r_0 \in (2, \frac{2N}{N-2 }] $ if $ N \ges 3 $. 
Let $ q_0 = 4 - \frac{4}{r_0} \in (2, r_0)$ and $ \al = 1 - \frac{2}{r_0} \in (0, 1)$, so that $ q_0 = 2 \al  + r_0 ( 1 - \al)$
and $ r_0( 1 - \al ) = 2$.
Let $ w \in H^1( \R^N)$.
Let $ B$ be any ball of radius $1$ in $ \R^N$. 
Using H\"older's inequality and the Sobolev embedding we get 
$$
\int_B | w |^{q_0} \, dx \les 
 \| w \|_{L^2( B) }^{ 2 \al } \| w \|_{L^{r_0}(B)}^{(1 - \al ) r_0} 
 \les C  \| w \|_{L^2( B) }^{ 2 \al } \int_{B} | w |^2 + |\nabla w |^2 \, dx . 
 $$
 Let $ (B_n)_{n \ges 1}$ be a collection of balls in $ \R^N$ such that $ \cup_{n \ges 1 } B_n = \R^N$ and each point of  $ \R^N$ belongs to at most $\ell $ of the $B_n$'s, where $ \ell $ is fixed.  
 Writing the above inequality for each $B_n$ and summing up we obtain
\beq
\label{ineq-lions}
\| w \|_{L^{q_0}}^{q_0} \les C \left(\sup _{y \in \R^N} \int_{B(y, 1 ) } |w |^2 \, dx  \right)^{\al } \| w \|_{H^1( \R^N)}^2. 
\eeq

From (\ref{ineq-lions}) we get 
$\| v_c  \|_{L^{q_0}}^{q_0} \les C s(c)^{\al } \| v_c \|_{H^1} \lra 0 $ as $ c \lra 0 $. 
Let $ p $ be as in Proposition \ref{assympto-zero}, $ p \neq q_0$. 
Let $ r \in (p, \infty)$ if $ N \les 4$, respectively   $ r \in (p, 2^{**})$ if $ N \geq 5$, so that $ v_c$ is bounded in $ L^r ( \R^N)$ by the Sobolev embedding.
By H\"older's inequality we have 
$ \| v_c \| _{L^p} \les \| v_c \|_{L^2}^{\al _1} \|v_c \|_{L^{q_0}}^{ 1 - \al _1}$ if $ p < q_0$, respectively 
$ \| v_c \| _{L^p} \les \| v_c \|_{L^r}^{\al _2} \|v_c \|_{L^{q_0}}^{ 1 - \al _2}$ if $ p > q_0$ for some $ \al _1, \al _ 2 \in (0, 1)$, 
hence $  \| v_c \| _{L^p} \lra 0 $ as $ c \lra 0 $.
 \hfill
 $\Box$

 \medskip
 
 From Proposition \ref{P-interesting} (i), Theorem \ref{Global}, Proposition \ref{Lagrange} (ii), Proposition \ref{PomToc} and Proposition \ref{assympto-zero} we get the following result on the behaviour of solutions of the problem $(\Po _m )$ as $ m \lra 0 $.

\begin{Corollary}
\label{asympt-zero}
Assume that $ 0 < \si < \max (1, \frac{4}{N+1}) $ and $ \si \les \frac 4N$. 
For any $m> 0 $ let $u_m$ be any solution of the minimisation problem $(\Po_m)$, as given by Theorem \ref{Global} and Proposition \ref{P-interesting}, 
and let $c_m =  c(u_m)$ be the Lagrange multiplier given by Proposition \ref{Lagrange}, so that $(u_m, c_m)$ solve (\ref{Eq}).
Denote $ v_m = \frac{ v_m }{\sqrt{m}} = \frac{ v_m }{\| u _m \|_{L^2}}$, so that $ \| v _m \|_{L^2} = 1$. 
Then we have 
\beq
\label{ass-zero}
\| \Delta v_m \|_{L^2} \lra 1, \quad 
\| \nabla v_m \|_{L^2} \lra 1, \quad 
\| (\Delta + 1)v_m \|_{L^2} \lra 0 \qquad \mbox{ as } m \lra 0 , 
\eeq
and $ \| v_m \|_{L^p} \lra 0 $ for any $ p \in (2, \infty)$ if $ N \ges 4$, respectively  for any $ p \in (2, 2^{**})$ if $ N \geq 5$.

\end{Corollary}

We will study the behaviour 
of minimum action solutions of (\ref{Eq}) as $ c \lra \infty$. 
To do this we use once again the scaling properties of functionals. 
Given $ c > 0 $ and $ v \in H^2( \R^N)$, we denote
\beq
\label{Kc}
K_c ( v ) = \int_{\R^N} |\Delta v |^2 \, dx - \frac{2}{\sqrt{ 1 + c }}  \int_{\R^N} |\nabla v |^2 \, dx + \int_{\R^N} | v |^2 \, dx 
\eeq
and $ K ( v ) = \int_{\R^N} |\Delta v |^2 \, dx + \int_{\R^N} | v |^2 \, dx $. 
We consider the minimisation problems 
$$
\label{Aoc}
\mbox{minimize } K_c (v) \mbox{ in  } H^2 ( \R^N) \mbox{ under the constraint } \int_{\R^N} | v |^{ 2 \si + 2 } \, dx = 1 , 
\leqno{(\Ao _c)}
$$
$$
\label{Ao}
\mbox{minimize } K(v) 
  \mbox{ in  } H^2 ( \R^N) \mbox{ under the constraint } \int_{\R^N} | v |^{ 2 \si + 2 } \, dx = 1 .
\leqno{(\Ao)}
$$

Let $ c > 0 $. Take $ b = ( 1 + c )^{ - \frac 14 }$. Using (\ref{scaling}) we see that for any $ u \in H^2( \R^N)$ we have
$$
T_c ( u_{a, b } ) = a^2 ( 1 + c )^{ 1 - \frac N4} K_c ( u )  \qquad \mbox{ and } \qquad
\| u _{a, b } \|_{L^{ 2 \si + 2}}^{2 \si + 2 } = a^{ 2 \si + 2} ( 1 + c )^{ - \frac N4}  \| u  \|_{L^{ 2 \si + 2}}^{2 \si + 2 }.
$$
Now choose $ a $ such that $  a^{ 2 \si + 2} ( 1 + c )^{ - \frac N4} = t(c)^{ \frac{ \si + 1}{\si }}$,
 that is $ a = ( 1 + c )^{ \frac{N}{8(\si + 1)}} t(c)^{ \frac{1}{2 \si}}$.
With choice of $a $ and $ b$ and using Proposition \ref{groundstates} we see that $ u_{a, b }$ is a minimum action solution for 
(\ref{Eq}) if and only if $ u $ is a minimizer for $({\Ao_c})$. 
Theorem \ref{Tc} and Proposition \ref{groundstates} give the existence of minimizers for $({\Ao_c})$ for any $ c > 0 $. 

\begin{Lemma}
\label{L3.13} 
(i) Let 
\beq
\label{I}
I = \inf \{ K(u) \; | \; u \in H^2( \R^N) \mbox{ and } \| u \|_{L^{ 2 \si + 2 } } = 1 \}. 
\eeq
Then $ I > 0 $. Moreover, given any sequence $ ( u_n )_{n \ges 1 } \subset H^2( \R^N)$ such that $  \| u_n\|_{L^{ 2 \si + 2 } } \lra 1 $ and $ K(u_n) \lra I $,  there exist a subsequence $ (u_{n_k})_{k \ges 1}$, a sequence $ ( x_k)_{k \ges 1} \subset \R^N $ and a solution $u $ of the minimization problem $(\Ao )$ such that $ u_{n_k } \lra u $ strongly in $ H^2( \R^N)$. 

(ii) Any solution $u$ of the  problem $(\Ao )$ is an optimal function for the inequality (\ref{GNS}).


\end{Lemma} 

{\it Proof. } It is clear that $ K^{\frac 12}$ is a norm on $ H^2( \R^N)$ equivalent to the usual norm. 
By the Sobolev inequality there is $ C_S > 0 $ such that $ \| v \|_{L^{2\si + 2}} \les C_s K^{\frac 12} ( v) $ for any $ v \in H^2( \R^N)$, and this implies that $I > 0 $. 
The proof of part (i) is very similar to the proof of Theorem \ref{Tc} and we omit it. 

(ii) Let $ v \in H^2( \R^N)$. Let $ v^t ( x) = t^{\frac{N}{2 \si + 2}} v( tx)$, as in (\ref{fiber}).
Then we have 
$$
 K( v^t ) = t^{4 - \frac{N\si}{\si + 1}} \int_{\R^N} |\Delta v | ^2 \, dx +  t^{ - \frac{N\si}{\si + 1}} \int_{\R^N} | v | ^2 \, dx .$$
A straightforward computation shows that 
$$ 
\min_{t > 0 } K(v^t ) = C(N, \si ) \| \Delta u \|_{L^2} ^{\frac{N \si}{2( \si + 1)}}  \|  u \|_{L^2} ^{2 - \frac{N \si}{2( \si + 1)}} 
,
\;  \mbox{ where } \;
C(N, \si ) =  \frac{4( \si + 1)}{N \si } \left(  \frac{4( \si + 1)}{N \si }  - 1 \right)^{\frac{ N \si}{4( \si + 1) }-1 }  .
$$
For $ v \in H^2( \R^N) \setminus \{ 0 \}$, let $ \Ro(v ) = \frac{ \| \Delta u \|_{L^2} ^{\frac{N \si}{2( \si + 1)}}  \|  u \|_{L^2} ^{2 - \frac{N \si}{2( \si + 1)}} }{ \| v \|_{L^{2 \si + 2}}^2}.$
If  $ Q $ is a minimizer for the problem $(\Ao)$, then for any $ v \in H^2( \R^N)$, $ v \neq 0 $ we have $ K ( Q) \les K \left( \frac{v}{\| v \|_{L^{2 \si + 2}}} \right)$ and  this gives 
$$
C(N, \si) \Ro ( Q ) = \min _{t > 0 } K( Q^t )  = K(Q) \les  \min _{t > 0 } K \left( \frac{v^t}{\| v \|_{L^{2 \si + 2}}} \right)
= C(N, \si) \Ro ( v ) .
$$
We infer that $ \Ro ( Q ) \les \Ro ( v) $ for any $ v \in H^2( \R^N)$, thus (Q) is an optimal function for (\ref{GNS}). 
We have also $ I = C(N, \si )  B(N, \si)^{-\frac{1}{\si + 1}}$, where $B(N, \si )$ is as in (\ref{maxGNS}).
\hfill
$\Box$


\begin{Proposition}
\label{convergence}

Let $ (c_n)_{n \ges 1} $ be any sequence of positive numbers such that $ c_n \lra \infty$.
Assume that for each $n$, $v_n$ is a minimizer for the problem $(\Ao_{c_n})$. 
There exists a subsequence $( c_{n_k})_{k \ges 1}$, a sequence of points $ ( x_k )_{k \ges 1 } \subset \R^N$ and a minimizer $ v $ for the problem $ ( \Ao )$ such that $ v_{n_k} \lra v $ strongly in $ H^2 ( \R^N)$. 

\end{Proposition}

{\it Proof. } It suffices to show that $( v_n)_{ n \ges 1}$ is a minimizing sequence for the problem $(\Ao )$. 
Then the conclusion of Proposition \ref{convergence} is a consequence of  Lemma \ref{L3.13} (i). 

Since $ \| v_n \|_{L^{2 \si + 2}} = 1 $ for any $ n$, all we have to do is to show that $ K( v_n ) \lra I$ as $ n \lra \infty $. 

From (\ref{Planch}) we have
$$
K_c( u ) < \left( 1 - \frac{1}{\sqrt{1 + c }} \right) K(u) \qquad \mbox{ for any } u \in H^2 ( \R^N), \, u \neq 0. 
$$
Let $Q$ be a minimizer for the problem $(\Ao )$ and let $ v_c $ be a minimizer for $(\Ao _c)$.
 Taking $Q$ as test function in $(\Ao _{c} )$ we get $K_c( v_c) \les K_c(Q )$ and taking $ v_c$ as test function in $(\Ao )$ we obtain 
$K(Q) \les K( v_c)$, hence
\beq
\label{comparison}
\left( 1 - \frac{1}{\sqrt{1 + c }} \right) K( Q )
\les \left( 1 - \frac{1}{\sqrt{1 + c }} \right) K(v_c) < K_{c} ( v_c ) \les K_{c } (Q) < K(Q)  . 
\eeq
Using (\ref{comparison}) we infer that $K( v_n)$ is bounded, thus $ (v_n)_{n \ges 1}$ is bounded in $ H^2( \R^N)$. 
Moreover, the above inequality implies that $ \ds \lim_{n \ra \infty } K( v_n ) = K(Q) = I$ 
 and the conclusion of Proposition \ref{convergence} follows. 
\hfill
$\Box $

\medskip

\begin{Corollary}
\label{estimates}
Let $ I$ be as in (\ref{I}). For any $ c >0 $ we have 
\beq
\label{tc-estimate}
\left( 1 - \frac{1}{\sqrt{1 + c }} \right) ( 1 + c) ^{ 1 -  \frac{ N \si}{4 ( \si + 1)} } I < t(c) <  ( 1 + c) ^{ 1 -  \frac{ N \si}{4 ( \si + 1)} } I.
\eeq
Moreover, if $u_c$ is any minimum action solution of (\ref{Eq}) we have
\beq
\label{converge-m}
( 1 + c )^{ \frac N4 - \frac{1}{\si} } \int_{\R^N} | u_c |^2 \, dx \lra  \frac{4(\si + 1) - N \si }{4( \si + 1)}  I^{\frac{\si + 1}{\si } }, 
\eeq
\beq
\label{converge-d}
( 1 + c )^{ \frac N4 - \frac{1}{\si} -1 } \int_{\R^N} | \Delta u_c |^2 \, dx \lra  \frac{ N \si }{4( \si + 1)}  I^{\frac{\si + 1}{\si } }
\qquad \mbox{ and } 
\eeq
\beq
\label{converge-si}
( 1 + c )^{ \frac N4 - \frac{1}{\si} -1 } \int_{\R^N} |  u_c |^{2\si + 2}  \, dx \lra   I^{\frac{\si + 1}{\si } }  \qquad \mbox{ as } c \lra \infty.
\eeq

\end{Corollary}

{\it Proof. } 
If $Q$ is a minimizer for $(\Ao )$ and $ Q^t ( x ) = t^{\frac{N}{2 \si + 2}} Q ( tx) $ is as in (\ref{fiber}), 
the mapping $ t \longmapsto K( Q^t)$ achieves its minimum on $(0, \infty )$ at $ t = 1$, hence $ \frac{t}{dt}_{ |t = 1} K( Q^t) = 0 $ and this gives 
$$
\left( 4 - \frac{ N \si}{\si + 1} \right) \int_{\R^N } | \Delta Q  | ^2 \, dx - \frac{ N \si}{\si + 1} \int_{\R^N } |  Q  | ^2 \, dx = 0.
$$
From this identity  and the fact that $ K(Q) = I$ we get 
$$
\int_{\R^N } | \Delta Q  | ^2 \, dx = \frac{N \si }{4( \si + 1)} I \qquad \mbox{ and } \qquad 
\int_{\R^N } |  Q  | ^2 \, dx = \frac{4(\si + 1) - N \si }{4( \si + 1)} I .
$$
Notice that the above identities hold for {\it any} minimizer of $(\Ao )$. 
For $ c > 0 $, let $ v_c $ be any  minimizer for the problem $(\Ao _c)$. 
Then Proposition \ref{convergence} and the previous identities imply that 
\beq
\label{converge-v}
\int_{ \R^N } \! | \Delta v_c  | ^2 \, dx \lra  \frac{N \si }{4( \si + 1)} I \quad \mbox{ and } \quad 
\int_{ \R^N } \! |  v_c  | ^2 \, dx = \frac{4(\si + 1) - N \si }{4( \si + 1)} I
\quad \mbox{ as } c \lra \infty. 
\eeq
Given $c> 0 $, let $ u_c $ be a minimum action solution of (\ref{Eq}). 
Let  $ a = ( 1 + c )^{ \frac{N}{8(\si + 1)}} t(c)^{ \frac{1}{2 \si}}$, $ b =  ( 1 + c )^{ - \frac 14 }$, and let 
$ v_c = ( u_c)_{a^{-1}, b^{-1} } = \frac 1a u_c (\b \cdot)$. 
We have already seen that $v_c $ is a minimizer for problem $(\Ao _c)$.
We have $ u_c = (v_c)_{a, b } $ and 
$$
 t(c )^{ \frac{ \si + 1}{\si }} = T_c ( u_c ) = a^2 ( 1 + c)^{ 1 - \frac N4 } K_c ( v_c)
= ( 1 + c ) ^{ 1 - \frac{ N \si}{4 ( \si + 1)} } t(c)^{\frac{1}{\si }} K_c ( v_c) . 
$$
From the above equality and (\ref{comparison}) we get (\ref{tc-estimate}).
We have also 
$$
\int_{\R^N} | u_c |^2 \, dx = a^2 b^N \int_{\R^N} | v_c|^2 \, dx = ( 1 + c ) ^{ -\frac{N \si }{4( \si + 1)}} t(c) ^{\frac{1}{\si}}  \int_{\R^N} | v_c|^2 \, dx \qquad \mbox{ and }
$$
$$
\int_{\R^N} | \Delta u_c |^2 \, dx = a^2 b^{N-4} \int_{\R^N} |\Delta v_c|^2 \, dx 
= ( 1 + c ) ^{ 1 - \frac{N\si }{4( \si + 1)} } t(c) ^{\frac{1}{\si}}  \int_{\R^N} | \Delta v_c|^2 \, dx .
$$
Then taking into account (\ref{converge-v}) we obtain (\ref{converge-m}) and (\ref{converge-d}).
Recall that $  \int_{\R^N} |  u_c |^{2\si + 2}  \, dx  = t(c)^{\frac{\si + 1}{\si }} $, and consequently (\ref{converge-si})
follows from (\ref{tc-estimate}).
\hfill
$\Box$

\begin{remark} 
\label{massasym}
\rm
We have $ 1 + \frac{1}{\si } - \frac N4 > 0 $ because $ 2 + 2 \si < 2^{**}$, and (\ref{converge-d}) implies that we have always 
$\| \Delta u _c \|_{L^2} \lra \infty $ as $ c \lra \infty$. 
On the contrary, from  (\ref{converge-m}) we see that 
$ \| u_c \|_{L^2} \lra \infty $ if $ N \si < 4 $ and $ \| u_c \|_{L^2} \lra 0 $ if $ N \si > 4 $.
In the case $ N \si = 4 $,  (\ref{converge-m})  implies that 
$ \| u_c \|_{L^2} \lra ( \si + 1) ^{\frac{1}{\si}} B(N, \si ) ^{-\frac{1}{\si}} = k_*$, 
where $ k_*$ is as in Proposition \ref{Emin} (vi).

\end{remark}

\begin{remark} 
\label{col}
\rm
For any $ \si > 0 $ such that $ 2\si + 2 < 2^{**}$, the functional $S_c$ has a mountain-pass geometry.
Indeed, we have 
$$
S_c ( u ) = T_c ( u ) - \frac{1}{\si + 1} \| u \|_{L^{2 \si + 2 }}^{ 2 \si + 2} \ges T_c ( u ) - \left(\frac{ T_c(u)}{t(c) } \right)^{\si + 1 } . 
$$
The mapping $ \ph (t) : = t - \frac{1}{\si + 1} \left( \frac{t}{t(c)} \right)^{ \si + 1 }$ is increasing on 
$[0, t(c) ^{\frac{\si + 1}{\si }} ]$,  decreasing on $[ t(c) ^{\frac{\si + 1}{\si }}, \infty )$, 
and $ \ph \left( t(c) ^{\frac{\si + 1}{\si }} \right) = \frac{\si }{\si + 1}  t(c) ^{\frac{\si + 1}{\si }} > 0 $. 
Denoting $ B_c := \{ u \in H^2 ( \R^N) \; | \; T_c ( u ) <  t(c) ^{\frac{\si + 1}{\si }} \}$, we have:

$\bullet \; $ $ S_c ( u ) \ges \ph \left( T_c ( u ) \right) \ges 0 $ for any $ u \in B_c $ and $ \ds \inf_{u \in B_c } S_c ( u ) = S_c ( 0 ) = 0 .$

$\bullet \; $ $ \inf \{ S_c ( u ) \; | \; u \in H^2 ( \R^N) \mbox{ and } T_c ( u ) = t(c) ^{\frac{\si + 1}{\si }} \} =  \frac{\si }{\si + 1}  t(c) ^{\frac{\si + 1}{\si }} > 0 .$

$\bullet \; $  $ \ds \lim_{ t \ra \infty } S_c ( tu ) = - \infty $ for any $ u \neq 0 $. 

Let $ \Gamma := \{ \g : [0, 1] \lra H^2 ( \R^N) \; | \; \g \mbox{ is continuous, } \g ( 0 ) = 0 \mbox{ and } S_c ( \g ( 1) ) < 0 \}. $
It is obvious that for any $ \g \in \Gamma $ there exists $ s \in (0, 1) $ such that $ T_c ( \g ( s )) =  t(c) ^{\frac{\si + 1}{\si }}  $ and therefore 
$$
i_c := \inf_{\g \in \Gamma } \Big( \sup _{s \in [0, 1]} S_c ( \g ( s)) \Big) 
\ges \ph \left (  t(c) ^{\frac{\si + 1}{\si }}  \right) = \frac{\si }{\si + 1}  t(c) ^{\frac{\si + 1}{\si }}  > 0 .
$$
On the other hand, let $ u_c $ be a minimum action solution of (\ref{Eq}), as given by Theorem \ref{Tc} and Proposition \ref{groundstates}. 
We have 
$$
S_c ( \tau^{\frac 12} u_c ) = \tau  t(c) ^{\frac{\si + 1}{\si }} - \frac{\tau^{\si + 1}}{\si + 1}  t(c) ^{\frac{\si + 1}{\si }} 
$$
and we see that $ \ds \max _{\tau > 0 } S_c ( \tau^{\frac 12 } u_c ) = S_c ( u_c ) = \frac{\si }{\si + 1 }  t(c) ^{\frac{\si + 1}{\si }} $. 
We conclude that necessarily $ i_c = \frac{\si }{\si + 1 }  t(c) ^{\frac{\si + 1}{\si }} $, that for $ a > 0 $ sufficiently large the mapping 
$ \tau \longmapsto \tau^{\frac 12 } a u_c  $ is an optimal path in $ \Gamma, $ and that 
$u_c$ is a "mountain-pass solution" of (\ref{Eq}). 

Conversely, if $u $ is any critical point of $ S_c $ at the mountain-pass level $i_c $ (that is, $S_c ( u ) = i_c$), 
by Proposition \ref{groundstates} we know that $ u $ is a minimum action solution of (\ref{Eq}).

\end{remark}

\section{Local minimization in the case $ N \si > 4 $ }
\label{Local}

Throughout this section we assume that $ \si > \frac 4N$   and $ 2 \si + 2 < 2^{**}$, that is 
$ \si < \infty $ if $ N \les 4 $ and $ \si < \frac{4}{N-4}$ if $ N \ges 5$. 
By Proposition \ref{Emin} (i) we have $ E_{min}(m) = - \infty$ for any $ m > 0 $. 
We will investigate the existence of {\it local } minimizers of $E$ when the $L^2-$norm is kept fixed. 
By {\it local minimizer } we mean a function $ u \in H^2( \R^N)$ such that there exists an open set $ \Uo \subset H^2( \R^N)$ 
such that $ u \in \Uo $ and $ E(u ) = \inf \{ E(v) \; | \; v \in \Uo \mbox{ and } \| v \|_{L^2} = \| u \|_{L^2}\}$.

For any $u \in H^2( \R^N)$ let $ u_t ( x) = t^{\frac N4} u ( t^{\frac 12} x) $ be as in (\ref{fiber}). We denote
$$
\ph _ u ( t ) = E( u_t) = t^2  \int_{\R^N} |\Delta u |^2 \, dx - 2t  \int_{\R^N} |\nabla u |^2 \, dx - \frac{ t^{\frac{N \si }{2}}}{ \si + 1}  \int_{\R^N} | u |^{2\si + 2 } \, dx, 
$$
and 
\beq
\label{defD}
D(u ) =  \int_{\R^N} |\Delta u |^2 \, dx - \frac{ N \si}{4(\si + 1)} \left( \frac{N\si}{2} -1 \right) \int_{\R^N} | u |^{2\si + 2 } \, dx.
\eeq
The behaviour of the function $ \ph _u$  inspired the local minimization approach developed below. 
For later use we state here the following elementary lemma.

\begin{Lemma}
\label{elem}

Let $ a, \, b, \, c > 0 $ and define $ f : [0, \infty ) \lra \R$ by $ f ( t ) = at ^2 - 2 bt - c t^{ \frac{N \si}{2}}. $
We have:

\medskip

(i) The second derivative $ f  '' $ is decreasing. 
There exists a unique $t_{infl} > 0 $ such that $ f '' (t_{infl}) = 0 $, and it is given by 
$ t_{infl} = \left( \frac{8a}{N\si (N\si - 2)c}  \right)^{\frac{2}{N\si - 4}}.$

\medskip

(ii) The derivative $ f  ' $ is increasing on $[0, t_{infl}] $ and decreasing on $[t_{infl}, \infty)$, 
and we have $ f  ' (t_{infl} ) > 0 $ if and only if   
$
{a^{1 - \frac{N\si}{2} }}{b^{\frac{N\si}{2} -2} c}< \frac{8}{N\si( N \si -2 )} \left( \frac{N \si - 4}{N \si - 2} \right)^{\frac{N \si}{2} -2 }.
$

\medskip

For the next statements we assume that $ f  ' (t_{infl} ) > 0 $.

\medskip

(iii) There exist a unique $ t_1 \in (0, t_{infl}) $ and a unique $ t_2 \in (t_{infl}, \infty)$  such that 
$ f ' (t_1 ) = 0 $ and  $ f' (t_2 ) = 0 $.
The map $ f $ is decreasing on $[0, t_1]$, increasing on $[t_1, t_2]$, decreasing on $[t_2, \infty )$ 
and reaches its minimum value on $[0, t_2]$ at $ t_1$.

\medskip

(iv)  For $  t_2 \les t' < t '' $ we have $ f ( t'') - f ( t ') \les \frac 12 ( t '' - t')^2 f''( t_2)$.

\medskip

(v)  We have $f( t_{infl}) - f( t_1) = h \left( \frac{t_{infl}}{t_1} \right) t_1 ^{\frac{N\si}{2}} c$, where
$$
\begin{array}{l}
h(s) = \frac 12 \left(\frac{N\si}{2} + 1 \right) \left(\frac{N\si}{2} - 2 \right) s^{\frac{N \si}{2}} 
- \frac{N\si}{2} \left( \frac{N\si}{2} -1 \right) s^{\frac{N\si}{2} -1} 
+ \frac{N\si}{4} \left( \frac{N\si}{2} -1 \right) s^{\frac{N\si}{2} -2}
+ \frac{N\si}{2} (s-1) + 1.
\end{array}
$$
The function $h$ satisfies  $ h(1 ) = h'(1) = h''(1) = 0 $ and
$$
\begin{array}{l}
h''(s) = \frac{N\si}{4} \left(\frac{N\si}{2} -1\right) \left(\frac{N\si}{2} -2 \right) s^{\frac{N\si}{2} - 4} (s-1) \left[ \left( \frac{N\si}{2} +1 \right) s - \left( \frac{N\si}{2} -3 \right) \right], 
\end{array}
$$
thus $h$ is positive, increasing and convex on $(1, \infty)$. 

\end{Lemma}

{\it Proof. } 
This is simple Calculus. We have 
$$
f'(t ) = 2at - 2b - \frac{N \si}{2} c t^{\frac{N \si}{2} -1} 
\qquad \mbox{ and } \qquad
f''( t ) = 2a -  \frac{N \si}{2} \left(  \frac{N \si}{2} -1 \right) c t^{\frac{N \si}{2} -2}   .
$$
Statements (i), (ii), (iii) are obvious.
For (iv) we use the fact that $f''$ is decreasing on $ [0, \infty)$ and $ f'< 0 $ on $ (t_2, \infty)$. We  have:
$$
f ( t'') - f( t') 
= \int_{t'}^{t''} \left(f'(t') + \int_{t'} ^s f''( \tau ) \, d \tau \right) \, ds
\les   \int_{t'}^{t''} \int_{t'} ^s f''( t_2 ) \, d \tau  \, ds
= \frac 12 ( t'' - t')^2 f''( t_2).
$$

(vi)  Let $ s = \frac{t_{infl}}{t_1}$. Recall that $ s > 1$ because $ t_1 < t_{infl}$. 
From the identity $ f''( t_{infl} ) = f''( t_1 s) = 0 $ we get 
$ a = \frac{N \si}{4} \left(  \frac{N \si}{2} -1 \right) c t_1 ^{\frac{N \si}{2} -2}  s  ^{\frac{N \si}{2} -2}. $
Replacing this into the identity $ f'( t_1 ) = 0 $ we obtain 
$ b = \frac{N \si}{4} c   t_1 ^{\frac{N \si}{2} -1} \left[ \left(  \frac{N \si}{2} -1 \right)   s  ^{\frac{N \si}{2} -2} -1 \right].$
Replacing these values of $ a $ and $b$ into $ f ( t_1 s ) - f( t_1)$
we get the announced identity. 
The properties of the function $h$ are obtained by direct computation. 
\hfill
$\Box$

\medskip

Recall that the functional $P_1$ has been introduced in (\ref{P1}). We have 
\beq
\label{deriv}
\begin{array}{l}
\ph _u '(t )  = 2t  \!  \int_{\! \R^N} |\Delta u |^2 \, dx  -  2 \! \int_{\! \R^N} |\nabla u |^2 \, dx 
- \ds \frac{ N \si t^{\frac{N \si }{2}-1}}{2( \si + 1)}  \! \! \int_{\! \R^N} | u |^{2\si + 2 } \, dx 
= \ds \frac{ 2 P_1( u_t)}{t}
\; \mbox{ and }
\\
\\
\ph_u '' (t) 
=  2 \int_{\R^N} |\Delta u |^2 \, dx - \ds \frac{ N \si t^{ \frac{N\si}{2 } -2}}{2(\si + 1)} \left( \frac{N\si}{2} -1 \right) \int_{\R^N} | u |^{2\si + 2 } \, dx  =\ds  \frac{2}{t^2} D( u_t) .
\end{array}
\eeq
For any $ u \neq 0 $ there exists a unique $ t_{u, infl} > 0 $ such that $ \ph_u '' (t_{u, infl}) =0 $. It is given by 
\beq
\label{tinfl}
t_{u, infl} = \left( \frac{8(\si + 1)}{N \si ( N\si - 2) }  \int_{\R^N} |\Delta u |^2 \, dx \right)^{\frac{2}{N \si -4 }} \left( \int_{\R^N} | u |^{2\si + 2 } \, dx \right)^{-\frac{2}{N \si -4 }}. 
\eeq
We have $ \ph_u '' > 0 $ on $( 0 , t_{u, infl}) $ and $ \ph_u '' < 0 $ on $( t_{u, infl}, \infty)$, 
hence $ \ph_u'$ is increasing on $( 0 , t_{u, infl}] $ and decreasing on  $[ t_{u, infl}, \infty)$, 
therefore reaches its maximum value at $  t_{u, infl}$. 
If $ \ph_u '(t_{u, infl}) \les 0 $, the mapping $ \ph _u $ is (strictly) decreasing on $(0, \infty)$ and consequently none of the functions 
$(u_t)_{t >0}$ can be a local minimizer of $E$ when the $L^2-$norm is kept fixed.
If $  \ph_u '(t_{u, infl}) > 0 $, it is easily seen that there exist a unique $ t_{u, 1} \in ( 0 , t_{u, infl}) $ and a unique 
$t_{u, 2} \in ( t_{u, infl}, \infty)$ such that $ \ph_u'( t_{u, 1}) = \ph_u'( t_{u, 2}) = 0$.
We have 
$ \ph _u ' < 0 $ on $(0, t_{u, 1}) \cup (t_{u, 2}, \infty)$ and $ \ph _u ' > 0 $ on $ ( t_{u, 1}, t_{u, 2})$, therefore
$ \ph _u $ is decreasing on $(0, t_{u, 1}]$, 
increasing on $[t_{u, 1}, t_{u, 2} ]$ and decreasing on $[ t_{u, 2}, \infty)$. 
It is now clear that among the functions $(u_t)_{t >0}$, the only one that could eventually be a local minimizer of $E$ when the $L^2-$norm is fixed is $ u_{t_{u, 1}}$. 
If $ u $ is  a local minimizer of $E$ at constant $L^2-$norm, we must have $ t_{u, 1} = 1$ and $ 1 < t_{u, infl} < t_{u, 2}$, thus necessarily $ D(u ) > 0 $. 
The above discussion indicates that it is natural  to look for local minimizers of $E$ at fixed $L^2-$norm in the set
\beq
\label{O}
\begin{array}{rcl}
\Oo & = & \{ u \in H^2( \R^N) \; | \; u \neq 0, \; t_{u, infl} > 1 \mbox{ and } \ph_u' (t_{u, infl}) > 0 \}
\\
\\
& = & 
\{ u \in H^2( \R^N) \; | \; u \neq 0, \; D( u ) > 0 \mbox{ and } \ph_u' (t_{u, infl}) > 0 \}. 
\end{array}
\eeq
It is clear that $ u \longmapsto t_{u, infl} $ and $  u \longmapsto P_1( u _{t_{u, infl}} ) $
are continuous on $ H^2(\R^N) \setminus \{ 0 \}$ (see (\ref{tinfl})), hence $\Oo$ is  open. 
Given any $ u \in H^2( \R^N)\setminus \{ 0 \}$, using Lemma \ref{elem} (ii) 
we see that $  \ph_u' (t_{u, infl}) > 0$ if and only if 
\beq
\label{cond-infl}
H(u):= \frac{ \left(  \int_{\R^N} |\Delta u |^2 \, dx  \right)^{\frac{ N \si }{2} - 1}}
{ \int_{\! \R^N} | u |^{2\si + 2} \, dx \cdot \left( \int_{\! \R^N} |\nabla u |^2 \, dx  \right)^{\frac{ \! N \si }{2} - 2}}
> \frac{ N \si}{4( \si + 1)} \! \left( \! \frac{ N \si }{2} - 1 \! \right) \! \left( \! \frac{ N \si - 2}{N \si - 4} \right)^{\! \! \frac{ N \si }{2} - 2}.
\eeq

Using (\ref{interpol}) (with strict inequality because $u \neq 0 $) and (\ref{GNS}) we have 
$$
 \int_{\R^N} | u |^{2\si + 2} \, dx \cdot \left( \int_{\R^N} |\nabla u |^2 \, dx  \right)^{\frac{ N \si }{2} - 2}
< B(N, \si ) \| \Delta u \|_{L^2} ^{ N \si - 2} \| u \|_{L^2}^{ 2 \si }. 
$$
Therefore 
$H(u ) >  \frac{1}{B(N, \si) \| u \|_{L^2}^{ 2 \si }}$. 
Denote 
\beq
\label{m0}
\mu _0 = B(N, \si ) ^{- \frac{1}{\si } } \left[  \frac{ N \si}{4( \si + 1)} \left( \frac{ N \si }{2} - 1 \right) \right] ^{- \frac{1}{\si }}
\left( \frac{ N \si - 2}{N \si - 4} \right)^{\frac{2}{\si}  - \frac{ N  }{2} }.
\eeq
We infer that (\ref{cond-infl}) holds for any $ u \in H^2( \R^N) \setminus \{ 0 \}$ satisfying $ \| u \|_{L^2}^2 \les \mu_0 $. 

Using (\ref{GNS}) we see that 
$ D(u) \ges \| \Delta u \|_{L^2}^2 - \frac{ N \si}{4(\si + 1)} \left( \frac{N\si}{2} -1 \right) B(N, \si) \| \Delta u \|_{L^2}^{\frac{N\si}{2}} \| u \|_{L^2}^{ 2 \si + 2 - \frac{N\si}{2}}$ for any $ u $, hence $ D(u) > 0  $ if $ u \neq 0 $ and
$ 1 >  \frac{ N \si}{4(\si + 1)} \left( \frac{N\si}{2} -1 \right) B(N, \si) \| \Delta u \|_{L^2}^{\frac{N\si}{2}-2} \| u \|_{L^2}^{ 2 \si + 2 - \frac{N\si}{2}}$.
Let 
$$
\begin{array}{l}
\Oo_1 = \left\{ u \in H^2( \R^N ) \; \Big| \; u \neq 0, \; 
\| \Delta u \|_{L^2}^{\frac{N\si}{2}-2 } \| u \|_{L^2}^{ 2 \si + 2 - \frac{N\si}{2}} < \frac{8(\si + 1)}{N \si(N \si -2 ) B(N, \si)} \mbox{ and } \| u \|_{L^2}^2 < \mu_0 \right\}. 
\end{array}
$$
Obviously, $ \Oo _1 \cup \{ 0 \} $ is an open neighbourhood of $0$ in $H^2( \R^N)$  
and  $ \Oo _1 \subset \Oo$. 
It follows immediately from the definition of $D$ and from (\ref{cond-infl}) that $ \Oo \cup\{ 0 \}$ is "star-shaped": for all $ u \in \Oo$ and for all $ a \in (0, 1)$ we have $ au \in \Oo.$

\medskip

For any $ m > 0 $ we denote 
\beq
\label{tildE}
\tilde{E}_{min} ( m ) = \inf \{ E(u) \; | \; u \in \Oo \; \mbox{ and }  \| u \|_{L^2}^2 = m \}. 
\eeq

It is obvious that $ (u_t)_s = u_{ts}$. If $ u \in \Oo$ and $ t > 0 $ we have $ u_t \in \Oo $ if and only if $ t < t_{u, infl}$, 
and $ t _{u_t, infl} = \frac{t_{u, infl}}{t}$, $ t_{u_t, i } = \frac{ t_{u, i}}{t}$ for $ i = 1, 2$. 
If $ u \in \Oo$  satisfies  $  \| u \|_{L^2}^2 = m $, 
the previous discussion shows that $ \min \{ E( u_t ) \; | \; 0 < t < t_{u, infl} \} = E( u_{t_{u, 1}}) $ 
and $  u_{t_{u, 1}} $ is the only function among $ ( u_t)_{ 0 < t < t_{u, infl}} $ where $ P_1$ vanishes. 
We have thus proved that 
\beq
\label{tildE-bis}
\tilde{E}_{min} ( m ) = \inf \{ E(u) \; | \; u \in \Oo \; \mbox{ and }  \| u \|_{L^2}^2 = m \; \mbox{ and }   P_1(u) = 0 \}. 
\eeq

\begin{remark}
\label{no-local-min}
\rm 
If $ \si > \frac 4N$ and $E$ is as in (\ref{sign-energy}) with $ \epsilon \les 0 $, there do not exist non-trivial minimizers of $E$ at fixed $L^2-$norm. 
Indeed, let $ u \in H^2( \R^N) \setminus \{ 0 \}$, let $ u_t = t ^{\frac N4} u ( t^{\frac N2} \cdot)$, as in (\ref{fiber}), and let $ \ph _u ( t ) = E( u_t)$ as above. 
There exists a unique $ t_{u, infl} > 0 $ such that $ \ph_u ''( t_{u, infl}) = 0 $ and it is given by (\ref{tinfl}). 
We have $ \ph_u '' > 0 $ on $ (0, t_{u, infl})$ and $ \ph_u '' < 0 $ on $(t_{u, infl}, \infty)$.
There exists  a unique $ t _u > 0 $ such that $ \ph_u '(t_u) = 0 $ and we have $ t_u > t_{u, infl}$, 
$ \ph_u ' > 0 $ on $ (0, t_{u})$ and $ \ph_u ' < 0 $ on $(t_{u}, \infty)$.
Therefore $ \ph_u $ is increasing on $ (0, t_{u})$, decreasing on $(t_{u}, \infty)$,  it achieves its global maximum at $ t = t_u$ and it has no local minimum on $(0, \infty)$. 
The previous discussion shows that no function $ u  \in H^2( \R^N) \setminus \{ 0 \}$ can be a local minimizer of the energy at fixed mass.

\end{remark}

\begin{Lemma}
\label{basic}
The following assertions hold true:

\medskip

(i) For any $ m >0$, the set $ \{ u \in \Oo \; | \; \| u \|_{L^2} ^2 = m \mbox{ and } P_1( u ) = 0 \}$ is not empty  (thus $\tilde{E}_{min} (m) < \infty$),
 and $\tilde{E}_{min} (m) \ges - \frac{ (N \si -2)^2}{N\si ( N \si - 4)}m .$

(ii) For all  $ m > 0 $ and all $d,  e \in \R$, the set 
$$ \left\{ u \in H^2( \R^N ) \; | \; D( u ) \ges d, \; \| u \|_{L^2}^2 \les m \; \mbox{ and } E(u ) \les e \right\}$$
 is bounded in $ H^2 ( \R^N)$.

\medskip 

(iii) $\tilde{E}_{min} $ is sub-additive:   $ \tilde{E}_{min} (m_1 + m_2) \les \tilde{E}_{min} (m_1) + \tilde{E}_{min} (m_2)$ for any $ m_1, m_2 > 0 $.

\medskip

(iv) $ \tilde{E}_{min} (m) \les -m $ for any $ m > 0 $.

\medskip

(v) $\tilde{E}_{min} $ is decreasing and continuous on $(0, \infty)$ and $\tilde{E}_{min}(m) \lra 0 $ as $ m \lra 0 $. 

\medskip

(vi) Let $ m > 0 $. Assume that  $(u_n)_{n \ges 1} $ is a bounded sequence in $ H^2( \R^N)$ such that $ \| u _n \|_{L^2 }^2 \lra m $ and 
$ E( u_n) \lra e $ as $ n \lra \infty$, where $ e \les - m $.
Then we have
$ \ds \liminf_{n \ra \infty } \| \Delta u _n \|_{L^2} ^2 > 0 $. 
In addition, if $ e <-  m $ then we have 
$ \ds \liminf_{n \ra \infty } \|  u _n \|_{L^{2\si + 2}} ^{2\si + 2} > 0 $.

\medskip

(vii) If $ u \in H^2( \R^N)$ satisfies $ D( u ) > 0 $ and $ P_1 ( u ) = 0 $,  we have
\beq
\label{ordered}
\frac{N\si -2 }{N \si - 4} \| u \|_{L^2}^2 > \| \nabla u\|_{L^2}^2 > \frac{N \si - 4}{N \si -2 } \| \Delta u\|_{L^2}^2
 >  \frac{ N \si}{4( \si + 1)} \left( \frac{ N \si}{2 } - 2 \right) \int_{\R^N} | u |^{ 2 \si + 2 } \, dx.
\eeq

\end{Lemma}

{\it Proof.  } (i) If $ m \les \mu_0$ (where $\mu_0$ is as in (\ref{m0})), 
we have seen that any $ u \in H^2( \R^N)$ with $ \| u \| _{L^2} ^2= m $ satisfies (\ref{cond-infl}), and then 
$u_{t_{u, 1}}\in \Oo$, $\| u_{t_{u, 1}}\|_{L^2}^2 = m $ and $ P_1( u_{t_{u, 1}}) = 0$. 
If $ m > \mu_0$, choose an integer $ n $ such that $ \frac mn < \mu_0$, and take $ v \in C_c^{\infty} ( \R^N) $ such that $ \| v \| _{L^2} ^2= \frac mn$. 
Let $ w = v_{t_{v, 1}}$, so that $ w \in C_c^{\infty} ( \R^N) $, $ \| w \|_{L^2}^2 = \frac mn$, $ P_1( w) = 0 $ and $ D(w) > 0 $. 
Choose $ R > 0  $ such that $ \mbox{supp}(w) \subset B(0, R)$, then  choose $ x_0 \in \R^N$ such that $ | x_0 | > 2R$. 
Let $ u = w + w( \cdot + x_0) + w( \cdot + 2 x_0) + \dots + w( \cdot + (n-1) x_0)$. 
Then we have $ \| u \|_{L^2}^2= n \| w \|_{L^2}^2 = m $, 
$ P_1( u ) = n P_1 ( w ) = 0 $ and $ D(u ) = n D(w) > 0 $. 
From  (\ref{deriv}) we see that $ \ph_u '(t) > 0 $ if $ t > 1$ and $t$ is close to $1$, 
hence $ \ph_u '(t_{u, infl}) > 0 $ and therefore $ u \in \Oo$. 

We use (\ref{tildE-bis}) to obtain a lower bound for $\tilde{E}_{min}(m)$. 
Let $ u \in H^2( \R^N)$ such that $ \| u \| _{L^2} ^2= m $ and $ P_1( u ) = 0 $.
From the identity $ P_1( u ) = 0 $ we obtain 
$$ \frac{1}{\si + 1} \int_{\R^N} |u |^{ 2 \si + 2} \, dx = \frac{4}{N\si} \int_{\R^N} |\Delta u |^2 \, dx -  \frac{4}{N\si} \int_{\R^N} |\nabla u |^2 \, dx .$$
Replacing this into $E(u)$ and using (\ref{interpol}) we get 
$$
\begin{array}{l}
E(u) = \left( 1 -  \frac{4}{N\si} \right) \int_{\R^N} |\Delta u |^2 \, dx 
-  \left( 2 -  \frac{4}{N\si} \right) \int_{\R^N} |\nabla u |^2 \, dx 
\\
\\
\ges 
\left( 1 -  \frac{4}{N\si} \right) \| \Delta u \|_{L^2} ^2 -  \left( 2 -  \frac{4}{N\si} \right)   \| \Delta u \|_{L^2} \| u \|_{L^2}
\ges 
{\ds \inf_{s > 0 } } \left\{ \left( 1 -  \frac{4}{N\si} \right) s^2 -  \left( 2 -  \frac{4}{N\si} \right)  m^{\frac 12} s  \right\}
\\
\\
= - \frac{ (N \si -2)^2}{N\si ( N \si - 4)}m .
\end{array}
$$
The above estimate is true for any $ u $ satisfying $ \| u \| _{L^2} ^2= m $ and $ P_1( u ) = 0 $, and (i) follows from (\ref{tildE-bis}).

\medskip

(ii)  From $ D(u) \ges d $ we get 
$ 
 \frac{1}{\si + 1} \int_{\R^N} |u |^{ 2 \si + 2} \, dx \les \frac{8}{N\si( N \si - 2)} \int_{\R^N} |\Delta u |^2 \, dx -  \frac{8d}{N\si( N \si - 2)} .
$
Using this inequality, the bound  $ E(u) \les e$, then (\ref{interpol}) and the fact that $ \| u \|_{L^2}^2 \les m $ we find
$$
\begin{array}{l}
e \ges E(u) = \left( 1 -  \frac{8}{N\si(N \si - 2)} \right) \int_{\R^N} |\Delta u |^2 \, dx 
-  2 \int_{\R^N} |\nabla u |^2 \, dx + \frac{8d}{N\si( N \si - 2)}
\\
\\
\ges 
\left( 1 - \frac{8}{N\si(N \si - 2) }\right) \| \Delta u \|_{L^2} ^2 -  2  m^{\frac 12} \| \Delta u \|_{L^2} + \frac{8d}{N\si( N \si - 2)}.
\end{array}
$$
Notice that $ 1 - \frac{8}{N\si(N \si - 2)}> 0 $ because $ N \si > 4$, and the above inequality implies that $ \| \Delta u \|_{L^2}$ is bounded.
Since $ \| u \|_{L^2}^2 \les m $,  we infer that $ \| u \|_{H^2} $ is bounded. 

\medskip

(iii) Fix $ m_1, m_2 > 0 $ and $ \e > 0 $. 
Using the density of $ C_c^{\infty}( \R^N)$ in $ H^2( \R^N)$, it is easily seen that for $ i \in \{ 1, 2 \}$  there exist 
$ u_i \in C_c^{\infty}(\R^N) \cap \Oo $ such that $ \| u_i \|_{L^2} ^2 = m_i $ and $ E( u_i )  < \tilde{E}_{min} (m_i ) + \frac{\e}{2}$. 
We may assume that $ P_1 ( u_i) = 0$ for $ i = 1, 2 $ (otherwise we replace $ u_i $ by $ (u_i)_{t_{u_i, 1}}$). 
Choose $R > 0 $ so large that $ \mbox{supp} ( u_i ) \subset B( 0, R)$ for $ i = 1, 2$. 
Choose $ x_0 \in \R^N$ such that $ | x_0 | > 2R$ and define $ u = u_1 + u_2 ( \cdot + x_0 ).$  
It is obvious that $ \| u \|_{L^2}^2 = \| u _1\|_{L^2}^2 +  \| u _2\|_{L^2}^2 = m_1 + m_2$, 
$ D( u ) = D( u_1 ) + D( u_2) > 0 $ and $ P_1( u ) = P_1( u_1 ) + P_1( u_2 ) = 0 $. 
This implies that $ P_1 ( u_t ) > 0 $ for $ t > 1 $ and $t$ close to $1$, and we infer that $\ph_u' (t_{u, infl})>0 $ 
and consequently $ u \in \Oo$. 
Then we have 
$$
 \tilde{E}_{min} (m_1 + m_2) \les E( u) = E( u_1 ) + E( u_2) \les \tilde{E}_{min} (m_1) + \tilde{E}_{min} (m_2)  + \e. 
$$
Since $\e$ is arbitrary, the conclusion follows.

\medskip

(iv) Let $ m > 0 $ and $ \e >0$. 
Let $u$ be the function constructed in the proof of Proposition \ref{Emin} (iii). 
Since $ \mbox{supp}(\wh{u}) \subset B(0, 1) \setminus B( 0, 1-\e)$, we have 
$
\| \Delta u \|_{L^2}^2 \les \|  u \|_{L^2}^2.
$
Then using (\ref{GNS}) and the fact that $ N \si > 4 $ we get 
$ \| u \|_{L^{ 2 \si + 2} }^{ 2 \si + 2 } \les B(N, \si )  \| \Delta u \|_{L^2}^{ \frac{N\si}{2} }\|u \|_{L^2}^{ 2 + 2 \si - \frac{N \si}{2}} 
\les  B(N, \si )  \| \Delta u \|_{L^2}^{ 2 }\|u \|_{L^2}^{  2 \si }, 
$
and consequently
$$
 D(u) \ges \| \Delta u \|_{L^2}^{ 2 } \left( 1 - \frac{N\si( N\si - 2) B(N, \si)}{8(\si + 1)} \| u \|_{L^2}^{ 2 \si} \right).
$$
Denote $ m_1 = \min \left( \mu_0, \left( \frac{8(\si + 1)}{N\si( N\si - 2) B(N, \si)} \right) ^{\frac{1}{2 \si}} \right)$.
If $ m < m_1$ we have $ D(u) > 0$ by the above inequality. It is obvious that $u$ satisfies (\ref{cond-infl}) because $ \|u \|_{L^2}^2 < \mu_0$, hence $ u \in \Oo.$
In the proof of Proposition \ref{Emin} (iii) we have shown that $ E(u) \les - \| u \|_{L^2}^2 + 4 \e ^2 m $, thus
$\tilde{E}_{min}(m) \les E(u) \les -m +  4 \e ^2 m $.
Since $ \e  > 0 $ is arbitrary, assertion (iv) is proven in the case $ m < m_1$. 

If $ m \ges m_1$, choose  $n\in \N^* $ such that $ \frac mn < m_1$. 
Using the sub-additivity of $\tilde{E}_{min}$ we get
 $$ \tilde{E}_{min}(m) \les n \tilde{E}_{min}\left( \frac mn \right) \les -m .$$

\medskip

(v)  Form (i) and (iv) we get $\tilde{E}_{min}(m) \lra 0 $ as $ m \lra 0 $. 
If $ 0 < m_1 < m_2$, by (iii) and (iv) we have $\tilde{E}_{min}(m_2) \les \tilde{E}_{min}(m_1) + \tilde{E}_{min}(m_2 - m_1) 
\les \tilde{E}_{min}(m_1) - (m_2 - m_1)$, thus $\tilde{E}_{min}$ is decreasing.

Fix $ M > 0 $. By (ii), the set $ \{ u \in \Oo \; | \; \| u \|_{L^2}^2 \les M \mbox{ and } E(u ) \les 0  \}$ is bounded in $H^2( \R^N)$.
Using the Sobolev embedding we see that there exists $ K = K(M, N , \si ) > 0 $ such that for any $ u $ in the above set we have 
$ \frac{1}{\si + 1} \| u \|_{L^{2 \si + 2}}^{2\si + 2} \les K . $
It is easily seen that for any $ u \in \Oo $ and any $ a \in (0, 1)$ we have $ au \in \Oo$. 
Let $ 0 < m_1 < m_2 \les M$ and denote $a = \left(\frac{m_1}{m_2} \right)^{-\frac 12}$. 
Let $ u \in \Oo $ such that $ \| u \| _{L^2}^2 = m_2$ and $ E(u ) < 0 $. 
We have $ au \in \Oo $, $ \| au \| _{L^2}^2 = m_1$ and  consequently
$$
\tilde{E}_{min} (m_1) \les E( au ) = a^2 E(u) + \frac{a^2 - a^{ 2 \si + 2}}{\si + 1} \| u \|_{L^{2 \si + 2}}^{2 \si + 2} \les a^2 E(u) + ( a^2 - a^{ 2 \si + 2} )K. 
$$
Taking the infimum in the above inequality we find
$$
\tilde{E}_{min} (m_1) \les \frac{ m_1}{m_2} \tilde{E}_{min}( m_2 ) + \left(  \frac{ m_1}{m_2} -  \frac{ m_1^{\si + 1}}{m_2^{\si + 1}} \right) K.
$$
Thus $0 < \tilde{E}_{min} (m_1) - \tilde{E}_{min} (m_2) \les 
\left( \frac{ m_1}{m_2} -1 \right) \tilde{E}_{min}( m_2 ) + \left(  \frac{ m_1}{m_2} -  \frac{ m_1^{\si + 1}}{m_2^{\si + 1}} \right) K.
$
Using (i) we infer that $  \tilde{E}_{min}$ is continuous on $(0, M)$. Since $M$ is arbitrary, (v) is proven.

\medskip

(vi)
Let $ \ell = {\ds \liminf_{n \ra \infty }} \| u _n \|_{L^{2 \si + 2}} ^{ 2 \si + 2 }$. 
If $ \ell = 0 $, 
there is a subsequence $ (u_{n_k})_{k \ges 1} $ such that $ \| u _n \|_{L^{2 \si + 2}} ^{ 2 \si + 2 }\lra 0 $ and 
using (\ref{Planch}) with $ \e = 0 $ we get $ \ds \limsup_{k \ra \infty } E( u_{n_k}) \ges -m$. 
Since $  E( u_{n_k})  \lra e \les m$ we infer that necessarily $ e = -m $. 
Moreover,  using again (\ref{Planch}) we have
$$
\int_{\R^N} \left( |\xi |^2 - 1 \right)^2 |\wh{u}_{n_k} ( \xi ) |^2 \, d \xi 
= ( 2 \pi )^{ N} \left(E(u_{n_k}) + \| u_{n_k} \|_{L^2}^2 + \frac{1}{\si + 1} \| u_{n_k} \|_{L^{2\si + 2}}^{2\si + 2}\right) \lra 0 
$$
as $ k \lra \infty.$ Using Plancherel's formula, the Cauchy-Schwarz inequality, the above convergence 
 and the boundedness of $(u_n)_{n\ges 1}$ in $H^2( \R^N)$ we get 
$$
\begin{array}{l}
\Big| \|\Delta u_{n_k} \|_{L^2}^2 -  \| u_{n_k} \|_{L^2}^2 \Big|
\les \frac{1}{(2 \pi)^N} \int_{\R^N} \Big| |\xi |^4 - 1 \Big| |\wh{u}_{n_k} ( \xi ) |^2 \, d \xi 
\\
\\
\les \frac{1}{(2 \pi)^N}  \left(  \int_{\R^N} \left( |\xi |^2 - 1 \right)^2 |\wh{u}_{n_k} ( \xi ) |^2 \, d \xi  \right)^{\frac 12} 
 \left(  \int_{\R^N} \left( |\xi |^2 + 1 \right)^2 |\wh{u}_{n_k} ( \xi ) |^2 \, d \xi  \right)^{\frac 12} \lra 0 
\end{array}
$$
 as $ k \lra \infty $
and we conclude that $ \ds \lim_{k \ra \infty} \|\Delta u_{n_k} \|_{L^2}^2 =  \lim_{k \ra \infty} \| u_{n_k} \|_{L^2}^2 =m$.

If $ \ell > 0 $,  
from (\ref{GNS}) and the fact that $ \| u _n \|_{L^2} $ is bounded it follows that there exist  $ \eta > 0 $ and $ n_0 \in \N$ such that 
$ \| \Delta u _n  \|_{L^2} \ges \eta $ for all  $ n\ges n_0$, thus $ \ds \liminf_{n \ra \infty }  \| \Delta u _n  \|_{L^2}^2 \ges \eta ^2$.

Obviously, our arguments hold for any subsequence of $(u_n)_{n \ges 1}$.  
We infer that there cannot be a subsequence 
$(u_{n_j})_{j \ges 1}$ satisfying $ \| \Delta u_{n_j}\|_{L^2} \lra 0 $ as $ j \lra \infty$, and this implies that 
$\ds \liminf_{n \ra \infty } \| \Delta u _n \|_{L^2}^2 > 0 $.

It follows from the above arguments that in the case $ e < -m$ we must have $ \ell > 0 $ and the second assertion in (vi) is now clear.

\medskip

(vii) Assume that $ u \in H^2( \R^N)$ satisfies $ D( u ) > 0 $ and $ P_1 ( u ) = 0 $.
From $ D(u ) > 0$  we get 
$ \int_{\R^N} | \Delta u |^2 \, dx > \frac{ N \si}{4( \si + 1)} \left( \frac{ N \si}{2 } - 1 \right) \int_{\R^N} | u |^{ 2 \si + 2 } \, dx $, 
which is the last inequality in  (\ref{ordered}).   Replacing this into $P_1 ( u ) = 0 $ we obtain the second inequality in (\ref{ordered}). 
Then the second inequality in (\ref{ordered}) and (\ref{interpol}) give
$ \| u \|_{L^2} >  \frac{N \si - 4}{N \si -2 }  \|\Delta  u \|_{L^2}. $
Combining this with (\ref{interpol})  we get the first inequality in (\ref{ordered}).
\hfill
$\Box$

\begin{Lemma}
\label{merdique}
(i) If $ N \ges 5 $ and $ \frac 4N < \si < 1$,  we have $  \tilde{E}_{min}(m) < -m$ for any $ m > 0 $.

\medskip

(ii) If $ \frac 4N < \si $ and $ \si \ges 1$, there exists $ m_0 > 0 $ such that $  \tilde{E}_{min}(m) = -m$ for any $ m \in ( 0, m_0] $. 
\end{Lemma}

{\it Proof. }  (i) Let $ m > 0 $. 
We use the same test functions as in the proof of Proposition \ref{P-interesting}, constructed in Example \ref{knapp}.
For small $ \e, \de > 0 $ let $ u_{\e, \de}$ be as in (\ref{u-epsilon}) and let
 $ w_{\e, \de}= \frac{\sqrt{m}}{\| u _{\e, \de}\|_{L^2}}  u_{\e, \de}$, so that $ \|  w_{\e, \de}\|_{L^2}^2 = m $.
 Fix $ \de _0 \in(0, \frac {1}{10}).$
 We have already seen in the proof of Proposition \ref{P-interesting}  that 
 $ E(w_{\e, \de_0 })  +  \| w_{\e, \de_0 }  \|_{L^{ 2}}^{2} < 0$ for all sufficiently small $ \e$ 
 (cf. (\ref{wedee})).
 
The conclusion of Lemma \ref{merdique} (i) follows if we  show that $ w_{\e, \de_0 } \in \Oo$ for all sufficiently small  $\e$. 
Since $ \mbox{supp}(\wh{w}_{\e, \de _0} ) = \mbox{supp}(\wh{u}_{\e, \de _)} ) \subset \ov{B}(0, 1 + \e) \setminus  B(0, 1 - \e )$, we have 
$ ( 1 - \e )^2 \| w_{\e, \de_0 } \|_{L^2}\les  \| \Delta  w_{\e, \de_0} \|_{ L^2} \les ( 1 + \e )^2 \| w_{\e, \de_0 } \|_{L^2}$
and 
$ ( 1 - \e ) \|  w_{\e, \de_0} \|_{L^2}\les  \| \nabla  w_{\e, \de_0} \|_{ L^2} \les ( 1 + \e ) \|  w_{\e, \de_0} \|_{L^2}$.

 By the Hausdorff-Young inequality (\ref{Hausdorff-Young}) and H\"older's inequality we have 
 $$
 \| w_{\e, \de _0}  \|_{L^{2\si + 2}}
 \les C \| \wh{w}_{\e, \de _0}  \|_{L^{\frac{2\si + 2}{2 \si + 1}}}
\les C   \| \wh{w}_{\e, \de _0}  \|_{L^2} \cdot | \mbox{supp} (\wh{w}_{\e, \de _0 } ) |^{\frac{2 \si + 1}{2 \si + 2} - \frac 12}
\les C\sqrt{m} \left(\e \de_0^{N-1} \right)^{\frac{ \si}{2 \si + 2}}.
$$
Since $ \de _0 $  is fixed, we have  
$ \| w_{\e, \de }  \|_{L^{2\si + 2}} ^{ 2 \si + 2 } \les C m^{ \si + 1} \e^{\si }$ 
and therefore 
$ D( w_{\e, \de_0 }) \ges ( 1 - \e )^2 m - C m^{\si + 1} \e^{\si } > 0 $ if $ \e $ is small enough. 
Moreover, if $H$ is given by (\ref{cond-infl}) we have $ H(w_{\e, \de_0 })  \ges C m^{- \si} \e^{- \si} \lra \infty $  as $ \e \lra 0 $ and we conclude that 
$  w_{\e, \de_0 } \in \Oo$ for all sufficiently small $ \e$.

\medskip

(ii) Recall that by (\ref{energy-rewritten}) we have 
\beq
\label{energy-rewritten-2}
E(u ) + \| u \|_{L^2}^2 = 
\| ( \Delta + 1 ) u \|_{L^2}^2 \left( 1 - \frac{\| u \|_{L^2}^{ 2 \si }}{\si + 1} Q_{\ka}(u )^{2 \si + 2} \right)
\quad \mbox{ for any } u \in H^2( \R^N) \setminus \{ 0 \}, 
\eeq
where  $ \ka = \frac{\si}{\si + 1}$ and $ Q_{\ka }$ is given in (\ref{Qka}).

Lemma \ref{basic} (vii) implies that there exists $ R_0 > 0 $ such that 
$ \|(  \Delta + 1)u \|_{L^2} \les R_0 \| u \|_{L^2}$ for any $ u \in H^2( \R^N)$ satisfying 
$ D(u)  > 0 $ and $ P_1( u) = 0 $. 
Since $ \frac 4N < \si $ and $ \si \ges 1$, condition (\ref{conditions-qualitat}) is satisfied with 
 $ s = 2$, $ p = 2 \si + 2 $ and $ \ka = \frac{\si}{\si + 1}$. Then Corollary \ref{C-quant}   implies that there exists  $M > 0 $ such that  $Q_{\ka }(u) \les M$ 
 for any $u$ as above. 
Using (\ref{tildE-bis}) and (\ref{energy-rewritten-2}) we infer that $ \tilde{E}_{min}(m) + m \ges 0$  
if $ 0 < m \les (\si + 1)^{\frac{1}{\si}}M^{- \frac{ 2 \si + 2}{\si}}$. 
The conclusion follows from this inequality and Lemma \ref{basic} (iv).
\hfill
$\Box$

\begin{Lemma}
\label{L2.4}
Let $(u_n)_{n \ges 1} \subset \Oo$ be a sequence satisfying

\medskip

(a) $ P_1( u_n ) \lra 0$, 

\medskip

(b) $ \| u _n \|_{L^2}^2 \lra m $ as $ n \lra \infty $ and $ m < \mu_0$, where $ \mu_0 $ is given by (\ref{m0}), and

\medskip

(c) there exists $ k > 0 $ such that $\| \Delta u _n \|_{L^2} \ges k $ for all $ n $. 

\medskip

Then $ \ds \liminf_{n \ra \infty } D( u_n ) > 0 $. Moreover, if $\| \Delta u _n \|_{L^2} $ is bounded then we have $ \ds \liminf_{n \ra \infty } t_{u_n, infl} > 1$. 
\end{Lemma}

{\it Proof. } 
We have $ D(u_n) > 0 $ for all $n $ because $ u_n \in \Oo$. 
We argue by contradiction and we assume that there is a subsequence, still denoted $(u_n)_{n \ges 1}$, such that 
$ D( u_n ) \lra 0 $. 
We have
\beq
\label{dn}
\frac{ N \si }{ 4( \si + 1 )} \int_{\R^N} | u _n |^{ 2 \si + 2} \, dx = \frac{2}{N \si - 2} \left( \int_{\R^N} |\Delta u _n |^2 \, dx - D ( u_n)\right). 
\eeq
Using this identity in the expression of $ P_1( u_n)$  
we get 
\beq
\label{P1n}
 P_1( u_n) = \frac{ N \si - 4}{N \si - 2} \int_{\R^N} |\Delta u _n |^2 \, dx - \int_{\R^N} |\nabla u _n |^2 \, dx + \frac{2}{N \si - 2} D( u_n). 
\eeq
From the equality above  and (\ref{interpol}) we obtain 
$$
\| \Delta u_n \|_{L^2}\|  u_n \|_{L^2}    >  \| \nabla u _n \|_{L^2}^2 = \frac{ N \si - 4}{N \si - 2} \| \Delta u_n \|_{L^2}^2 + \frac{2}{N \si - 2} D( u_n) - P_1 ( u_n). 
$$
The last inequality, assumptions (a) and (b) and the fact that $ D( u_n ) \lra 0 $ imply that $ \| \Delta u_n \|_{L^2}$ is bounded. 
We rewrite the last inequality  in the form 
\beq
\label{2.9}
\| \Delta u_n \|_{L^2} <  \frac{ N \si - 2}{N \si - 4} \left( \|u _n \|_{L^2} -  \frac{2}{N \si - 2} \frac{D( u_n)}{\| \Delta u_n \|_{L^2} } + \frac{P_1 ( u_n)}{\| \Delta u_n \|_{L^2} }  \right).
\eeq
Using the definition of $D$ (see (\ref{defD})) and (\ref{GNS}) we get 
$$
\| \Delta u_n \|_{L^2}^2 - D( u_n) =\frac{ N \si ( N \si -2)}{8 ( \si + 1)}  \| u _n \|_{L^{2 \si + 2 }}^{ 2 \si + 2 } 
\les 
\frac{ N \si ( N \si -2)B(N, \si )  }{8 ( \si + 1)}  \| \Delta u_n \|_{L^2}^{ \frac{N \si}{2}} \| u _n \|_{L^2}^{ 2 \si + 2 - \frac{N \si}{2}}.
$$
Dividing by $ \| \Delta u_n \|_{L^2}^2 $ and using (\ref{2.9}) we discover
$$
\begin{array}{l}
1 - \frac{ D(u_n)}{ \| \Delta u_n \|_{L^2}^2} \les 
\frac{ N \si ( N \si -2)B(N, \si )  }{8 ( \si + 1)}  \| \Delta u_n \|_{L^2}^{ \frac{N \si}{2}-2} \| u _n \|_{L^2}^{ 2 \si + 2 - \frac{N \si}{2}}
\\
\\
< \frac{ N \si ( N \si -2)B(N, \si )  }{8 ( \si + 1)} \left[    \frac{ N \si - 2}{N \si - 4} \left( \|u _n \|_{L^2} -  \frac{2}{N \si - 2} \frac{D( u_n)}{\| \Delta u_n \|_{L^2} } + \frac{P_1 ( u_n)}{\| \Delta u_n \|_{L^2} } \right)    \right] ^{ \frac{N \si}{2}-2} \| u _n \|_{L^2}^{ 2 \si + 2 - \frac{N \si}{2}}.
\end{array}
$$
Letting $ n \lra \infty $ in the above inequality and using assumptions (b), (c) and (\ref{m0}) we obtain $ 1 \les \frac{ m^{\si }}{\mu_0^{\si }}$, 
contradicting the fact that $ m < \mu_0$. 
We have thus proved that $ \ds \liminf_{n \ra \infty } D( u_n ) > 0 $. 

Since $ u_n \in \Oo$ we have $ t_{u_n ,infl } > 1 $ for each $n$. 
We argue again by contradiction for the second part of Lemma \ref{L2.4} and we assume that there is a subsequence, still denoted 
$(u_n)_{n \ges 1}$, such that $ t_{u_n, infl } \lra 1$ as $ n \lra \infty$. 
By (\ref{tinfl})  we have   $ D( u_n ) = \left( 1 - t_{u_n, infl}^{- \frac{N \si}{2} + 2} \right) \| \Delta u_n \|_{L^2}^2 $ and the boundedness of $ \| \Delta u _n \|_{L^2}$ implies that
$ D( u_n ) \lra 0 $,  contradicting assumption (c).
\hfill
$\Box$

\begin{Lemma}
\label{L2.5}
Assume that $m< \mu_0$, where $\mu_0$ is given by (\ref{m0}).
Suppose that the sequence  
$ (u_n)_{n \ges 1} \subset H^2( \R^N) $ 
 satisfies $ \| u_n \|_{L^2}^2 \lra m $ and $ D( u_n ) \lra 0 $ as $ n \lra \infty$. 

Then we have $ \ds \liminf_{n \ra \infty} E( u_n ) \ges   \tilde{E}_{min}(m). $

Moreover, if  $  \tilde{E}_{min}(m) < -m$ we have $ \ds \liminf_{n \ra \infty} E( u_n ) >   \tilde{E}_{min}(m). $
\end{Lemma}

{Notice that in Lemma \ref{L2.5} we do not assume that $ (u_n)_{n \ges 1} \subset \Oo $.}

{\it Proof. } 
The sequence $(u_n)_{n \ges 1}$ is bounded in $H^2( \R^N)$ by Lemma \ref{basic} (ii).
By  Lemma \ref{basic} (iv) we have  $ \tilde{E}_{min}(m ) \les -m $, so 
the conclusion of Lemma \ref{L2.5} is obvious if $  \ds \liminf_{n \ra \infty} E( u_n )\ges -m$.
Form now on we only consider the case when $  \ds \liminf_{n \ra \infty} E( u_n ) < -m$.
Passing to a subsequence 
we may assume that $E( u_n)\lra  e < -m$ as $ n \lra \infty $ and that $ \| u_n\|_{L^2}^2 < \mu_0 $ for all $n\ges 1$, so that $ \ph_{u_n}'(t_{u_n, infl}) > 0 $ and $ t_{u_n, 1}$,   $ t_{u_n, 2}$ do exist. 

It follows from Lemma \ref{basic} (vi) that there exist  $ \eta > 0 $ and $ n_0 \in \N$ such that 
$ \| \Delta u _n  \|_{L^2} \ges \eta $ for all  $ n\ges n_0$. 
Using (\ref{dn}) and (\ref{GNS}) we get for $ n \ges n_0$,
$$
\frac{N\si }{4(\si + 1)}  B(N, \si ) \| \Delta u_n \|_{L^2}^{ \frac{N \si}{2} - 2} \| u _n \|_{L^2}^{ 2 \si + 2 -  \frac{N \si}{2} } 
\ges 
\frac{2}{N \si -2} \left( 1 - \frac{1}{\eta ^2} D( u_n ) \right)
$$
Then using (\ref{P1n}), (\ref{interpol})  and the above inequality we obtain for $ n \ges n_0$
$$
\begin{array}{l}
P_1( u_n) \ges \| \Delta u_n \|_{L^2} \left( \frac{ N \si - 4}{N \si - 2} \| \Delta u _n \|_{L^2} - \| u _n \|_{L^2} \right)+ \frac{2}{N \si - 2} D( u_n)
\\
\\
\ges \! \| \Delta u_n \|_{L^2} \left( \! \!  \frac{ N \si - 4}{N \si - 2} \left( \frac{8(\si + 1) }{N \si ( N \si - 2 ) B(N, \si ) }
 \left( \! 1 - \frac{1}{\eta ^2} D( u_n ) \right) \right)^{\frac{2}{N \si - 4}} \| u_n \|_{L^2}^{ - \frac{4 \si + 4 - N \si }{N \si - 4}} 
 \! \! - \| u _n \|_{L^2} \! \! \right)+ \frac{2 D( u_n)}{N \si - 2}.
\end{array}
$$
Letting $ n \lra \infty $ and using the fact that $ D(u_n) \lra 0 $  we discover
\beq
\label{borneinf}
\liminf_{n \ra \infty} P_1( u_n) \ges \eta 
\left(  \frac{ N \si - 4}{N \si - 2}  \left( \frac{8(\si + 1) }{N \si ( N \si - 2 ) B(N, \si ) } \right)^{\frac{2}{N \si - 4}}m ^{ - \frac{2 \si + 2 - \frac{N \si}{2} }{N \si - 4}} 
 - m^{\frac 12}   \right).
\eeq
Since $ 0 < m < \mu_0$, where $ \mu_0 $ is given by (\ref{m0}), the right-hand side of (\ref{borneinf}) is equal to 
$ \eta m^{\frac 12} \left( \left(\frac{\mu_0}{m} \right)^{\frac{2 \si}{N \si - 4}} - 1 \right) $ and this quantity is positive. 
We conclude that there exists $ n_1 \in \N$ such that $ P_1 ( u_n ) > 0 $ for all $ n \ges n_1$. 
This means that $ t_{u_n ,1} < 1 <  t_{u_n ,2}$ for all $ n \ges n_1$.

We denote $ v_n = (u_n)_{t_{u_n, 1}}$, so that $ v_n \in \Oo$, $ \| v_n \|_{L^2} = \| u _n \|_{L^2}$, $ P_1( v_n ) = 0 $  and $ E ( v_n ) \les E( u_n)$ for each $n\ges n_1$ (recall that $ t \longmapsto E( u_t)$ is increasing on $[t_{u, 1}, t_{u, 2}]$).  
By Lemma \ref{basic} (ii), the sequence   $(v_n)_{n \ges 1}$ is bounded in $ H^2( \R^N)$. 
Since $ \tilde{E}_{min} ( \| v_n \|_{L^2}^2 ) \les E( v_n) \les E(u_n)$, passing to the limit and using the continuity of $ \tilde{E}_{min}$ (see Lemma \ref{basic} (v)) we get 
\beq
\label{2.10}
 \tilde{E}_{min} (m) \les \liminf_{n \ra \infty} E( v_n ) \les \limsup_{n \ra \infty} E( v_n ) \les \lim_{n \ra \infty} E( u_n) = e . 
\eeq

We show that if $  \tilde{E}_{min} (m) < -m$, then at least one inequality in  (\ref{2.10}) is strict. This  clearly implies the conclusion of  
 Lemma \ref{L2.5}. 
We assume that equality occurs in the first two inequalities in (\ref{2.10}), 
which means precisely that $ E( v_n ) \lra  \tilde{E}_{min} (m) < - m $. 
We show that in this case the last inequality in (\ref{2.10}) must be strict. 
Denote $ \ds \ell  := \liminf_{n \ra \infty } \| v_n \|_{L^{ 2 \si + 2} }^{ 2 \si + 2}  $. 
Using Lemma \ref{basic} (vi) we see that $ \ell > 0 $ and  there exist $ \eta _1 > 0 $ such that 
$ \| \Delta v_n \|_{L^2} \ges \eta _1 $ for all sufficiently large $ n $. 
Now we may apply Lemma \ref{L2.4} to $( v_n)_{n \ges 1}$ and we infer that $ \ds \liminf_{n \ra \infty} D( v_n) > 0 $ and $ \ds \liminf_{n \ra \infty } t_{v_n, infl } > 1$. 

Denote $ s_n = ( t_{u_n, 1 })^{-1}$, so that $ u_n = ( v_n )_{s_n}$. 
Recall that $ t_{u_n, 1 } < 1$, hence $ s_n  > 1 $ for all $ n \ges n_1$.
We have 
$$
D( u_n ) = D( (v_n)_{s_n}) = s_n ^2 \left( \int_{\R^N} |\Delta v_n |^2 \, dx - s_n ^{\frac{N\si}{2} - 2} \frac{N\si ( N \si - 2)}{8( \si + 1)} 
\int_{\R^N} |v_n |^{ 2 \si + 2} \, dx \right) .
$$
Since $D(u_n ) \lra 0$,  the second factor in the expression of $D( (v_{n})_{s_{n}}) $ here above must tend to $ 0 $, and from (\ref{tinfl}) and the fact that 
$ \ell > 0 $ we infer that $ r_n := \frac{ s_{n}}{t_{v_{n},infl}} \lra 1 $ as $ k \lra \infty$. 
Using the boundedness of $ ( u_n )_{n \ges 1}$ in $ H^2( \R^N)$ 
we obtain   then
\beq
\label{2.11}
E( u_{n} ) - E( (v_{n})_{t_{v_{n}, infl}}) 
= E( u_{n} )- E( (u_{n} )_{r_n ^{-1}})
\lra 0 \qquad \mbox{ as } n \lra \infty. 
\eeq
From Lemma \ref{elem} (v) we have 
\beq
\label{2.12}
 E( (v_{n})_{t_{v_{n}, infl}})- E( v_{n}) 
= \frac{h(t_{v_{n}, infl})}{\si + 1}  \int_{\R^N} |v_{n}|^{ 2 \si + 2} \, dx . 
\eeq
Fix $ t_* $ such that $ 1 < t_* < \ds \liminf_{n \ra \infty } t_{v_n, infl } $. From (\ref{2.11}) and (\ref{2.12}) we see that 
$E( u_{n} ) - E( v_{n}) > \frac{\ell}{2 \si + 2} h(t_*)$ for all sufficiently large $n$. 
Therefore the last inequality in (\ref{2.10}) is strict and  the conclusion of Lemma \ref{L2.5} follows.
\hfill
$\Box$

\begin{Lemma}
\label{decrease}
Let $\mu_0$ be as in (\ref{m0}). Denote $ m_0 = \ds \inf \{ m \in (0, \mu_0] \; | \; \tilde{E}_{min}(m) < -m \}$. Then:

(i) The mapping $ m \longmapsto \frac{ \tilde{E}_{min} (m)}{m} $ is non-increasing on $(0, \mu_0]$, and it is decreasing on $(m_0, \mu_0]$. 

(ii) If $ m \in (0, \mu_0]$ satisfies $\tilde{E}_{min}(m) < -m $, then for any $ m' \in (0, m)$ we have 
$$
\tilde{E}_{min}(m) <  \tilde{E}_{min}({m'}) + \tilde{E}_{min}(m - {m'}) .
$$

\end{Lemma}

{\it Proof. } 
It is easy to see that for any  $ u \in \Oo$  and  any $ a \in (0, 1)$ we have $ au \in \Oo $. 

Assume  that $ u \in H^2( \R^N)$ satisfies $ D( u ) > 0 $,  $ P_1 ( u ) = 0 $ and $ \| u \|_{L^2}^2 < \mu_0$. 
We show that for any $ a \in \left[1, \ds \frac{\mu_0}{\| u \|_{L^2}^2} \right]$ we have $ a^{\frac 12 } u \in \Oo.$
Since $ \| a^{\frac 12 } u \|_{L^2} ^2 = a \| u \|_{L^2}^2\les \mu_0$, the function $ a^{\frac 12 } u$ automatically satisfies 
(\ref{cond-infl}) and we only  have to prove that $ D( a^{\frac 12} u ) >0 $. 
Using (\ref{GNS}) and the fact that $ \| u \|_{L^2} >  \frac{N \si - 4}{N \si -2 }  \|\Delta  u \|_{L^2} $ (see (\ref{ordered})) we have
$$
\begin{array}{l}
D(  a^{\frac 12 } u ) 
= a \| \Delta u \|_{L^2}^2 -  \frac{ N \si}{4( \si + 1)} \left( \frac{ N \si}{2 } - 1 \right) a^{ \si + 1} \| u \|_{L^{2 \si + 2}}^{ 2 \si + 2}
\\
\\
\ges a \| \Delta u \|_{L^2}^2 \left( 1 - \frac{ N \si}{4( \si + 1)} \left( \frac{ N \si}{2 } - 1 \right) a^{ \si} B(N, \si ) \| \Delta u \|_{L^2}^{ \frac{ N \si }{2}-2} \| u \|_{L^2} ^{ 2 \si + 2 - \frac{N \si}{2}} \right)
\\
\\
> a \| \Delta u \|_{L^2}^2 \left( 1 - \frac{ N \si}{4( \si + 1)} \left( \frac{ N \si}{2 } - 1 \right) B(N, \si )  \left( \frac{ N \si - 2}{N \si -4} \right)^{ \frac{ N \si }{2}-2}  a^{ \si}  \| u \|_{L^2} ^{ 2 \si } \right) .
\end{array}
$$
The last expression is non-negative if $ a  \| u \|_{L^2} ^2 \les \mu_0$ by  (\ref{m0}), hence  $  a^{\frac 12} u \in \Oo$.
Thus we have 
\beq
\label{compare}
\tilde{E}_{min}( a \| u \|_{L^2}^2) 
\les E( a^{\frac 12} u ) = a E(u) + \frac{a - a^{ \si + 1}}{\si + 1} \int_{\R^N} | u |^{ 2 \si + 2}\, dx   \mbox{ for any } 
 a \in \left( \left. 0,  \ds \frac{\mu_0}{\| u \|_{L^2}^2}  \right] \right. .
\eeq

Let  $ m \in (0, \mu_0)$. 
Take a minimising sequence $ ( u_n)_{n \ges 1 } \subset \Oo $ such that $ \| u _n \|_{L^2} = m $, 
$ P_1 ( u_n ) = 0 $ and $ E( u_n ) \lra \tilde{E}_{min} (m) $. 
By (\ref{Planch}) we have $ \frac{1}{\si + 1 } \| u _n \|_{L^{2 \si + 2}}^{ 2 \si + 2 } \ges - ( E( u_n) + \| u _n \|_{L^2}^2)$ for each $n$.
Using this in (\ref{compare}) and letting $ n \lra \infty $ we obtain 
$$
\tilde{E}_{min}( am ) \les a \tilde{E}_{min} (m) + ( a ^{ \si + 1 } - a ) ( \tilde{E}_{min}(m) + m)
 \mbox{ for any } 
 a \in \left[ 1,  \ds \frac{\mu_0}{m}  \right]
$$
or equivalently
\beq
\label{decreasing}
\frac{\tilde{E}_{min}( am )}{am } \les \frac{ \tilde{E}_{min} (m) }{m}+ ( a^{ \si } - 1) \left(  \frac{ \tilde{E}_{min} (m) }{m} + 1 \right)
 \mbox{ for any } 
 a \in \left[ 1,   \frac{\mu_0}{m}  \right].
\eeq
Since $ \tilde{E}_{min} ( m ) \les -m $ (see Lemma \ref{basic} (iv)),  conclusion (i)  of Lemma \ref{decrease} follows easily from (\ref{decreasing}). 

(ii) It follows from the continuity of $ \tilde{E}_{min} $ that $ \tilde{E}_{min}  ( m' ) < -m'$ for $m'$ in a neighbourhood of $ m$, 
and using part (i) we infer  that $ \frac{\tilde{E}_{min} (m) }{m } <  \frac{\tilde{E}_{min} (m') }{m '} $ for all $ m' \in (0, m)$. 
In particular, for ${m'} \in (0, m )$ we have 
 $ \frac{\tilde{E}_{min} (m) }{m } <  \frac{\tilde{E}_{min} ({m'}) }{{m'}} $ 
and 
 $ \frac{\tilde{E}_{min} (m) }{m } <  \frac{\tilde{E}_{min} (m - {m'}) }{m - {m'}} $.
Combining the last two inequalities we get (ii).
\hfill
$\Box $

\begin{Theorem}
\label{T2.6}
Assume that $ 0 < m < \mu_0 $ and $ \tilde{E}_{min}(m) < -m$. 
Then $ \tilde{E}_{min}(m) $ is achieved. 

Moreover, for any sequence $ (u_n)_{n \ges 1} \subset \Oo $ satisfying $ \| u _n \|_{L^2} ^2 \lra m $ and $ E( u_n ) \lra  \tilde{E}_{min}(m)$ there exist a subsequence $ (u_{n_k})_{k \ges 1} $, a sequence of points $ ( x_k ) \subset \R^N$ and $ u \in \Oo $ such that 
$u_{n_k} ( \cdot + x _k ) \lra u $ strongly in $ H^2( \R^N)$. 
(Then, obviously, $ \| u \|_{L^2}^2 = m $ and $ E(u ) =  \tilde{E}_{min}(m)$.)

\end{Theorem}

{\it Proof. } 
Let  $(u_n)_{n \ges 1}$ be a sequence as in Theorem \ref{T2.6}.
It follows from Lemma \ref{basic} (ii) that $(u_n)_{n \ges 1}$ is bounded in $ H^2( \R^N)$. 
Lemma \ref{basic} (vi) implies that
 there exist $ \de > 0 $ and $ \ell > 0 $ such that $ \| \Delta u_n \|_{L^2} \ges \de $ 
and $ \| u_n \|_{L^{2 \si + 2}}^{ 2 \si + 2} \ges \ell $ for all sufficiently large $ n $. 

Proceeding exactly as in the proof of Theorem \ref{Global} we see that there exists a subsequence, still denoted $ (u_n)_{n \ges 1} $, 
 there exist points $ x_n \in \R^N$ and there is $ u \in H^2(\R^N)$, $ u \neq 0 $ such that after replacing $ u_n$ by $ u_{n}( \cdot + x_n)$, 
(\ref{conv}) holds. 
Then the weak convergence $ u_n \rightharpoonup u $ gives (\ref{conv1}), while Brezis-Lieb Lemma 
and the fact that $ u_n \lra u $ a.e. give  (\ref{conv2}).

We denote $ v_n = (u_{n })_{t_{u_n, 1}}$. 
Then we have $ v_n \in \Oo$ for all $ n$, $ \| v _n \|_{L^2} = \| u _n \|_{L^2}$, $ P_1 ( v_n ) = 0 $, 
and $ \tilde{E}_{min} ( \| v_n \|_{L^2}^2) \les E ( v_n ) \les E( u_n)$ for all $n$, thus  
\beq
\label{v-min}
 E( v_n ) \lra  \tilde{E}_{min} (m ) < -m \qquad \mbox{  as } n \lra \infty.
\eeq
Lemma \ref{basic} (ii) implies that $ ( v_n )_{n \ges 1}$ is bounded in $ H^2( \R^N)$, and by Lemma \ref{basic} (vi) 
there exist $ \tilde{\de } > 0 $ and $ \tilde{\ell} > 0 $ such that 
$ \| \Delta v_n \|_{L^2} \ges \tilde{\de } $
and $ \| u_n \|_{L^{2 \si + 2}}^{ 2 \si + 2} \ges \tilde{\ell} $ for all sufficiently large $ n $. 
We have $ \| \Delta v _n \| _{L^2}= t_{u_n, 1}\| \Delta u _n \| _{L^2}$,  and 
$\| \Delta v _n \| _{L^2} $ as well as $ \| \Delta u _n \| _{L^2}$ are bounded and stay away from zero, 
thus the sequence $(  t_{u_n, 1} )_{n \ges 1}$ is bounded and stays away from zero. 
We infer that there exist $ t_1 \in (0, \infty)$ such that after passing to a subsequence of $(u_n)_{n \ges 1}$, still denoted the same, 
we have  $  t_{u_n, 1}  \lra t_1 $ as $n \lra \infty$. 
It is easy to see that $ (u_n)_{t_{u_n, 1}} \rightharpoonup u_{t_1}\neq 0$ as $ n \lra \infty$. 
Let $ v =  u_{t_1}$. Then $ v \neq 0 $ and  $ v_n \rightharpoonup v $ weakly in $ H^2( \R^N)$.
Passing eventually to further subsequences of $ ( u_n )_{n \ges 1}$ and of $ ( v_n )_{n \ges 1}$, still denoted the same, 
we may assume in addition that $ v_n \lra v $ in $ L_{loc}^p ( \R^N )$ for any $ 1 \les p < 2^{**}$, and almost everywhere. 
It is then clear that (\ref{conv1}) and (\ref{conv2}) hold with $ v_n$ and $v$ instead of $ u_n$ and $u$, respectively.

Our strategy is as follows. 
We will show firstly that $ \| v \|_{L^2}^2 = m $. Then we prove that $ v_n \lra v $ strongly in $ H^2( \R^N)$ and that 
$ v \in \Oo$. Finally we show that necessarily $ t_ 1 = 1$ (thus $ u = v $) and that $ u_n \lra u $ strongly in $ H^2( \R^N)$.

To carry out the first step of our plan we argue by contradiction and we assume that $ \| v \|_{L^2}^2 < m $.
Then (\ref{conv1}) implies that $ \| v_n - v \|_{L^2}^2 \lra m - \| v \|_{L^2}^2 \in (0, m)$ as $ n \lra \infty$. 
By (\ref{conv1}) and (\ref{conv2}) we have 
\beq
\label{cut-e}
E( v_n) = E(v) + E( v_n - v ) + o(1),
\eeq
\beq
\label{cut-d}
D( v_n) = D(v) + D( v_n - v ) + o(1), \quad \mbox{ and }
\eeq
\beq
\label{cut-p}
0 = P_1( v_n) = P_1(v) + P_1( v_n - v) + o(1) \qquad \mbox{ as } n \lra \infty . 
\eeq
Passing to a further subsequence if necessary, we may assume that 
$$
\int_{\R^N} |\Delta v_n - \Delta v |^2 \, dx \lra a, \quad
\|\nabla v_n - \nabla v \|_{L^2}^2  \lra b \quad \mbox{ and } \quad
\frac{1}{\si + 1} \int_{\R^N} | v_n - v |^{ 2 \si + 2 } \, dx \lra c 
$$
as $ n \lra \infty$, where $a, \, b, \, c \ges 0 $. 
Notice that $ \ds \liminf_{n \ra \infty } D( v_n ) > 0 $ by Lemma \ref{L2.4}, and using (\ref{cut-d}) and passing to the limit we infer that
$ D(v ) +  a - \frac{ N \si}{4} \left( \frac{N \si}{2} - 1 \right) c  > 0 $. 
Thus at least one of the quantities $D(v) $ and $ a - \frac{ N \si}{4} \left( \frac{N \si}{2} - 1 \right) c  $ must be positive.
There are several possibilities and we analyse all of them, showing that in each case we get a contradiction.

\medskip

{\it Case 1. } $ D(v ) > 0 $ and $ a - \frac{ N \si}{4} \left( \frac{N \si}{2} - 1 \right) c  > 0 $. 
In this case we have $ v \in \Oo $ and $ D( v_n - v ) > 0 $ (thus $ v_n - v \in \Oo $) for all sufficiently large $ n $, hence
$ E(v ) \ges \tilde{E}_{min} ( \| v \|_{L^2}^2) $ and $ E(v- v_n) \ges \tilde{E}_{min} ( \| v - v_n\|_{L^2}^2) $.
Using (\ref{cut-e}), (\ref{v-min}) and the continuity of $ \tilde{E}_{min}$ we get 
\beq
\label{superadd}
 \tilde{E}_{min} (m) \ges  \tilde{E}_{min}( \| v \|_{L^2}^2 ) +  \tilde{E}_{min} ( m - \| v \|_{L^2}^2)
\eeq
and this  contradicts Lemma \ref{decrease} (ii).

\medskip

{\it Case 2. } $ D(v ) > 0 $ and $ a - \frac{ N \si}{4} \left( \frac{N \si}{2} - 1 \right) c  = 0 $. 
We have $ v \in \Oo $ and $ D( v_n - v ) \lra 0 $, and Lemma \ref{L2.5} implies
 $\ds \liminf_{n \ra \infty } E( v_n - v ) \ges \tilde{E}_{min}(m - \| v \|_{L^2}^2)$. 
Proceeding as in the first case we get (\ref{superadd}), and this is in contradiction with Lemma \ref{decrease} (ii). 

\medskip

{\it Case 3. } $ D(v ) = 0 $ and $ a - \frac{ N \si}{4} \left( \frac{N \si}{2} - 1 \right) c  > 0 $. 
As in Case 1, for all sufficiently large $n $ we have $ D( v_n - v ) > 0 $, hence $ v_n - v \in \Oo $ and we find 
 $\ds \liminf_{n \ra \infty } E( v_n - v ) \ges \tilde{E}_{min}(m - \| v \|_{L^2}^2)$. 
We have $ t_{v, infl} = 1$, hence $ E(v) > E( v_{t_{v,1}}) \ges \tilde{E}_{min}( \| v \|_{L^2}^2)$ and 
using (\ref{cut-e}) we get (\ref{superadd}) (with strict inequality),  contradicting  again Lemma \ref{decrease} (ii). 

\medskip

{\it Case 4. }  $ a - \frac{ N \si}{4} \left( \frac{N \si}{2} - 1 \right) c  < 0 $. 
In this case we have necessarily $ D(v ) > 0 $, hence $ v \in \Oo $ and $ E(v) \ges \tilde{E}_{min}( \| v \|_{L^2}^2). $
We distinguish two subcases: 

{\it Subcase A. } There is a subsequence $(v_{n_k})_{k \ges 1}$ such that $ P_1( v_{n_k} - v ) \ges 0 $, 
that is $ t_{(v_{n_k} - v) , 1 } \les 1 \les  t_{(v_{n_k} - v) , 2 }$.
For all $ k $ sufficiently large we have $ \| v_{n_k} - v  \|_{L^2}^2 < \mu_0 $ and then $(v_{n_k} - v )_{t_{(v_{n_k} - v) , 1 } } \in \Oo$, hence
$$
E(  v_{n_k} - v ) \ges E \left( (v_{n_k} - v )_{t_{(v_{n_k} - v) , 1 } } \right) \ges \tilde{E}_{min} \left( \| v_{n_k} - v  \|_{L^2}^2  \right). 
$$
Letting $ k \lra \infty $ we discover $ \ds \liminf_{k \ra \infty } E(  v_{n_k} - v )  \ges \tilde{E}_{min} \left( m - \| v \|_{L^2}^2 \right).$
Then using (\ref{v-min}) and (\ref{cut-e})  for the subsequence $( v_{n _k})_{k \ges 1}$ 
we infer that (\ref{superadd}) holds and we reach a contradiction as in the previous cases. 

{\it Subcase B. }  $ P_1( v_n  - v ) < 0 $ for all sufficiently large $n$. 
To simplify notation, let $ w_n = v_n - v$. 
Since $ v_n $ satisfies (\ref{conv}), we have $ w_n \rightharpoonup 0 $ weakly in $ H^2( \R^N)$ and 
$ w_n \lra 0 $ strongly in $ L^p ( \R^N)$ for all $ p \in [1, 2^{**})$ and almost everywhere, and it is clear that 
for any fixed $ t > 0 $, the sequence $( (w_n)_t)_{n \ges 1}$ has the same properties. 
We fix  $ t > 1 $ (and $t$ sufficiently close to $1$) such that  
$$
D(v) + a t^2 - \frac{ N \si}{4} \left( \frac{N \si}{2} - 1 \right) c  t^{\frac{N \si}{2}} > 0 .
$$
Such $t$ exists because $ D(v ) +  a - \frac{ N \si}{4} \left( \frac{N \si}{2} - 1 \right) c  > 0 $. 
Let $ \tilde{v}_n = v + ( w_n )_t. $ 
The weak convergence $ ( w_n )_t  \rightharpoonup 0 $ weakly in $ H^2( \R^N)$ gives
$$
\| \tilde{v}_n \|_{L^2}^2 = \| v \|_{L^2}^2 +  \| (w_n)_t \|_{L^2}^2  + o (1) = \| v \|_{L^2}^2 +  \| w_n\|_{L^2}^2  + o (1) 
= \| v _ n \|_{L^2}^2 + o (1) = m + o(1), 
$$
$$
\| \nabla \tilde{v}_n \|_{L^2}^2 = \| \nabla v \|_{L^2}^2 +  \| \nabla (w_n)_t \|_{L^2}^2  + o (1) , 
$$
$$
\| \Delta \tilde{v}_n \|_{L^2}^2 = \| \Delta v \|_{L^2}^2 +  \| \Delta (w_n)_t \|_{L^2}^2  + o (1) .
$$
Since $ ( w_n )_t \lra 0 $ a.e. and $ \| ( w_n )_t  \|_{L^{ 2 \si + 2}}^{ 2 \si + 2 } $ is bounded, using Brezis-Lieb  Lemma we have
$$
\| \tilde{v}_n \|_{L^{2 \si + 2 }}^{ 2 \si + 2 } = \|v \|_{L^{2 \si + 2 }}^{ 2 \si + 2 } + \| (w_n)_t\|_{L^{2 \si + 2 }}^{ 2 \si + 2 } + o (1) 
=  \|v \|_{L^{2 \si + 2 }}^{ 2 \si + 2 } + t^{\frac{ N \si}{2} }  \| w_n\|_{L^{2 \si + 2 }}^{ 2 \si + 2 } + o (1).
$$
In particular, we infer that 
\beq
\label{cut-ee}
E(\tilde{ v} _n) = E(v) + E((w_n)_t ) + o(1) \qquad \mbox{ as } n \lra \infty.
\eeq

It follows from the above that 
$ D ( \tilde{v}_n ) \lra D(v) + a t^2 - \frac{ N \si}{4} \left( \frac{N \si}{2} - 1 \right) c  t^{\frac{N \si}{2}} > 0 $ as $ n \lra \infty$, 
hence $ D( \tilde{v}_n) > 0 $  and consequently $ \tilde{v}_n \in \Oo $ for all sufficiently large $n$. 
This implies $ E( \tilde{v}_n ) \ges \tilde{E}_{min}( \| \tilde{v} \| _{L^2}^2 ) $, and letting $ n \lra \infty $ and using the continuity of $\tilde{E}_{min} $ we get 
\beq
\label{contr-4B-1}
\liminf_{n \ra \infty }  E( \tilde{v}_n ) \ges \tilde{E}_{min}(m ) . 
\eeq
On the other hand, from (\ref{cut-e}) and (\ref{cut-ee}) we get 
$$
E( \tilde{v}_n ) - E( v_n ) = E( (w_n )_t ) - E( w_n ) + o (1).
$$
For all sufficiently large $n$ (so that $ P_1 ( w_n ) \les 0 $),   using (\ref{deriv}) and Lemma \ref{elem} (iv) we obtain
$$
E( (w_n )_t ) - E( w_n ) \les 2( t-1) P_1( w_n ) -  ( t-1 )^2 D( w_n ) \les  ( t-1 )^2 D( w_n ).
$$
We infer that 
\beq
\label{contr-4B-2}
\begin{array}{l}
\ds \limsup_{n \ra \infty } E( \tilde{v}_n ) \les \lim_{ n \ra \infty } E( v_n ) + ( t-1 )^2  \lim_{ n \ra \infty } D( w_n )
\\
\\
= \tilde{E}_{min}(m )+  ( t-1 )^2 \left(  a - \frac{ N \si}{4} \left( \frac{N \si}{2} - 1 \right) c \right)
<  \tilde{E}_{min}(m )
\end{array}
\eeq
and this is in contradiction with (\ref{contr-4B-1}).

\medskip

{\it Case 5. } $ D(v ) < 0 $. 
This case is very similar to Case 4, and a bit simpler. 
We have necessarily $ a - \frac{ N \si}{4} \left( \frac{N \si}{2} - 1 \right) c  > 0 $, thus 
$ D( v_n - v ) > 0 $ and consequently $ v _n - v \in \Oo$ for all sufficiently large $n$, 
 and we find $ E ( v_n - v ) \ges \tilde{E}_{min} ( \| v _n - v \|_{L^2}^2 ) $, which implies 
$ \ds \liminf_{n \ra \infty } E ( v_n - v) \ges \tilde{E}_{min} ( m - \| v \|_{L^2}^2)$. 

If $ P_1 ( v) \ges 0 $ we have $  t_{1, v} \les 1 \les t_{2, v } $, hence   $ E( v ) \ges E( v_{t_{v, 1}}) \ges \tilde{E}_{min} ( \| v \|_{L^2}^2)$ and we find that (\ref{superadd})  holds. 

If $ P_1 ( v) < 0 $ we may choose $ t > 1 $ such that $ D( v_t ) +  a - \frac{ N \si}{4} \left( \frac{N \si}{2} - 1 \right) c  > 0 $.
Denoting $ v_n ^{\sharp} = v_t + v_n - v $ we see that $ \| v_n ^{ \sharp } \|_{L^2}^2 \lra m $ and 
$ D (  v_n ^{ \sharp }  ) > 0 $ if $n$ is sufficiently large.  Thus  $ v_n ^{ \sharp } \in \Oo$  for all large $n$ and then it is easy to see that 
$ \ds \liminf_{n \ra \infty } E (  v_n ^{ \sharp }  ) \ges \tilde{E}_{min}(m). $ On the other hand, 
$$
\ds \lim_{n \ra \infty }  E(  v_n ^{ \sharp } ) = E ( v_t ) +  \lim_{n \ra \infty } E( v_n - v ) < E(v ) +  \lim_{n \ra \infty } E( v_n - v ) =
 \lim_{n \ra \infty } E( v_n ) =  \tilde{E}_{min}(m), 
$$ 
which is  a contradiction. 

Cases 1-5 here above cover all possible situations and all of them lead to a contradiction. 
We have thus proved that $ \| v \|_{L^2}^2 = m $. 
Now let us prove  that $ v_n \lra v $ in $ H^2( \R^N)$ and  $ v \in \Oo$.  
The weak convergence $ v_n \rightharpoonup v $ in $ L^2 ( \R^N)$ and the convergence of norms 
$ \| v_n \|_{L^2} ^2 \lra m = \| v \|_{L^2} ^2 $ imply that $ v _n \lra v $ strongly in $ L^2 ( \R^N)$. 
Then (\ref{interpol}), (\ref{GNS}) and the boundedness of $ v_n $ in $H^2 ( \R^N)$ imply that 
$ v_n \lra v $ in $ L^{ 2 \si + 2} ( \R^N)$ and $ \nabla v _n \lra \nabla v $ in $ L^2 ( \R^N)$. 
Since $ \Delta v _n \rightharpoonup \Delta v $ in $ L^2 ( \R^N)$ we have 
$ \| \Delta v _n \|_{L^2}^2 =  \| \Delta v  \|_{L^2}^2  + \| \Delta v _n - \Delta v \|_{L^2}^2 + o ( 1 ) $, therefore
$ E( v_n ) = E(v ) +  \| \Delta v _n - \Delta v \|_{L^2}^2 + o ( 1 ) $. 
We have $ E( v_n ) \lra \tilde{E}_{min}(m)$ and we infer that $  \| \Delta v _n - \Delta v \|_{L^2}^2  $ converges in $ \R$. 
It is clear that 
$
D( v_n ) = D(v ) +  \| \Delta v _n - \Delta v \|_{L^2}^2 + o ( 1 ) .
$
Recall that by Lemma \ref{L2.4} we have $ \ds \liminf_{n \ra \infty } D( v_n) > 0$, hence 
$ D(v ) +  \ds \lim_{n \ra \infty } \| \Delta (v _n -v) \|_{L^2}^2 > 0 $. 
We may thus choose $ t \in (0, 1) $ such that $ D(v ) + t^2 \| \Delta( v _n -v ) \|_{L^2}^2 > 0 $ for all $ n$ sufficiently large and we denote 
$ \tilde{v}_n = v + t( v_n - v )$. 
Since $ v_n \lra v $ in $L^2 \cap L^{ 2 \si + 2} ( \R^N)$ we have $ \tilde{v}_n \lra v $ in $L^2 \cap L^{ 2 \si + 2} ( \R^N)$.
Similarly we get $ \nabla \tilde{v}_n \lra \nabla v $ in $L^2( \R^N )$. 
Since $ \Delta v_n - \Delta v \rightharpoonup 0 $ we get 
$ \| \Delta \tilde{v}_n \|_{L^2}^2 = \| \Delta v \|_{L^2}^2 + t^2 \| \Delta ( v_n - v ) \|_{L^2}^2 + o (1) $ and then 
$ D( \tilde{v}_n ) = D(v) + t^2  \| \Delta (v _n -v) \|_{L^2}^2 + o (1)$. 
Therefore $ \| \tilde{v}_n \|_{L^2}^2 < \mu_0 $ and $ D( \tilde{v}_n ) > 0 $ for all sufficiently large $n$, 
which implies that $ \tilde{v}_n \in \Oo $ and consequently $ E( \tilde{v}_n ) \ges \tilde{E}_{min}(\| \tilde{v}_n \|_{L^2}^2 ) $ for all large $n$. 
Letting $ n \lra \infty $ we get 
$$
 \ds \liminf_{n \ra \infty } E( \tilde{v}_n ) \ges \tilde{E}_{min}(m).
$$
On the other hand we have 
$$
E( \tilde{v}_n ) = E( v ) + t^2 \| \Delta( v _n -v ) \|_{L^2}^2  + o (1) = E( v_n )+ ( t^2 -1) \| \Delta( v _n -v ) \|_{L^2}^2  + o (1) 
$$
and letting $ n \lra \infty $ we find
$$
 \tilde{E}_{min}(m) \les  \ds \liminf_{n \ra \infty } E( \tilde{v}_n ) =  \tilde{E}_{min}(m)  +  ( t^2 -1) \lim_{n \ra \infty } \| \Delta( v _n -v ) \|_{L^2}^2.
$$
We conclude that necessarily $ \| \Delta v_n - \Delta v \|_{L^2} \lra 0 $. 
Since $ \|v_n - v \|_{L^2} \lra 0 $, this implies  that $ v_n \lra v $ in $H^2( \R^N)$, as desired. 
Then $ D(v) = \ds \lim_{n \ra \infty } D (v_n)  $ and Lemma \ref{L2.4} implies that $ D(v ) > 0$,  hence $ v \in \Oo$. 
Moreover, we have $ E(v ) =  \ds \lim_{n \ra \infty } E (v_n) =\tilde{E}_{min}(m)$, hence $v$ minimizes $E$ in the set $ \{ w \in \Oo \; | \; \| w \|_{L^2}^2 = m \}$,
 and $ P_1 ( v ) = \ds \lim_{n \ra \infty }  P_1(v_n) =0 $.

Recall that $ v_n = (u_{n })_{t_{u_n, 1}}$ and $ t_{u_n, 1} \lra t_1 \in (0, \infty )$ as $ n \lra \infty$. 
Then we have $ u_n = (v_n)_{t_{u_n, 1}^{-1}} $.
Since $ v_n \lra v $ in $H^2( \R^N)$,  it is easy to show that $ u_n \lra v_{t_1}^{-1}$ in $H^2( \R^N)$.
This implies that $ E( u_n ) \lra E(v_{t_1}^{-1})$, that is $ E(v_{t_1}^{-1}) = \tilde{E}_{min}(m)$.
We have $ D( u_n ) > 0 $ for all $n$ and we infer that $ D(  v_{t_1}^{-1} ) \ges 0 $; in other words, $ {t_1}^{-1} \les t_{v, infl}$. 
Therefore $ 0 < {t_1}^{-1} \les t_{v, infl} $ and $ E( v_{t_1}^{-1} ) = E(v) =  \tilde{E}_{min}(m)$.
Since $ t \longmapsto E( v_t)$ reaches its minimum on $ [0, t_{v, infl}]$ only at $ t = 1$, we infer that necessarily $ t_1 = 1$, thus $ u = v $ and $ u_n \lra u $ strongly in $ H^2( \R^N)$. 
This completes the proof of Theorem \ref{T2.6}.
\hfill
$\Box$

\begin{remark}
\label{nolocmin}
 \rm 

 If there exists $ m_0 > 0 $ such that $ \tilde{E}_{min}(m) = - m $ on $(0, m_0]$, it is easily seen that
 $ \tilde{E}_{min}(m) $ is not achieved for $ m \in (0, m_0)$. 
 Indeed, if $u\in \Oo $ is a minimizer for  $ \tilde{E}_{min}(m) $  then $ \sqrt{a}u \in \Oo $ for 
 $ a > 1 $ and $ a $ close to $1$ and we get   $ \tilde{E}_{min}(am) \les E( \sqrt{a} u ) < a E(u) =  - a m$, contradicting the fact that $ \tilde{E}_{min}(am)= - am $. 
 \end{remark}
 
 \begin{remark}
\label{loc-Lagrange}
 \rm 
Let  $u$ be a minimizer for $ \tilde{E}_{min}(m)$, as given by Theorem \ref{T2.6}.
It is obvious that $ P_1 ( u ) = 0 $ and $u$ satisfies (\ref{ordered}).
In particular, we have $ \| u \|_{H^2} ^2 \les C m = C \| u \|_{L^2}^2$, where $ C$ depends only on $N $ and on $ \si $. 

Since $u$ minimizes $E$ at constant $L^2-$norm in the open set $ \Oo \subset H^2( \R^N)$, 
it is standard to see that there exists a Lagrange multiplier $ \la _ u $ such that (\ref{Eqla}) holds. 
Taking the $H^{-2} - H^2 $ duality product of (\ref{Eqla}) with $u$ we see that $u$  satisfies (\ref{idla}) and this integral identity can be written as
$$
E(u ) - \frac{\si}{\si + 1} \| u \|_{L^{2 \si + 2}}^{ 2 \si + 2 } = \la _u \| u \|_{L^2}^2. 
$$
Since $ 0 < \frac{\si}{\si + 1} \| u \|_{L^{2 \si + 2}}^{ 2 \si + 2 } < \frac{8( N \si - 2)}{N ( N \si - 4 )^2}  \| u \|_{L^2}^2$
(see (\ref{ordered}), 
we infer that 
$$
 -1 \ges  \frac{ \tilde{E}_{min}(m) }{m}  > \la  _u >   \frac{ \tilde{E}_{min}(m) }{m}  -  \frac{8( N \si - 2)}{N ( N \si - 4 )^2} . 
$$ 
Denoting $ \la _u = - 1 - c(u)$ and using the above estimate and Lemma \ref{basic} (i) we see that $ u $ satisfies (\ref{Eq}) with 
$ 0 < c (u) < -1 + \frac{(N \si - 2)^2}{N \si ( N \si - 4)} +  \frac{8( N \si - 2)}{N ( N \si - 4 )^2}$. 
Thus we have an explicit bound on  Lagrange multipliers associated to local minimizers provided by Theorem \ref{T2.6}.

Using (\ref{Planch}) with $ \e = 0 $ and (\ref{idla}) we get $ c(u )   \| u \|_{L^2}^2  <  \| u \|_{L^{2 \si + 2}}^{ 2 \si + 2 } $.
Then using (\ref{GNS}) and (\ref{ordered}) we see that there is $ C > 0 $, depending  only on $N$ and $ \si$, 
such that $ \| u \|_{L^{2 \si + 2}}^{ 2 \si + 2 } \les C \| u \|_{L^2}^{ 2 \si + 2}$. 
These estimates give $ c(u) \les C \| u \|_{L^2}^{ 2 \si } = C m^{\si }$ and we conclude that necessarily $ c( u ) \lra 0 $ as $ m \lra 0 $.

\end{remark}

\begin{remark} 
\label{double}
\rm 
Let $ u_c $ be a minimum action solution of (\ref{Eq}) as provided by Theorem \ref{Tc} and Proposition \ref{groundstates}.
We have already seen (cf. (\ref{converge-m}) and Remark \ref{massasym}) that in the case $ N \si > 4 $ we have $ \| u _c \|_{L^2}^2 \lra 0 $ as $ c \lra \infty$. 
Using (\ref{converge-d}) and (\ref{converge-si}) we see that as $ c \lra \infty, $
$$
( 1 + c )^{ \frac N4 - \frac{1}{\si} -1 } D ( u_c ) \lra \frac{N \si}{4 ( \si + 1)} \left( 2 - \frac{ N \si }{2} \right) I^{ \frac{ \si + 1}{\si }} < 0 .
$$

Therefore for sufficiently large $c$ we have $ \| u _c \|_{L^2}^2 < \mu _0 $ and $ D( u_c ) < 0$, hence $ u_c \not\in \Oo.$
We conclude that if $ c $ is large enough, $u_c $ cannot be a local minimizer of $ E$ when the $L^2-$norm is kept fixed. 
(See  Remark \ref{col} for another interesting variational characterization of $u_c$.)

We have thus proved that there are at least two types of small $L^2-$ norm solutions for equation (\ref{Eq}): 

(a) the minimizers in $ \Oo $ of the energy $E$ at fixed $L^2-$ norm. 
They exist for any $ L^2-$norm smaller than $ \mu _0$. 
Their Lagrange multipliers are bounded, and their $H^2-$norm is controlled by their $ L^2-$norm because  they  satisfy  (\ref{ordered}).

(b) the minimum action solutions $u_c$ for large values of the Lagrange multiplier $c$.
These solutions have large $H^2-$norm:
although $ \| u _c \|_{L^2} \lra 0 $, we have $ \| \Delta u _c \|_{L^2} \lra \infty $ (see Remark \ref{massasym}). 
We were no able to show that the set $ \{ \| u _c \|_{L^2}^2 \; | \; c > 0 \}$ contains an interval of the form $ (0, a)$ with $ a > 0 $.

\end{remark}

If $ N, \si $ and $ m $ satisfy the assumptions of Theorem \ref{Global}, we denote 
$$
\Sr ( m ) = \{ u \in H^2( \R^N) \; \big| \; \| u \|_{L^2 }^2 = m \mbox{ and } E(u) = E_{min}(m) \}. 
$$
If $ N, \si $ and $ m $ satisfy the assumptions of Theorem \ref{T2.6} (thus, in particular, $ N \si > 4$), we denote 
$$
\tilde{\Sr } ( m ) = \{ u \in \Oo \; \big| \; \| u \|_{L^2 }^2 = m \mbox{ and } E(u) =\tilde{E}_{min}(m) \}. 
$$
The sets $ \Sr (m) $ and $ \tilde{\Sr }(m) $ possess some invariances: 

$ \bullet $ For any $ u \in \Sr ( m ) $ and any $ y \in \R^N$, we have $ u ( \cdot + y) \in \Sr ( m ) $.

$ \bullet $ For any $ u \in \Sr ( m ) $ and any orthogonal matrix $ A \in O(N)$, we have $ u ( A\,  \cdot ) \in \Sr ( m ) $.

$ \bullet $ For any $ u \in \Sr ( m ) $ and any $a \in \R$, we have $ e^{i a }u  \in \Sr ( m ) $ and $ \ov{u}  \in \Sr ( m ) $.
\\
Obviously, the same is true with $ \tilde{\Sr } ( m ) $ instead of $ \Sr(m)$. 
Theorems \ref{Global} and \ref{T2.6}, respectively, show that the sets $ \Sr(m)$ and $ \tilde{\Sr } ( m ) $ are compact modulo translations. Using the arguments of T. Cazenave and P.-L. Lions in \cite{CaL}, we get the following (weak) orbital stability result:

\begin{Proposition}
\label{orbital}
Suppose that the assumptions of Theorem \ref{Global} hold. For any $ \e > 0 $ there exists $ \de > 0 $ such that for any 
$ u_0 \in H^2( \R^N)$ satisfying $dist ( u_0, \Sr ( m) ) < \de$, the solution $ t \longmapsto u(t) $ of the 
equation (\ref{BNLSg}) with initial data $ u(0 ) = u_0 $ exists globally in time and satisfies $ dist ( u(t), \Sr(m) ) < \e$ for any $ t \in \R$. Here $ dist ( v, \Sr ( m) ) = \inf \{ \| v - w \|_{H^2} \; \big| \; w \in \Sr ( m) \}$.

The same conclusion holds  under the assumptions of Theorem \ref{T2.6} if $\Sr(m) $ is replaced by $\tilde{\Sr } ( m ) $.

\end{Proposition}

{\it Proof. } We only prove the statement for $\tilde{\Sr } ( m ) $. The proof for $\Sr ( m) $ is similar. 

It follows from Proposition 4.1 p. 204 in \cite{Paus07} that the Cauchy problem for equation (\ref{BNLSg}) is locally well-posed: 
for any $ u_0 \in H^2( \R^N)$ there exists  $ T > 0 $ and a unique solution $ u \in  C([0, T], H^2( \R^N)) $ such that $ u( 0 ) = u_0$. 
If $[0, T^*)$ is the maximal interval of existence of the solution $u(\cdot ) $ then we have
either $ T^* = \infty $ or $ \ds \lim_{t \nearrow T^* } \| u(t ) \|_{H^2} = \infty$, and the $L^2-$norm as well as the energy are conserved:
$ \| u( t ) \|_{L^2 } = \| u _0 \|_{L^2}$ and $ E( u(t) )  = E( u(0)) $ for any $ t \in [0, T^*)$. 

The $ L^2-$norm as well as the functionals $E$ and $D$ are invariant by translations in $ \R^N$.
The set  $\tilde{\Sr} (m)$ is
compact in $H^2( \R^N)$ modulo translations, $ D( u ) > 0 $ and $ \| u \|_{L^2}^2 = m < \mu _0 $   for any $ u \in \tilde{\Sr} (m)$.
It is then easy to see that there exists $ \e_* > 0 $ such that $ D( v ) > 0 $ and $ \| v \|_{L^2}^2 < \mu _0 $  for any $ v \in H^2( \R^N)$ satisfying $ dist ( v,  \tilde{\Sr} (m) ) \les  \e _*. $ In particular, any $v\in H^2( \R^N)$ such that $ dist ( v,  \tilde{\Sr} (m) ) \les  \e _*$
belongs to $ \Oo$. 

We argue by contradiction to prove Proposition \ref{orbital}. 
If the statement is not true, there exist $ \e _0 \in( 0, \e _*) $ and  a sequence $( u_0 ^n )_{n \ges 1 } \subset H^2( \R^N)$ such that $ \de _n := dist( u_0 ^n , \tilde{\Sr } (m)) \lra 0 $ and the solution $ u_n ( \cdot ) $ of  (\ref{BNLSg}) with initial data $ u_n (0) = u_0^n $ either blows up in finite time, or there is some $ t > 0 $ such that $ dist ( u_n (t), \tilde{\Sr } (m) ) \ges \e _0 $. 
In either case, we denote 
$$
t_n = \inf \{ t > 0 \; \big| \; dist ( u_n (t), \tilde{\Sr} (m) ) \ges \e _0 . 
$$
We may assume that   $ \de _n < \e _ 0 $ 
for all $n$, and then we have $ t_n \in (0, \infty)  $ 
and $ dist (u_n ( t_n), \tilde{\Sr } (m)) = \e _0 $ for any $n$. In particular, this implies that $ u_n ( t_n) \in \Oo $ for all $n$. 
Since $ \de _n \lra 0 $, it is standard to prove that $ \| u_0 ^n \|_{L^2}^2 \lra m $ and $ E( u_0 ^n) \lra \tilde{E}_{min}(m)$ as $ n \lra \infty$.
By the conservation of mass and of energy we have 
$$
\| u_{n} ( t_{n }) \|_{L^2}^2 = \| u_0 ^{n} \|_{L^2 } ^2 \lra m 
\quad \mbox{ and } \quad 
E( u_{n} ( t_{n }) )={ E}( u_0 ^{n }  ) \lra \tilde{E}_{min} ( m ) .
$$
Since $ u_{n} ( t_{n }) \in \Oo $, using Theorem \ref{T2.6}  we infer that 
there exist a subsequence $ ( u_{n_k} ( t_{n_k }))_{k \ges 1}$, a sequence $(x_k )_{k \ges 1 } \subset \R^N$ and $ w \in \tilde{\Sr} ( m )$ such that $ \| u_{n_k}( t_{n_k}) ( \cdot + x_k ) - w \|_{H^2} \lra 0 $ as $ k \lra \infty$, 
contradicting the fact that $ dist (u_n ( t_n), \tilde{\Sr } (m)) = \e _0 $ for all $n$.
\hfill
$\Box $

\begin{remark} \rm 
Some related  results have been obtained in \cite{LY}. 
The authors have worked in the space of radial functions 
$H_{rad}^2 = \{ u \in H^2( \R^N) \; | \; u \mbox{ is radially symmetric} \}$ 
and for $m$  sufficiently small  they proved the existence of two solutions of (\ref{Eq}) with $L^2-$norm equal to $m$.
The first one is a local minimizer, and the second one is a mountain-pass type solution. 
The associated Lagrange multipliers are not explicit (they are part of the problem). 

It is an open question whether the minimizers provided by Theorems \ref{Global}, \ref{Tc} and \ref{T2.6} are or not radially symmetric
(some partial results if $\si $ is an integer can be found in \cite{BLSS}).
We could have worked in $H_{rad}^2$, too. 
All our arguments are valid when working in this space, and most proofs become much simpler. 
In this way we get the analogues of Theorems \ref{Global}, \ref{Tc} and \ref{T2.6} in $H_{rad}^2$, which give the existence of radial 
solutions to (\ref{Eq}). 
We do not know whether the energies of  solutions in $H_{rad}^2$ are higher or not than the energies of the corresponding solutions in $H^2( \R^N)$. 
Our main motivation is to understand the existence and the properties of standing waves to a fourth-order non-linear Schr\"odinger equation. 
Since $E$ and the $L^2-$norm are conserved quantities by that equation, the set of travelling waves that we obtain is orbitally stable. 
When working in $H_{rad}^2$ one can get stability only with respect to radial perturbations.

\medskip

\noindent
{\bf Acknowledgement. } The work of MM has been partly funded by the Institut Universitaire de France (IUF Maris 2015).
RM was funded by the Deutsche Forschungsgemeinschaft (DFG, German Research Foundation) - Project-ID 258734477 - SFB 1173.
We are very grateful to the anonymous referee for useful remarks.

\end{remark}

\bigskip

\vspace{1cm}

\noindent Antonio J. \textsc{Fern\'andez}  \smallbreak

Department of Mathematical Sciences, University of Bath, 

Bath BA2 7AY, United Kingdom

{\tt ajf77@bath.ac.uk}

\bigskip

\noindent
Louis \textsc{Jeanjean} \smallbreak

Laboratoire de Math\' ematiques (UMR 6623), Universit\' e Bourgogne Franche-Comt\'e,

16 Route de Gray, 25030 Besan\c con Cedex, France

{\tt louis.jeanjean@univ-fcomte.fr}

\bigskip

\noindent
Rainer  \textsc{Mandel}  \smallbreak

Institute for Analysis, Karlsruhe Institute of Technology,

Englerstra{\ss}e 2 D-76131 Karlsruhe, Germany

{\tt rainer.mandel@kit.edu}

\bigskip

\noindent
Mihai \textsc{Mari\c{s}}  \smallbreak

Institut de Math\'ematiques de Toulouse (UMR 5219), Universit\'e de Toulouse , UPS  \& 

Institut Universitaire de France 

118, Route de Narbonne, 31062 Toulouse, France

{\tt mihai.maris@math.univ-toulouse.fr}

\end{document}